\documentclass[a4paper, 10pt]{amsart}
\usepackage{amscd, amssymb, color, booktabs}
\usepackage{myown01}
\title[Howe Correspondence of Unipotent Characters]
{Howe Correspondence of Unipotent Characters for a Finite Symplectic/Even-orthogonal Dual Pair}
\author{Shu-Yen Pan}
\address{Department of Mathematics,
National Tsing Hua University, Hsinchu 300, Taiwan}
\email{sypan@math.nthu.edu.tw}

\thanks{This work is partially supported by Taiwan MOST-grant 107-2115-M-007-010-MY2.}


\date{\today}

\begin{document}

\begin{abstract}
In this paper we give a complete and explicit description of the Howe correspondence of
unipotent characters for a finite reductive dual pair of a symplectic group and an even
orthogonal group in terms of the Lusztig parametrization when the characteristic of the base
field is not equal to $2$.
That is, the conjecture by Aubert-Michel-Rouquier is confirmed.
\end{abstract}

\maketitle
\tableofcontents


\section{Introduction}

\subsection{}
Let $\omega^\psi_{\Sp_{2N}}$ be the character of the Weil representation (\cf.~\cite{gerardin})
of a finite symplectic group $\Sp_{2N}(q)$ with respect to a nontrivial additive character $\psi$ of a finite field
$\bff_q$ of characteristic $p\neq 2$.
Let $(\bfG,\bfG')$ be one of the following three basic types of reductive dual pairs in $\Sp_{2N}$:
\begin{enumerate}
\item[(1)] two general linear groups $(\GL_n,\GL_{n'})$;

\item[(2)] two unitary groups $(\rmU_n,\rmU_{n'})$;

\item[(3)] one symplectic group and one orthogonal group $(\Sp_{2n},\rmO^\epsilon_{n'})$
\end{enumerate}
where $\epsilon=\pm$.
Now $\omega^\psi_{\Sp_{2N}}$ is regarded as a character of $G\times G'$ and denoted by
$\omega^\psi_{\bfG,\bfG'}$ via the homomorphisms $G\times G'\rightarrow G\cdot G'\hookrightarrow \Sp_{2N}(q)$.
Then $\omega^\psi_{\bfG,\bfG'}$ is decomposed as a sum of irreducible characters
\[
\omega^\psi_{\bfG,\bfG'}=\sum_{\rho\in\cale(G),\ \rho'\in\cale(G')}m_{\rho,\rho'}\rho\otimes\rho'
\]
where each $m_{\rho,\rho'}$ is a non-negative integer,
and $\cale(G)$ denotes the set of irreducible characters of $G$.
Then it establishes a relation
\[
\Theta_{\bfG,\bfG'}=\{\,(\rho,\rho')\in\cale(G)\times\cale(G')\mid m_{\rho,\rho'}\neq 0\,\}
\]
between $\cale(G)$ and $\cale(G')$ which is called the \emph{Howe correspondence}
(or \emph{$\Theta$-correspondence}) for the dual pair $(G,G')$.
The main task is to describe the correspondence explicitly.

\subsection{}
Recall that $\cale(G)$ is partitioned as a disjoint union
\[
\cale(G)=\bigsqcup_{(s)\subset (G^*)^0}\cale(G)_s
\]
of \emph{Lusztig series} $\cale(G)_s$ indexed by the conjugacy classes $(s)$ of
semisimple elements in the connected component $(G^*)^0$ of the dual group $G^*$ of $G$.
Elements in $\cale(G)_1$ are called \emph{unipotent}.
Lusztig shows that there exists a bijection
\[
\grL_s\colon\cale(G)_s\rightarrow\cale(C_{G^*}(s))_1
\]
where $C_{G^*}(s)$ is the centralizer in $G^*$ of $s$.

Now suppose that $\rho\otimes\rho'$ occurs in $\omega^\psi_{\bfG,\bfG'}$
where $\rho\in\cale(G)_s$ and $\rho'\in\cale(G')_{s'}$ for some $s,s'$.
For such a semisimple element $s$ we can define three groups $G^{(1)},G^{(2)},G^{(3)}$ so that
there is a natural bijection $\cale(C_{G^*}(s))_1\simeq\cale(G^{(1)}\times G^{(2)}\times G^{(3)})_1$
(\cf.~\cite{pan-Lusztig-correspondence} subsection 6.2).
Then we have a (modified) Lusztig correspondences
\begin{align*}
\grL_s\colon\cale(G)_s &\rightarrow\cale(G^{(1)}\times G^{(2)}\times G^{(3)})_1, \\
\grL_{s'}\colon\cale(G')_{s'} &\rightarrow\cale(G'^{(1)}\times G'^{(2)}\times G'^{(3)})_1
\end{align*}
such that $G^{(1)}\simeq G^{(2)}$, $G^{(2)}\simeq G^{(2)}$,
and $(G^{(3)},G'^{(3)})$ is a reductive dual pair
of either two general linear groups, two unitary groups,
or one symplectic group and one even orthogonal group.
Therefore we can write $\grL_s(\rho)=\rho^{(1)}\otimes\rho^{(2)}\otimes\rho^{(3)}$ for
$\rho\in\cale(G)_s$ and $\rho^{(j)}\in\cale(G^{(j)})_1$.

Recall that a class function on a classical group $G$ is called \emph{uniform}
if it is a linear combination of the Deligne-Lusztig virtual characters $R_{T,\theta}$.
For a class function $f$ on $G$,
let $f^\sharp$ denote its projection on the subspace of uniform class functions.
If $(G,G')$ consists of one symplectic group and one orthogonal group,
using the decomposition of $\omega_{\bfG,\bfG'}^\sharp$,
it is proved in \cite{pan-Lusztig-correspondence} theorems 6.9 and 7.9
that the following diagram is commutative up to a twist of the sign character of a product of
even orthogonal groups:
\begin{equation}\label{0103}
\begin{CD}
\rho @> \Theta_{\bfG,\bfG'} >> \rho' \\
@V \grL_s VV @VV \grL_{s'} V \\
\rho^{(1)}\otimes\rho^{(2)}\otimes\rho^{(3)} @> {\rm id}\otimes{\rm id}\otimes\Theta_{\bfG^{(3)},\bfG'^{(3)}} >> \rho'^{(1)}\otimes\rho'^{(2)}\otimes\rho^{(3)}
\end{CD}
\end{equation}
i.e., the Howe correspondence and the Lusztig correspondence commute (up to a twist of a sign character).
If the pair $(G,G')$ consists of two general linear groups or two unitary groups,
then all the irreducible characters of $G$ and $G'$ are uniform and so
the above commutative diagram can be read off from the result in \cite{amr} th\'eor\`eme 2.6
(\cf.~\cite{pan-chain01} theorem 3.10).

Therefore we can reduce the Howe correspondence $\Theta_{\bfG,\bfG'}$ of
general irreducible characters to the correspondence $\Theta_{\bfG^{(3)},\bfG'^{(3)}}$
of irreducible unipotent characters.
In \cite{amr} th\'eor\`eme 5.5, th\'eor\`eme 3.10 and conjecture 3.11,
Aubert, Michel and Rouquier give an explicit description (in terms of partitions or bi-partitions)
of the correspondence of unipotent characters for a dual pair of the first two cases
(general linear, or unitary) and provide a conjecture on the third case (symplectic/even orthogonal).
The purpose of this paper is to prove their conjecture.
Unlike the cases of general linear groups or unitary groups,
most of the irreducible characters of symplectic groups or orthogonal groups are not uniform.
This is the main difference and difficulty for studying the correspondence for symplectic/orthogonal dual pairs.
So we need new insight and technique.

\subsection{}
So now we focus on the correspondence of irreducible unipotent characters
for dual pairs of symplectic groups and even orthogonal groups.
First we review some results on the classification of the irreducible unipotent characters
by Lusztig in \cite{lg}, \cite{lg-symplectic} and \cite{lg-orthogonal}.
Let
\[
\Lambda=\binom{A}{B}=\binom{a_1,a_2,\ldots,a_{m_1}}{b_1,b_2,\ldots,b_{m_2}}
\]
denote a (reduced) \emph{symbol} where $A,B$ are finite subsets of non-negative integers such that $0\not\in A\cap B$.
Note that we always assume that $a_1>a_2>\cdots>a_{m_1}$ and $b_1>b_2>\cdots>b_{m_2}$.
The \emph{rank} and \emph{defect} of a symbol $\Lambda$ are defined by
\begin{align*}
\rank(\Lambda)
&= \sum_{a_i\in A}a_i+\sum_{b_j\in B}b_j-\left\lfloor\biggl(\frac{|A|+|B|-1}{2}\biggr)^2\right\rfloor,\\
{\rm def}(\Lambda) &=|A|-|B|
\end{align*}
where $|X|$ denotes the cardinality of a finite set $X$.

Define
\begin{align*}
\cals_{\Sp_{2n}} &=\{\,\Lambda\mid{\rm rank}(\Lambda)=n,\ {\rm def}(\Lambda)\equiv 1\pmod 4\,\};\\
\cals_{\rmO^+_{2n}} &=\{\,\Lambda\mid{\rm rank}(\Lambda)=n,\ {\rm def}(\Lambda)\equiv 0\pmod 4\,\};\\
\cals_{\rmO^-_{2n}} &=\{\,\Lambda\mid{\rm rank}(\Lambda)=n,\ {\rm def}(\Lambda)\equiv 2\pmod 4\,\}.
\end{align*}
Then Lusztig gives a parametrization of irreducible unipotent characters $\cale(G)_1$
by the symbols $\cals_\bfG$.
The irreducible character parametrized by $\Lambda$ will be denoted by $\rho_\Lambda$.

\subsection{}
For a symbol $\Lambda=\binom{a_1,a_2,\ldots,a_{m_1}}{b_1,b_2,\ldots,b_{m_2}}$,
we associate it a \emph{bi-partition}
\[
\Upsilon(\Lambda)
=\sqbinom{a_1-(m_1-1),a_2-(m_1-2),\ldots,a_{m_1-1}-1,a_{m_1}}{b_1-(m_2-1),b_2-(m_2-2),\ldots,b_{m_2-1}-1,b_{m_2}}.
\]
Let $(\bfG,\bfG')=(\Sp_{2n},\rmO^\epsilon_{2n'})$.
For $\Lambda\in\cals_\bfG$, $\Lambda'\in\cals_{\bfG'}$,
we write $\Upsilon(\Lambda)=\sqbinom{\lambda}{\mu}$ and $\Upsilon(\Lambda')=\sqbinom{\lambda'}{\mu'}$.
Then we define a relation on $\cals_\bfG\times\cals_{\bfG'}$:
\begin{align*}
\calb_{\Sp_{2n},\rmO^+_{2n'}}
&=\{\,(\Lambda,\Lambda')\in\cals_{\Sp_{2n}}\times\cals_{\rmO^+_{2n'}}
\mid\mu^\rmT\preccurlyeq\lambda'^\rmT,\ \mu'^\rmT\preccurlyeq\lambda^\rmT,
\ {\rm def}(\Lambda')=-{\rm def}(\Lambda)+1\,\}; \\
\calb_{\Sp_{2n},\rmO^-_{2n'}}
&=\{\,(\Lambda,\Lambda')\in\cals_{\Sp_{2n}}\times\cals_{\rmO^-_{2n'}}
\mid\lambda'^\rmT\preccurlyeq\mu^\rmT,\ \lambda^\rmT\preccurlyeq\mu'^\rmT,\
{\rm def}(\Lambda')=-{\rm def}(\Lambda)-1\,\}
\end{align*}
where $\lambda^\rmT$ denotes the \emph{dual partition} of a partition $\lambda$,
and for two partitions $\lambda=[\lambda_i]$ and $\mu=[\mu_i]$
(with $\{\lambda_i\},\{\mu_i\}$ are written in decreasing order),
we denote $\lambda\preccurlyeq\mu$ if $\mu_i-1\leq\lambda_i\leq\mu_i$ for each $i$.

It is proved in \cite{pan-uniform} that
\[
\omega_{\bfG,\bfG',1}^\sharp
=\sum_{(\Lambda,\Lambda')\in\calb_{\bfG,\bfG'}}\rho_\Lambda^\sharp\otimes\rho_{\Lambda'}^\sharp.
\]
where $\omega_{\bfG,\bfG',1}$ denotes the unipotent part of $\omega_{\bfG,\bfG'}^\psi$.
In this article, we can go a step further to remove the uniform projection
and obtain an explicit description in terms of Lusztig's symbols of the Howe correspondence of
unipotent characters for a symplectic/even-orthogonal dual pair:
\begin{thm*}
Let $(\bfG,\bfG')=(\Sp_{2n},\rmO_{2n'}^\epsilon)$.
Then
\[
\omega_{\bfG,\bfG',1}
=\sum_{(\Lambda,\Lambda')\in\calb_{\bfG,\bfG'}}\rho_\Lambda\otimes\rho_{\Lambda'}.
\]
\end{thm*}
Combining the theorem and the commutativity between Howe correspondence and Lusztig correspondence
in (\ref{0103}),
we obtain a complete description of the whole Howe correspondence of irreducible characters
for a finite reductive dual pair.
Some applications of our main result can be found in \cite{pan-Lusztig-correspondence} and
\cite{pan-eta}.

\subsection{}
The contents of the paper are organized as follows.
In Section 2, we recall the definition and basic properties of symbols introduced by
Lusztig.
Then we discuss the relations $\cald_{Z,Z'}$ and $\calb^\epsilon_{Z,Z'}$ which play the
important roles in our main results.
In Section 3,
we recall the Lusztig's parametrization of irreducible unipotent characters of
a symplectic group or an even orthogonal group.
Then we state our main theorems in Subsection~\ref{0322}.
In Subsection~\ref{0323} we provide a comparison between our theorem and the conjecture
by Aubert-Michel-Rouquier.
In Section 4, we provide several properties of cells of a symplectic group or an
even orthogonal group.
These properties will be used in Section~\ref{1004}.
In last two sections, we prove our main result: Theorem~\ref{0334}.


\section{Symbols and Bi-partitions}
In the first part of this section we recall the notion of ``symbols'' and ``bi-partitions'' from
\cite{lg} \S 3.

\subsection{Symbols}\label{0233}

A \emph{symbol} is an ordered pair
\[
\Lambda=\binom{A}{B}=\binom{a_1,a_2,\ldots,a_{m_1}}{b_1,b_2,\ldots,b_{m_2}}
\]
of two finite subsets $A,B$ (possibly empty) of non-negative integers.
We always assume that elements in $A,B$ are written respectively in strictly decreasing order, i.e.,
$a_1>a_2>\cdots>a_{m_1}$ and $b_1>b_2>\cdots>b_{m_2}$.
A symbol is called \emph{degenerate} if $A=B$, and it is called \emph{non-degenerate} otherwise.
The \emph{size}, \emph{rank} and \emph{defect} of a symbol $\Lambda=\binom{A}{B}$ are defined by
\begin{align}
\begin{split}
{\rm size}(\Lambda)
&=(|A|,|B|) \\
\rank(\Lambda)
&=\sum_{a_i\in A}a_i+\sum_{b_j\in B}b_j-\left\lfloor\biggl(\frac{|A|+|B|-1}{2}\biggr)^2\right\rfloor, \\
{\rm def}(\Lambda)
&=|A|-|B|
\end{split}
\end{align}
where $|X|$ denotes the cardinality of the finite set $X$.
For a symbol $\Lambda$,
let $\Lambda^*$ (resp.~$\Lambda_*$) denote the first row (resp.~second row) of $\Lambda$, i.e.,
$\Lambda=\binom{\Lambda^*}{\Lambda_*}$.
For a symbol $\Lambda=\binom{A}{B}$,
we define its \emph{transpose} $\Lambda^\rmt=\binom{B}{A}$.
A symbol $\binom{A}{B}$ is called \emph{reduced} if $0\not\in A\cap B$.
Let $\cals_{n,\beta}$ denote the set of reduced symbols of rank $n$ and defect $\beta$.
In the remaining part of this article,
a symbol is always assumed to be reduced unless specified otherwise.

\subsection{Special symbols}\label{0232}
A symbol
\begin{equation}\label{0217}
Z=\binom{a_1,a_2,\ldots,a_{m+1}}{b_1,b_2,\ldots,b_m}
\end{equation}
of defect $1$ is called \emph{special} if
\[
a_1\geq b_1\geq a_2\geq b_2\geq\cdots\geq a_m\geq b_m\geq a_{m+1};
\]
similarly a symbol
\begin{equation}\label{0218}
Z=\binom{a_1,a_2,\ldots,a_m}{b_1,b_2,\ldots,b_m}
\end{equation}
of defect $0$ is called \emph{special} if
\[
a_1\geq b_1\geq a_2\geq b_2\geq\cdots\geq b_{m-1}\geq a_m\geq b_m.
\]

For a special symbol $Z=\binom{A}{B}$,
define $Z_\rmI=\binom{A\smallsetminus(A\cap B)}{B\smallsetminus(A\cap B)}$.
Elements in $Z_\rmI$ are called \emph{singles} of $Z$.
The \emph{degree} of a special symbol $Z$ is defined to be
\[
\deg(Z)=\begin{cases}
\frac{|Z_\rmI|-1}{2}, & \text{if $Z$ has defect $1$;}\\
\frac{|Z_\rmI|}{2}, & \text{if $Z$ has defect $0$}.
\end{cases}
\]
A special symbol $Z$ is called \emph{regular} if $Z=Z_\rmI$.
Entries in $Z_{\rm II}:=A\cap B$ are called ``doubles'' in $Z$.
If $a\in A$ and $b\in B$ such that $a=b$,
then the pair $\binom{a}{b}$ is called a \emph{pair of doubles}.

Now for a subset $M$ of $Z_\rmI$,
define
\begin{equation}\label{0201}
\Lambda_M=(Z\smallsetminus M)\cup M^\rmt,
\end{equation}
i.e.,
$\Lambda_M$ is obtained from $Z$ by switching the row position of elements in $M$
and keeping other elements unchanged.
For a special symbol $Z$ and an integer $\delta$,
we define
\begin{align}
\cals_{Z,\delta} &=\{\,\Lambda_M\mid M\subset Z_\rmI,\ {\rm def}(\Lambda_M)=\delta\,\}.
\end{align}
If ${\rm def}(Z)=1$, we define
\begin{equation}\label{0224}
\cals_Z = \bigcup_{\beta\equiv 1\pmod 4}\cals_{Z,\beta};
\end{equation}
if ${\rm def}(Z)=0$, we define
\begin{align}\label{0234}
\cals^+_Z = \bigcup_{\beta\equiv 0\pmod 4}\cals_{Z,\beta},\qquad
\cals^-_Z = \bigcup_{\beta\equiv 2\pmod 4}\cals_{Z,\beta}.
\end{align}
For simplicity, if ${\rm def}(Z)$ is not specified,
let $\cals_Z$ denote one of $\cals_Z$ (when ${\rm def}(Z)=1$) or $\cals^+_Z,\cals^-_Z$ (when ${\rm def}(Z)=0$).
For $\Lambda_M,\Lambda_{M'}\in\cals_Z$,
we define
\[
\Lambda_M+\Lambda_{M'}=\Lambda_N
\]
where $N=(M\cup M')\smallsetminus(M\cap M')$.

\subsection{Bi-partitions}\label{0220}
For a partition
\[
\lambda=[\lambda_i]=[\lambda_1,\lambda_2,\ldots,\lambda_k],\quad\text{ with } \lambda_1\geq\lambda_2\geq\cdots\geq\lambda_k\geq 0,
\]
define $|\lambda|=\lambda_1+\lambda_2+\cdots+\lambda_k$.
For a partition $\lambda=[\lambda_i]$, define its \emph{transpose} (or \emph{dual})
$\lambda^\rmT=[\lambda_j^*]$ by $\lambda_j^*=|\{\,i\mid\lambda_i\geq j\,\}|$ for $j\in\bbN$.

Let $\lambda=[\lambda_1,\ldots,\lambda_k]$ and $\mu=[\mu_1,\ldots,\mu_l]$ be two partitions with
$\lambda_1\geq\cdots\geq\lambda_k$ and $\mu_1\geq\cdots\geq\mu_l$.
We may assume that $k=l$ by adding several $0$'s if necessary.
Then we denote
\[
\lambda\preccurlyeq\mu\quad\text{ if }\mu_i-1\leq\lambda_i\leq\mu_i\text{ for each }i.
\]

Let $\calp_2(n)$ denote the set of bi-partitions of $n$,
i.e., the set of $\sqbinom{\lambda}{\mu}$ where $\lambda,\mu$ are partitions and $|\lambda|+|\mu|=n$.
To each symbol we can associate a bi-partition by:
\begin{equation}\label{0219}
\Upsilon\colon\binom{a_1,a_2,\ldots,a_{m_1}}{b_1,b_2,\ldots,b_{m_2}}\mapsto
\sqbinom{a_1-(m_1-1),a_2-(m_1-2),\ldots,a_{m_1-1}-1,a_{m_1}}{b_1-(m_2-1),b_2-(m_2-2),\ldots,b_{m_2-1}-1,b_{m_2}}.
\end{equation}
It is easy to check that $\Upsilon$ induces a bijection
\[
\Upsilon\colon\cals_{n,\beta}\longrightarrow
\begin{cases}
\calp_2(n-(\frac{\beta-1}{2})(\frac{\beta+1}{2})), & \text{if $\beta$ is odd};\\
\calp_2(n-(\frac{\beta}{2})^2), & \text{if $\beta$ is even}.
\end{cases}
\]

\subsection{The sets $\calb^\epsilon_{Z,Z'}$ and $\cald_{Z,Z'}$}
Let $Z,Z'$ be special symbols of defect $1,0$ respectively.
Define relations $\calb^\epsilon_{Z,Z'}$ on $\cals_Z\times\cals^\epsilon_{Z'}$ by
\begin{align}\label{0223}
\begin{split}
\calb^+_{Z,Z'} &=\{\,(\Lambda,\Lambda')\in\cals_Z\times\cals^+_{Z'}\mid\mu^\rmT\preccurlyeq\xi^\rmT,\ \nu^\rmT\preccurlyeq\lambda^\rmT,
\ {\rm def}(\Lambda')=-{}{\rm def}(\Lambda)+1\,\};\\
\calb^-_{Z,Z'} &=\{\,(\Lambda,\Lambda')\in\cals_Z\times\cals^-_{Z'}\mid\xi^\rmT\preccurlyeq\mu^\rmT,\ \lambda^\rmT\preccurlyeq\nu^\rmT,
\ {\rm def}(\Lambda')=-{\rm def}(\Lambda)-1\,\}
\end{split}
\end{align}
where $\sqbinom{\lambda}{\mu}=\Upsilon(\Lambda)$ and $\sqbinom{\xi}{\nu}=\Upsilon(\Lambda')$.
Then we define
\begin{align}\label{0227}
\begin{split}
\cald_{Z,Z'} &= \calb^+_{Z,Z'}\cap(\cals_{Z,1}\times\cals_{Z',0});\\
\cald_{\Sp_{2n},\rmO^\epsilon_{2n'}}=\cald_{n,n'} &=\bigsqcup_{Z,Z'}\cald_{Z,Z'};\\
\calb_{\Sp_{2n},\rmO^\epsilon_{2n'}}=\calb^\epsilon_{n,n'} &=\bigsqcup_{Z,Z'}\calb^\epsilon_{Z,Z'};\\
\calb^\epsilon &=\bigsqcup_{n,n'\geq 0}\calb^\epsilon_{n,n'}
\end{split}
\end{align}
where the disjoint union $\bigsqcup_{Z,Z'}$ is taken over all special symbols $Z,Z'$ of rank $n,n'$
and defect $1,0$ respectively.
The following two results are from \cite{pan-uniform}:
\begin{lemma}\label{0213}
Let $Z,Z'$ be special symbols of size $(m+1,m), (m',m')$ respectively for some non-negative integers $m,m'$.
If\/ $\cald_{Z,Z'}\neq\emptyset$,
then either $m'=m$ or $m'=m+1$.
\end{lemma}

\begin{lemma}\label{0210}
Let $Z,Z'$ be two special symbols of sizes $(m+1,m),(m',m')$ respectively
where $m'=m,m+1$.
Let
\[
\Lambda=\binom{a_1,a_2,\ldots,a_{m_1}}{b_1,b_2,\ldots,b_{m_2}}\in\cals_Z,\qquad
\Lambda'=\binom{c_1,c_2,\ldots,c_{m'_1}}{d_1,d_2,\ldots,d_{m'_2}}\in\cals_{Z'}.
\]
\begin{enumerate}
\item[(i)] Then $(\Lambda,\Lambda')\in\calb^+_{Z,Z'}$ if and only if
\[
\begin{cases}
m_1'=m_2,\ \text{ and }a_i>d_i,\ d_i\geq a_{i+1},\ c_i\geq b_i,\ b_i>c_{i+1}\text{ \ for each $i$}, & \text{if $m'=m$};\\
m_1'=m_2+1,\ \text{ and }a_i\geq d_i,\ d_i>a_{i+1},\ c_i>b_i,\ b_i\geq c_{i+1}\text{ \ for each $i$}, & \text{if $m'=m+1$}.
\end{cases}
\]

\item[(ii)] Then $(\Lambda,\Lambda')\in\calb^-_{Z,Z'}$ if and only if
\[
\begin{cases}
m_1'=m_2-1,\ \text{ and }d_i\geq a_i,\ a_i>d_{i+1},\ b_i>c_i,\ c_i\geq b_{i+1}\text{ \ for each $i$}, & \text{if $m'=m$};\\
m_1'=m_2,\ \text{ and }d_i>a_i,\ a_i\geq d_{i+1},\ b_i\geq c_i,\ c_i>b_{i+1}\text{ \ for each $i$}, & \text{if $m'=m+1$}.
\end{cases}
\]
\end{enumerate}
\end{lemma}


\section{Finite Howe Correspondence of Unipotent Characters}
In the first part of this section we review the parametrization of
irreducible unipotent characters by Lusztig in \cite{lg-symplectic} and \cite{lg-orthogonal}.
The comparison of our main results and the conjecture in \cite{amr} is in the final subsection.

\subsection{Deligne-Lusztig virtual characters}\label{0339}
The set of irreducible characters of a finite group $G$ is denoted by $\cale(G)$
where an irreducible character means the character of an irreducible representation of $G$.
The space $\calv(G)$ of (complex-valued) class functions on $G$ is an inner product space
and $\cale(G)$ forms an orthonormal basis of $\calv(G)$.

If $\bfG$ is a connected classical group,
let $R_{\bfT,\theta}=R_{\bfT,\theta}^\bfG$ denote the \emph{Deligne-Lusztig (virtual) character} of $G$
with respect to a rational maximal torus $\bfT$ and an irreducible character $\theta\in\cale(T)$ where $T=\bfT^F$.
If $\bfG=\rmO_n^\epsilon$, then we define
\[
R_{\bfT,\theta}^{\rmO_n^\epsilon}=\Ind_{\SO_n^\epsilon(q)}^{\rmO_n^\epsilon(q)}R_{\bfT,\theta}^{\SO_n^\epsilon}.
\]
Let $\calv(G)^\sharp$ denote the subspace of $\calv(G)$ spanned by all Deligne-Lusztig characters of $G$.
For $f\in\calv(G)$, the orthogonal projection $f^\sharp$ of $f$ over $\calv(G)^\sharp$ is called the \emph{uniform projection} of $f$.
A class function $f$ is called \emph{uniform} if $f^\sharp=f$.

If $\bfG$ is connected, it is well-known that the regular character ${\rm Reg}_\bfG$ of $\bfG$ is uniform
(\cf.~\cite{carter-finite} corollary~7.5.6).
Because ${\rm Reg}_{\rmO^\epsilon}=\Ind_{\SO^\epsilon}^{\rmO^\epsilon}({\rm Reg}_{\SO^\epsilon})$,
we see that ${\rm Reg}_{\rmO^\epsilon}$ is also uniform.
Therefore,
we have
\begin{equation}\label{0330}
\rho(1)
=\langle\rho,{\rm Reg}_\bfG\rangle
=\langle\rho^\sharp,{\rm Reg}_\bfG\rangle
=\rho^\sharp(1).
\end{equation}
This means that $\rho^\sharp\neq0$ for any $\rho\in\cale(G)$.

An irreducible character $\rho$ of $G$ is called \emph{unipotent} if
$\langle \rho, R_{\bfT,1}^\bfG\rangle_\bfG\neq 0$ for some $\bfT$.
The set of irreducible unipotent characters of $G$ is denoted by $\cale(G)_1$.

\subsection{Unipotent characters of $\Sp_{2n}(q)$}\label{0335}
Define
\begin{equation}\label{}
\cals_{\Sp_{2n}} = \bigsqcup_{\beta\equiv 1\pmod 4}\cals_{n,\beta}.
\end{equation}
From \cite{lg} theorem 8.2,
there is a parametrization of unipotent characters $\cals_{\Sp_{2n}}\rightarrow\cale(\Sp_{2n})_1$
denoted by $\Lambda\mapsto\rho_\Lambda$.
Every symbol $\Lambda$ of rank $n$ and defect $\beta\equiv 1\pmod 4$ is in $\cals_Z$ for some unique
special symbol $Z$ of rank $n$ and defect $1$, i.e.,
$\cals_{\Sp_{2n}}=\bigsqcup_Z\cals_Z$ where $Z$ runs over all special symbols of rank $n$ and defect $1$.

It is know that there is a one-to-one correspondence between the set $\calp_2(n)$ of bipartitions of $n$ and
the set of irreducible characters of the Weyl group $W_n$ of $\Sp_{2n}$ (\cf.~\cite{Geck-Pfeiffer} theorem~5.5.6).
Then for a symbol $\Sigma\in\cals_{n,1}$,
we can associate it a uniform function $R_\Sigma$ on $\Sp_{2n}(q)$ given by $R_\Sigma=R_\chi$
where $R_\chi=R_\chi^\bfG$ is defined in \cite{pan-uniform} subsection~3.2, 
and $\chi$ is the irreducible character of $W_n$
associated to $\Upsilon(\Sigma)$ where $\Upsilon$ is the bijection $\cals_{n,1}\rightarrow\calp_2(n)$
given in (\ref{0219}).

Let $\calv_Z$ denote the (complex) vector space spanned by $\{\,\rho_\Lambda\mid\Lambda\in\cals_Z\,\}$.
It is known that $\{\,\rho_\Lambda\mid\Lambda\in\cals_Z\,\}$ is an orthonormal basis for $\calv_Z$,
and $\{\,R_\Sigma\mid\Sigma\in\cals_{Z,1}\,\}$ is an orthonormal basis for the uniform projection $\calv_Z^\sharp$ of the space $\calv_Z$.

The following proposition is modified from \cite{lg-symplectic} theorem 5.8:

\begin{proposition}[Lusztig]\label{0306}
Let $\bfG=\Sp_{2n}$, $Z$ a special symbol of rank $n$, defect $1$ and degree $\delta$.
For $\Sigma\in\cals_{Z,1}$,
we have
\[
\langle R_\Sigma,\rho_{\Lambda}\rangle_\bfG
=\begin{cases}
(-1)^{\langle\Sigma,\Lambda\rangle}2^{-\delta}, & \text{if $\Lambda\in\cals_Z$};\\
0, & \text{otherwise}
\end{cases}
\]
where
$\langle,\rangle\colon\cals_{Z,1}\times\cals_Z\rightarrow\bff_2$ is given by
$\langle\Lambda_N,\Lambda_M\rangle=|N\cap M|\pmod 2$.
\end{proposition}

It is known that $\Sp_{2n}(q)$ has a unique irreducible unipotent cuspidal character if and only if $n=m(m+1)$ for some
nonnegative integer $m$.
In this case, the unique unipotent cuspidal character $\rho_\Lambda$ is associated to the
symbol
\begin{equation}\label{0320}
\Lambda=\begin{cases}
\binom{2m,2m-1,\ldots,1,0}{-}, & \text{if $m$ is even};\\
\binom{-}{2m,2m-1,\ldots,1,0}, & \text{if $m$ is odd}
\end{cases}
\end{equation}
of defect $(-1)^m(2m+1)$ under our convention.

\subsection{Unipotent characters of $\rmO_{2n}^\epsilon(q)$}\label{0331}
Define
\begin{align}\label{}
\begin{split}
\cals_{\rmO_{2n}^+} &= \bigsqcup_{\beta\equiv 0\pmod 4}\cals_{n,\beta}=\bigsqcup_{Z\ \text{special},\ {\rm def}(Z)=0,\ {\rm rank}(Z)=n}\cals^+_Z;\\
\cals_{\rmO_{2n}^-} &= \bigsqcup_{\beta\equiv 2\pmod 4}\cals_{n,\beta}=\bigsqcup_{Z\ \text{special},\ {\rm def}(Z)=0,\ {\rm rank}(Z)=n}\cals^-_Z.
\end{split}
\end{align}
From \cite{lg} theorem 8.2,
we know that there exists a one-to-one correspondence between $\cals_{\rmO^\epsilon_{2n}}$
and the unipotent characters $\cale(\rmO^\epsilon_{2n}(q))_1$.
For $\Lambda\in\cals_{\rmO_{2n}^\epsilon}$,
the associated irreducible character of $\rmO_{2n}^\epsilon(q)$ is denoted by $\rho_\Lambda$.
It is known that $\rho_{\Lambda^\rmt}=\rho_\Lambda\cdot\sgn$.

Let $Z$ be a special symbol of defect $0$.
As in the symplectic case,
let $\calv^\epsilon_Z$ denote the (complex) vector space spanned by $\{\,\rho_\Lambda\mid\Lambda\in\cals^\epsilon_Z\,\}$.
Then $\{\,\rho_\Lambda\mid\Lambda\in\cals^\epsilon_Z\,\}$ is an orthonormal basis for $\calv^\epsilon_Z$.
If $\Sigma$ is degenerate,
then $Z=\Sigma$, $\cals^-_Z=\emptyset$, $\cals^+_Z=\{Z\}$ and $\rho_\Sigma^{\rmO^+}=R_\Sigma^{\rmO^+}$.
If $\Sigma$ is non-degenerate,
let $\bar\cals_{Z,0}$ denote a complete set of representatives of cosets $\{\Sigma,\Sigma^\rmt\}$ in $\cals_{Z,0}$,
then $\{\,\frac{1}{\sqrt 2}R^{\rmO^\epsilon}_\Sigma\mid\Sigma\in\bar\cals_{Z,0}\,\}$ is an orthonormal basis for $(\calv_Z^\epsilon)^\sharp$.

The following proposition is a modification for $\rmO^\epsilon_{2n}(q)$
from \cite{lg-orthogonal} theorem 3.15:

\begin{proposition}\label{0319}
Let $Z$ be a non-degenerate special symbol of defect $0$ and degree $\delta\geq 1$.
For any $\Sigma\in\cals_{Z,0}$,
we have
\[
\langle R^{\rmO^\epsilon}_\Sigma,\rho^{\rmO^\epsilon}_\Lambda\rangle_{\rmO^\epsilon}
=\begin{cases}
(-1)^{\langle\Sigma,\Lambda\rangle}2^{-(\delta-1)}, & \text{if $\Lambda\in\cals^\epsilon_Z$};\\
0, & \text{otherwise}
\end{cases}
\]
where $\langle,\rangle\colon\cals_{Z,0}\times\cals^\epsilon_Z\rightarrow\bff_2$ by
$\langle\Lambda_M,\Lambda_N\rangle=|M\cap N|\pmod 2$.
\end{proposition}

The trivial character of the trivial group $\rmO^+_0(q)$ is regarded as unipotent cuspidal
and is associated to the symbol $\binom{-}{-}$ of rank $0$ and defect $0$.
It is known that $\rmO_{2n}^\epsilon(q)$ for $n$ positive has two irreducible unipotent cuspidal characters
if and only if $n=m^2$ and $\epsilon=(-1)^m$ for some nonzero integer $m$.
These two characters $\zeta^{\rm I},\zeta^{\rm II}$ will associated to the symbols $\Lambda,\Lambda^\rmt$ respectively
where
\begin{equation}\label{0332}
\Lambda=\begin{cases}
\binom{2m-1,2m-2,\ldots,1,0}{-}, & \text{if $m$ is even}; \\
\binom{-}{2m-1,2m-2,\ldots,1,0}, & \text{if $m$ is odd}
\end{cases}
\end{equation}
under our convention.
More properties of the parametrization will be given in Subsection~\ref{1006}.

\subsection{The main results}\label{0322}
Let $(\bfG,\bfG')=(\Sp_{2n},\rmO^\epsilon_{2n'})$.
It is known that unipotent characters are preserved in the Howe correspondence for the dual pair $(\bfG,\bfG')$.
Let $\omega_{\bfG,\bfG',1}$ denote the unipotent part of the Weil character $\omega^\psi_{\bfG,\bfG'}$.
The following result is proved in \cite{pan-uniform}:

\begin{proposition}\label{0310}
Let $(\bfG,\bfG')=(\Sp_{2n},\rmO^\epsilon_{2n'})$.
Then
\[
\omega_{\bfG,\bfG',1}^\sharp
=\sum_{(\Lambda,\Lambda')\in\calb_{\bfG,\bfG'}}\rho_\Lambda^\sharp\otimes\rho_{\Lambda'}^\sharp
\]
\end{proposition}

All the efforts in this article are to remove the uniform projection of both sides
of the above identity:

\begin{theorem}\label{0334}
Let $(\bfG,\bfG')=(\Sp_{2n},\rmO_{2n'}^\epsilon)$.
Then
\[
\omega_{\bfG,\bfG',1}
=\sum_{(\Lambda,\Lambda')\in\calb_{\bfG,\bfG'}}\rho_\Lambda\otimes\rho_{\Lambda'}.
\]
\end{theorem}

The proof of our main theorem is divided into two stages (Section~\ref{1004} and Section~\ref{0618}):
\begin{itemize}
\item

To recover the relation between $\omega_{\bfG,\bfG',1}$ and
$\sum_{(\Lambda,\Lambda')\in\calb_{\bfG,\bfG'}}\rho_\Lambda\otimes\rho_{\Lambda'}$
from the uniform projection,
we will use the technique learned from \cite{Kable-Sanat}, pp.436--438.
That is, we reduce the problem into a system of linear equations.
To write down these equations,
we need the theory of ``\emph{cells}'' by Lusztig from \cite{lg-symplectic} theorem 5.6 and \cite{lg-orthogonal} proposition 3.13.
The variables of the equations are the multiplicities of those $\rho_\Lambda\otimes\rho_{\Lambda'}$
occurring in $\omega_{\bfG,\bfG',1}$.
The solutions must be non-negative integers,
that is the reason that we are almost able to solve the equations.
Due to the disconnectedness of $\rmO_{2n}^\epsilon$,
irreducible characters $\rho_{\Lambda'},\rho_{\Lambda'^\rmt}$ are not distinguishable by
Deligne-Lusztig virtual characters.
So in the first stage we can only conclude that
$\rho_\Lambda\otimes\rho_{\Lambda'}$ or $\rho_\Lambda\otimes\rho_{\Lambda'^\rmt}$
occur in $\omega_{\bfG,\bfG',1}$ if and only if
$(\Lambda,\Lambda')$ or $(\Lambda,\Lambda'^\rmt)$ occur in $\calb_{\bfG,\bfG'}$.

\item

Because the Howe correspondence and the parametrization $\Lambda\rightarrow\rho_\Lambda$ are
both compatible with parabolic induction,
the ambiguity in the first stage can be removed once the correspondence of unipotent
cuspidal characters is fixed.
The proof of the theorem is in Subsection~\ref{1006} for ${\rm def}(\Lambda')>0$, and
in Subsection~\ref{1118} for ${\rm def}(\Lambda')=0$.
\end{itemize}

\subsection{The conjecture by Aubert-Michel-Rouquier}\label{0323}
In Theorem~\ref{0334},
we describe the Howe correspondence of unipotent characters in terms of ``symbols'';
the conjecture in \cite{amr} p.383 describes the correspondence in terms of ``bi-partitions''.
The main difference between these two descriptions is that a bi-partition does not contain
the information of the ``defect'' of a symbol which is controlled by the unipotent cuspidal characters.
Therefore, the description in \cite{amr} p.383 needs to specify the correspondence of unipotent
cuspidal characters first.
Now we want to make the comparison more explicit.

In our convention,
we always assume that the defect of a symbol for a symplectic group (resp.~split even orthogonal group, non-split even orthogonal group)
is $1\pmod 4$ (resp.~$0\pmod 4$, $2\pmod 4$).
Our convention is different from the original one in \cite{lg} p.134 where the defect of
a symbol is always assumed to be non-negative.
We can see that the relation between our bi-partition $\sqbinom{\phi}{\psi}$ and the
bi-partition $\phi\boxtimes\psi$ in \cite{amr} p.379 is given by
\[
\sqbinom{\phi}{\psi}\mapsto
\begin{cases}
\phi\boxtimes\psi, & \text{if ${\rm def}(\Lambda)>0$};\\
\psi\boxtimes\phi, & \text{if ${\rm def}(\Lambda)\leq 0$}
\end{cases}
\]
where $\sqbinom{\phi}{\psi}=\Upsilon(\Lambda)$.

Let $\Lambda,\Lambda'$ be two symbols and write $\sqbinom{\phi}{\psi}=\Upsilon(\Lambda)$ and
$\sqbinom{\phi'}{\psi'}=\Upsilon(\Lambda')$,
and let $x_i,x_i^*$ be the notation used in \cite{amr} p.383.
Let $\lambda_k$ be the unique unipotent cuspidal character of $\Sp_{2k(k+1)}(q)$,
and let $\lambda_k^{\rm I},\lambda_k^{\rm II}$ be the two unipotent cuspidal characters of $\rmO_{2k^2}^{\epsilon_k}(q)$ ($k>0$)
where $\epsilon_k=\sgn((-1)^k)$ such that $\lambda_k\otimes\lambda_k^{\rm II}$ and $\lambda_k\otimes\lambda_{k+1}^{\rm I}$
occur in the Howe correspondence.
Let $\Theta_{\lambda_k,\lambda_k^{\rm II}}$ and $\Theta_{\lambda_k,\lambda_{k+1}^{\rm I}}$ be the mappings between
bi-partitions defined in \cite{amr} p.383.

First we consider the correspondence $\Theta_{\lambda_k,\lambda_k^{\rm II}}$.
\begin{enumerate}
\item[(1)] Suppose that $k$ is even.
Then $\epsilon=\sgn((-1)^k)=+$.
Now ${\rm def}(\Lambda)>0$ from (\ref{0320}) and ${\rm def}(\Lambda')\leq 0$ from (\ref{0223}).
Then the condition: $\psi'^\rmT\preccurlyeq\phi^\rmT$ and $\psi^\rmT\preccurlyeq\phi'^\rmT$
is equivalent to the condition: $\psi'\boxtimes\phi'=x_i^*(\phi)\boxtimes x_j(\psi)$ for some $i,j\geq 0$.

\item[(2)] Suppose that $k$ is odd.
Then $\epsilon=\sgn((-1)^k)=-$.
Now ${\rm def}(\Lambda)<0$ and ${\rm def}(\Lambda')>0$.
Then the condition:
$\phi^\rmT\preccurlyeq\psi'^\rmT$ and $\phi'^\rmT\preccurlyeq\psi^\rmT$ is equivalent to the condition:
$\phi'\boxtimes\psi'=x_i^*(\psi)\boxtimes x_j(\phi)$ for some $i,j\geq 0$.
\end{enumerate}

Next we consider the correspondence $\Theta_{\lambda_k,\lambda_{k+1}^{\rm I}}$.
\begin{enumerate}
\item[(3)] Suppose that $k$ is even.
Then $\epsilon=\sgn((-1)^{k+1})=-$.
Now ${\rm def}(\Lambda)>0$ and ${\rm def}(\Lambda')\leq 0$.
Then the condition:
$\phi^\rmT\preccurlyeq\psi'^\rmT$ and $\phi'^\rmT\preccurlyeq\psi^\rmT$ is equivalent to the condition:
$\psi'\boxtimes\phi'=x_i(\phi)\boxtimes x_j^*(\psi)$ for some $i,j\geq 0$.

\item[(4)] Suppose that $k$ is odd.
Then $\epsilon=\sgn((-1)^{k+1})=+$.
Now ${\rm def}(\Lambda)<0$ and ${\rm def}(\Lambda')>0$.
Then the condition: $\psi'^\rmT\preccurlyeq\phi^\rmT$ and $\psi^\rmT\preccurlyeq\phi'^\rmT$
is equivalent to the condition: $\phi'\boxtimes\psi'=x_i(\psi)\boxtimes x_j^*(\phi)$ for some $i,j\geq 0$.
\end{enumerate}

Therefore we see that Theorem~\ref{0334} is equivalent to conjecture 3.11 in \cite{amr}:

\begin{corollary}
Keep the notations in \cite{amr} p.383.
Then
\begin{itemize}
\item[(i)] the mapping $\Theta_{\lambda_k,\lambda_{k+1}^\rmI}$ is given by $\phi\boxtimes\psi\mapsto X\phi\boxtimes X^*(\psi)$; and

\item[(ii)] the mapping $\Theta_{\lambda_k,\lambda_k^{\rm II}}$ is given by $\phi\boxtimes\psi\mapsto X^*(\phi)\boxtimes X\psi$.
\end{itemize}
\end{corollary}


\section{Cells for Symplectic Groups and Even Orthogonal Groups}
In this section, we provide several technical lemmas which are needed in the next
two sections.

\subsection{Consecutive pairs}
Let $Z$ be a special symbol of degree $\delta$.
A set $\Psi$ of some pairs (possibly empty) in $\Phi$ is called a \emph{subset of pairs} of $\Phi$
and is denoted by $\Psi\leq\Phi$.
A pair $\binom{s}{t}$ in $Z_\rmI$ is called \emph{consecutive} if there is no other element in $Z$
lying between $s$ and $t$ i.e., there is no element $x$ in $Z$ such that $s<x<t$ or $t<x<s$.

For a set of (disjoint) consecutive pairs $\Psi_0$ in $Z_\rmI$,
we define several subsets of $\cals_Z$ as follows.
\begin{enumerate}
\item[(1)] Suppose $Z$ is of defect $1$ and degree $\delta$.
Define
\begin{align*}
\cals_{Z,\Psi_0} &=\{\,\Lambda_M\in\cals_Z\mid M\leq\Psi_0\,\}\subset\cals_{Z,1};\\
\cals_Z^{\Psi_0} &=\{\,\Lambda_M\in\cals_Z\mid M\subset Z_\rmI\smallsetminus\Psi_0\,\}.
\end{align*}
Clearly, $\cals_{Z,\Psi_0}=\{Z\}$ and $\cals_Z^{\Psi_0}=\cals_Z$ if $\Psi_0=\emptyset$.
If $\delta_0$ is the number of pairs in $\Psi_0$,
then it is clear that $|\cals_{Z,\Psi_0}|=2^{\delta_0}$ and $|\cals_Z^{\Psi_0}|=2^{2(\delta-\delta_0)}$.

\item[(2)] Suppose $Z$ is of defect $0$.
Define
\begin{align*}
\cals_{Z,\Psi_0} &=\{\,\Lambda_M\in\cals^+_Z\mid M\leq\Psi_0\,\}\subset\cals_{Z,0};\\
\cals_Z^{\epsilon,\Psi_0} &=\{\,\Lambda_M\in\cals^\epsilon_Z\mid M\subset Z_\rmI\smallsetminus\Psi_0\,\}.
\end{align*}
\end{enumerate}
Note that if ${\rm def}(Z)$ is not specified,
we will just use the notation $\cals_Z^{\Psi_0}$ to denote $\cals_Z^{\Psi_0}$ (when ${\rm def}(Z)=1$) or
$\cals_Z^{+,\Psi_0},\cals_Z^{-,\Psi_0}$ (when ${\rm def}(Z)=0$).

\begin{lemma}\label{0521}
Let $\Psi_0$ be a set of consecutive pairs in $Z_\rmI$.
Suppose that $\Lambda_1,\Lambda_1'\in\cals_Z^{\Psi_0}$ and $\Lambda_2,\Lambda_2'\in\cals_{Z,\Psi_0}$.
If $\rho_{\Lambda_1+\Lambda_2}\simeq\rho_{\Lambda_1'+\Lambda_2'}$,
then $\Lambda_1=\Lambda_1'$ and $\Lambda_2=\Lambda'_2$.
\end{lemma}
\begin{proof}
If $\rho_{\Lambda_1+\Lambda_2}\simeq\rho_{\Lambda_1'+\Lambda_2'}$,
then $\Lambda_1+\Lambda_2=\Lambda'_1+\Lambda'_2$ and hence
$\Lambda_1+\Lambda'_1=\Lambda_2+\Lambda'_2$.
Note that both $\cals_{Z,\Psi_0}$ and $\cals^{\Psi_0}_Z$ are closed under addition and
$\cals_{Z,\Psi_0}\cap\cals_Z^{\Psi_0}=\{Z\}$.
Therefore $\Lambda_1+\Lambda_1'=Z$, i.e., $\Lambda_1=\Lambda'_1+Z=\Lambda'_1$.
Similarly, $\Lambda_2=\Lambda'_2$.
\qed
\end{proof}

\subsection{Cells}\label{0609}
We first recall the notion of ``a cell'' by Lusztig from \cite{lg-symplectic} and \cite{lg-orthogonal}.
Let $Z$ be a special symbol of degree $\delta$.

\begin{enumerate}
\item[(1)] If $Z$ is of defect $1$, then an \emph{arrangement} of $Z_\rmI$ is a partition $\Phi$
of the $2\delta+1$ singles in $Z_\rmI$ into $\delta$ (disjoint) pairs and
one isolated element such that each pair contains one element in the first row and one element in the
second row of $Z_\rmI$.

\item[(2)] If $Z$ is of defect $0$, then an \emph{arrangement} of $Z_\rmI$ is a partition $\Phi$
of the $2\delta$ singles in $Z_\rmI$ into $\delta$ pairs such that each pair contains one element
in the first row and one element in the second row of $Z_\rmI$.
\end{enumerate}

For a subset of pairs $\Psi$ of an arrangement $\Phi$ of $Z_\rmI$,
as usual let $\Psi^*$ (resp.\ $\Psi_*$) denote the set of entries in the first
(resp.\ second) row in $\Psi$.
Recall that the following uniform class function on $G$ is defined in \cite{lg-symplectic}:
\begin{equation}\label{0514}
R_{\underline c}=R_{\underline c(Z,\Phi,\Psi)}=\sum_{\Psi'\leq\Phi}(-1)^{|(\Phi\smallsetminus\Psi)\cap\Psi'^*|}
R_{\Lambda_{\Psi'}}
\end{equation}
where
$\Lambda_{\Psi'}=(Z\smallsetminus\Psi')\cup\Psi'^\rmt$
is defined as in (\ref{0201}).

\begin{remark}
Our notation is slightly different from that in \cite{lg-symplectic} and \cite{lg-orthogonal}.
More precisely, the uniform class function $R_{\underline c(Z,\Phi,\Psi)}$ in (\ref{0514}) is denoted by
$R(\underline c(Z,\Phi,\Phi\smallsetminus\Psi))$ in \cite{lg-symplectic} and \cite{lg-orthogonal}.
\end{remark}

For a subset of pairs $\Psi$ of an arrangement $\Phi$,
we define
\begin{equation}\label{0302}
C_{\Phi,\Psi}=\{\,\Lambda_M\in\cals_Z
\mid |M\cap\Psi'|\equiv|(\Phi\smallsetminus\Psi)\cap\Psi'^*|\pmod 2\text{ for all }\Psi'\leq\Phi\,\}.
\end{equation}
Such a subset $C_{\Phi,\Psi}$ of $\cals_Z$ will be also called a \emph{cell}.
If we need to emphasize the special symbol $Z$,
the cell $C_{\Phi,\Psi}$ is also denoted by $C_{Z,\Phi,\Psi}$.
From the definition it is easy to see that a symbol $\Lambda_M\in\cals_Z$ is in $C_{\Phi,\Psi}$ if and only if
$M$ satisfies the following two conditions:
\begin{itemize}
\item $M$ contains either none or two elements of each pair in $\Psi$; and

\item $M$ contains exactly one element of each pair in $\Phi\smallsetminus\Psi$.
\end{itemize}
In particular,
it is clear from the definition that
\[
C_{\Phi,\Phi}=\{\,\Lambda_M\in\cals_Z\mid M\leq\Phi\,\}.
\]

\begin{remark}\label{0608}
Suppose that $Z$ is of defect $1$ and $\Lambda_M$ is in $C_{\Phi,\Psi}$ for some $\Phi,\Psi$.
The requirement $|M^*|\equiv |M_*|\pmod 2$ implies that $M$ must contain the isolated element in the arrangement $\Phi$
if $\Psi$ consists of odd number of pairs; and $M$ does not contain the isolated element if $\Psi$ consists of even number of pairs.
\end{remark}

\begin{example}
Suppose that $Z=\binom{4,2,0}{3,1}$, $\Phi=\{\binom{4}{-},\binom{2}{3},\binom{0}{1}\}$
and $\Psi=\{\binom{0}{1}\}$.
There are four possible $M$ that satisfies condition in (\ref{0302}),
namely, $\binom{4,2,0}{3}$, $\binom{4}{1}$, $\binom{4,0}{-}$ and $\binom{4,2}{3,1}$
(note that we require that $|M^*|\equiv|M_*|\pmod 2$ to satisfy ${\rm def}(\Lambda_M)\equiv 1\pmod 4$),
and resulting $\Lambda_M$ are
$\binom{3}{4,2,1,0}$, $\binom{2,1,0}{4,3}$, $\binom{2}{4,3,1,0}$ and $\binom{3,1,0}{4,2}$ respectively,
i.e.,
\[
\textstyle
C_{\Phi,\Psi}=\{\binom{3}{4,2,1,0}, \binom{2,1,0}{4,3}, \binom{2}{4,3,1,0}, \binom{3,1,0}{4,2}\}.
\]
\end{example}

\begin{lemma}\label{0304}
Suppose that $\Lambda_{M_1},\Lambda_{M_2}$ are both in $C_{\Phi,\Psi}$ for some
arrangement $\Phi$ of $Z_\rmI$ and $\Psi\leq\Phi$.
Then
\[
|M_1\cap\Psi'|\equiv|M_2\cap\Psi'|\pmod 2
\]
for any $\Psi'\leq\Phi$.
\end{lemma}
\begin{proof}
Suppose that $\Lambda_{M_1},\Lambda_{M_2}\in C_{\Phi,\Psi}$.
Then by (\ref{0302}) we have
\[
|M_1\cap\Psi'|\equiv|(\Phi\smallsetminus\Psi)\cap\Psi'^*|\equiv|M_2\cap\Psi'|\pmod 2
\]
for all subsets $\Psi'$ of pairs of $\Phi$.
\qed
\end{proof}

\begin{lemma}\label{0515}
Suppose that $\Phi$ is an arrangement of $Z_\rmI$,
$\Psi$ is a subset of pairs and $\Psi_0$ is a subset of consecutive pairs such that
$\Psi_0\leq\Psi\leq\Phi$.
Then
\[
C_{\Phi,\Psi}=(C_{\Phi,\Psi}\cap\cals_Z^{\Psi_0})+\cals_{Z,\Psi_0}
=\{\,\Lambda_1+\Lambda_2\mid\Lambda_1\in(C_{\Phi,\Psi}\cap\cals_Z^{\Psi_0}),\ \Lambda_2\in\cals_{Z,\Psi_0}\,\}.
\]
\end{lemma}
\begin{proof}
Suppose $Z$ is of defect $1$.
Let $\Lambda_M$ be an element in $C_{\Phi,\Psi}$ for some $M\subset Z_\rmI$.
Then $M=M_1\sqcup M_2$ where $M_1=(M\cap(Z_\rmI\smallsetminus\Psi_0))$ and $M_2=(M\cap\Psi_0)$.
And we have $\Lambda_M=\Lambda_{M_1}+\Lambda_{M_2}$.
From the requirement of the set $M$ before Remark~\ref{0608},
we know that $M_2$ is a subset of pairs in $\Psi_0$, i.e., $M_2\leq\Psi_0$ and hence $\Lambda_{M_2}\in\cals_{Z,\Psi_0}$;
and clearly $\Lambda_{M_1}\in C_{\Phi,\Psi}\cap\cals_{Z,\Psi_0}$.
On the other hand, if $\Lambda_{M_1}\in C_{\Phi,\Psi}\cap\cals_Z^{\Psi_0}$ and $\Lambda_{M_2}\in\cals_{Z,\Psi_0}$,
then it is obvious that $\Lambda_{M_1}+\Lambda_{M_2}=\Lambda_{M_1\cup M_2}\in C_{\Phi,\Psi}$.

The proof for the case that $Z$ is of defect $0$ is similar.
\qed
\end{proof}

\begin{lemma}\label{0517}
Let $\Phi$ be an arrangement of $Z_\rmI$, $\Psi_0$ a subset of consecutive pairs,
and $\Lambda\in\cals_Z^{\Psi_0}$.
Suppose that $\Lambda\in C_{\Phi,\Psi}$ for some $\Psi\leq\Phi$.
Then $\Psi_0\leq\Psi$.
\end{lemma}
\begin{proof}
Suppose that $\Lambda=\Lambda_M$ for some $M\subset Z_\rmI\smallsetminus\Psi_0$, i.e., $M\cap\Psi_0=\emptyset$.
From the rule before Remark~\ref{0608},
the assumption $\Lambda\in C_{\Phi,\Psi}$ implies that $M$ contains exactly one element from each pair in $\Phi\smallsetminus\Psi$.
Therefore we must have $\Psi_0\leq\Psi$.
\qed
\end{proof}

\subsection{Cells for symplectic groups}\label{0604}
In this subsection,
let $\bfG=\Sp_{2n}$, and $Z$ a special symbol of rank $n$ and defect $1$.

\begin{lemma}\label{0303}
Let $Z$ be a special symbol of defect $1$ and degree $\delta$, $\Phi$ a fixed arrangement of $Z_\rmI$,
$\Psi,\Psi'$ subsets of pairs of\/ $\Phi$.
Then
\begin{enumerate}
\item[(i)] $|C_{\Phi,\Psi}|=2^\delta$;

\item[(ii)] $C_{\Phi,\Psi}\cap C_{\Phi,\Psi'}=\emptyset$ if\/ $\Psi\neq\Psi'$;

\item[(iii)] $\cals_Z=\bigcup_{\Psi\leq\Phi}C_{\Phi,\Psi}$.
\end{enumerate}
\end{lemma}
\begin{proof}
Let $z_0$ denote the isolated element in $\Phi$.
Suppose that $\Lambda_M$ is an element of $C_{\Phi,\Psi}$.
From the conditions before Remark~\ref{0608},
we can write $M=M_1\sqcup M_2$ where $M_1$ consists of exactly one element from each pair
of $\Phi\smallsetminus\Psi$ and possibly $z_0$ such that $|M_1|$ is even,
and $M_2$ consists of some pairs from $\Psi$.
The requirement that $|M_1|$ is even is equivalent to the condition
$|M^*|\equiv|M_*|\pmod 2$.

Suppose that $\Psi$ consists of $\delta'$ pairs for some $\delta'\leq\delta$.
So we have $2^{\delta'}$ possible choices for $M_2$.
We have $2^{\delta-\delta'}$ choices when we chose one element from each pair in $\Phi\smallsetminus\Psi$
and we have two choices to choose $z_0$ or not.
However, the requirement that $|M_1|$ is even implies that the possible choices of $M_1$ is exactly $2^{\delta-\delta'}$.
Thus the total choices for $M$ is $2^{\delta'}\cdot 2^{\delta-\delta'}=2^\delta$ and hence (i) is proved.

Suppose $\Psi\neq\Psi'$ and $\Lambda_M\in C_{\Phi,\Psi}\cap C_{\Phi,\Psi'}$ for some $M\subset Z_\rmI$.
Because $\Psi\neq\Psi'$, there is a pair $\binom{s}{t}\in\Phi$ such that $\binom{s}{t}\in\Psi$ and $\binom{s}{t}\not\in\Psi'$.
By the two conditions before Remark~\ref{0608}, $\Lambda_M\in C_{\Phi,\Psi}$ implies that $|M\cap\{s,t\}|=0,2$,
and $\Lambda_M\in C_{\Phi,\Psi'}$ implies that $|M\cap\{s,t\}|=1$.
We get a contradiction and hence (ii) is proved.

We know that $|\cals_Z|=2^{2\delta}$ from Subsection~\ref{0232},
and we have $2^\delta$ choices of $\Psi$ for a fixed arrangement $\Phi$.
Therefore (iii) follows from (i) and (ii) directly.
\qed
\end{proof}

\begin{proposition}\label{0601}
Let $\bfG=\Sp_{2n}$, $Z$ a special symbol of rank $n$ and defect $1$,
$\Phi$ an arrangement of $Z_\rmI$ and $\Psi\leq\Phi$.
Then
\[
R_{\underline c(Z,\Phi,\Psi)}=\sum_{\Lambda\in C_{\Phi,\Psi}}\rho_\Lambda.
\]
In particular, the class function $\sum_{\Lambda\in C_{\Phi,\Psi}}\rho_\Lambda$ is uniform.
\end{proposition}
\begin{proof}
Let $\Lambda_M$ be an element in $C_{\Phi,\Psi}$.
From Proposition~\ref{0319},
we know that
\[
\langle\Lambda_M,\Lambda_{\Psi'}\rangle=|M\cap\Psi'|\pmod 2.
\]
Hence by Proposition~\ref{0306},
we have
\[
\langle\rho_{\Lambda_M},R_{\Lambda_{\Psi'}}\rangle=\frac{1}{2^\delta}(-1)^{|M\cap\Psi'|}
\]
where $\delta$ is the degree of $Z$.
Then by (\ref{0514}) and (\ref{0302}),
we have
\[
\langle\rho_{\Lambda_M},R_{\underline c}\rangle
=\frac{1}{2^\delta}\sum_{\Psi'\leq\Phi}(-1)^{|(\Phi\smallsetminus\Psi)\cap\Psi'^*|}(-1)^{|M\cap\Psi'|}
=\frac{1}{2^\delta}\sum_{\Psi'\leq\Phi}1
=1.
\]
This means $\rho_\Lambda$ occurs with multiplicity one in $R_{\underline c}$
for each $\Lambda\in C_{\Phi,\Psi}$.
From \cite{lg-symplectic} theorem 5.6 we know that $R_{\underline c}$ is a sum of $2^\delta$ distinct
irreducible characters of $G$ and $C_{\Phi,\Psi}$ has also $2^\delta$ elements.
So all components of $R_{\underline c}$ are exactly from all elements in $C_{\Phi,\Psi}$.
\qed
\end{proof}

\begin{remark}
Note that in the original theorem~5.6 in \cite{lg-symplectic} the cardinality $q$ of the base field
is assume to be large, however, according the comment by the end of \cite{lg-orthogonal}
the restriction is removed by a result of Asai.
\end{remark}

If $\Psi=\Phi$,
then $(-1)^{|(\Phi\smallsetminus\Psi)\cap\Psi'^*|}=1$ for any $\Psi'\leq\Phi$,
and the identity in Proposition~\ref{0601} becomes
\begin{equation}
\sum_{\Psi'\leq\Phi}R_{\Lambda_{\Psi'}}=\sum_{\Psi'\leq\Phi}\rho_{\Lambda_{\Psi'}}.
\end{equation}

\begin{lemma}\label{0602}
Suppose $Z$ is a special symbol of defect $1$ with singles
$Z_\rmI=\binom{s_1,s_2,\ldots,s_{\delta+1}}{t_1,t_2,\ldots,t_\delta}$.
Let $\Phi_1,\Phi_2$ be two arrangements of $Z_\rmI$ given by
\[
\textstyle
\Phi_1=\{\binom{s_1}{t_1},\binom{s_2}{t_2},\cdots,\binom{s_\delta}{t_\delta},\binom{s_{\delta+1}}{-}\}
\text{\quad and\quad }
\Phi_2=\{\binom{s_1}{-},\binom{s_2}{t_1},\binom{s_3}{t_2},\cdots,\binom{s_{\delta+1}}{t_\delta}\}.
\]
Then for any $\Psi_1\leq\Phi_1$ and $\Psi_2\leq\Phi_2$,
we have $|C_{\Phi_1,\Psi_1}\cap C_{\Phi_2,\Psi_2}|=1$.
\end{lemma}
\begin{proof}
Let $\Psi_1\leq\Phi_1$ and $\Psi_2\leq\Phi_2$.
Suppose that $\Lambda_M$ is in the intersection
$C_{\Phi_1,\Psi_1}\cap C_{\Phi_2,\Psi_2}$
for $M\subset Z_\rmI$ such that $|M^*|\equiv |M_*|\pmod 2$.
From the two conditions before Remark~\ref{0608}, we have the following:
\begin{enumerate}
\item[(1)] if $s_i\in M$ and $\binom{s_i}{t_i}\leq\Psi_1$, then $t_i\in M$;

\item[(2)] if $s_i\not\in M$ and $\binom{s_i}{t_i}\leq\Psi_1$, then $t_i\not\in M$;

\item[(3)] if $s_i\in M$ and $\binom{s_i}{t_i}\not\leq\Psi_1$, then $t_i\not\in M$;

\item[(4)] if $s_i\not\in M$ and $\binom{s_i}{t_i}\not\leq\Psi_1$, then $t_i\in M$;

\item[(5)] if $t_i\in M$ and $\binom{s_{i+1}}{t_i}\leq\Psi_2$, then $s_{i+1}\in M$;

\item[(6)] if $t_i\not\in M$ and $\binom{s_{i+1}}{t_i}\leq\Psi_2$, then $s_{i+1}\not\in M$;

\item[(7)] if $t_i\in M$ and $\binom{s_{i+1}}{t_i}\not\leq\Psi_2$, then $s_{i+1}\not\in M$;

\item[(8)] if $t_i\not\in M$ and $\binom{s_{i+1}}{t_i}\not\leq\Psi_2$, then $s_{i+1}\in M$
\end{enumerate}
for $i=1,\ldots,\delta$.
This means that for any fixed $\Psi_1,\Psi_2$,
the set $M$ is uniquely determined by the ``initial condition'' whether $s_1$ belongs to $M$ or not.
So now there are two possible choices of $M$ one of which contains $s_1$ and the other does not.
Moreover, from (1)--(8) above,
it is easy to see that both possible choices of $M$ are complement subsets to each other in $Z_\rmI$.
Therefore, only a unique $M$ satisfies the condition $|M^*|\equiv |M_*|\pmod 2$ and hence the lemma is proved.
\qed
\end{proof}

\begin{lemma}\label{0610}
Let $\Phi_1,\Phi_2$ be the two arrangements of $Z_\rmI$ given in Lemma~\ref{0602}.
For any given $\Lambda\in\cals_Z$,
there exist $\Psi_1\leq\Phi_1$ and $\Psi_2\leq\Phi_2$ such that
$C_{\Phi_1,\Psi_1}\cap C_{\Phi_2,\Psi_2}=\{\Lambda\}$.
\end{lemma}
\begin{proof}
Let $\Lambda\in\cals_Z$,
and let $\Phi_1,\Phi_2$ be the two arrangements given in Lemma~\ref{0602}.
By (3) of Lemma~\ref{0303}, there is a subset of pairs $\Psi_i$ of $\Phi_i$
such that $\Lambda\in C_{\Phi_i,\Psi_i}$ for $i=1,2$.
Therefore $\Lambda\in C_{\Phi_1,\Psi_1}\cap C_{\Phi_2,\Psi_2}$.
Then the lemma follows from Lemma~\ref{0602} immediately.
\qed
\end{proof}

\begin{lemma}\label{0611}
Let $\Lambda_1,\Lambda_2$ be two distinct symbols in $\cals_Z$.
There exists an arrangement $\Phi$ of $Z_\rmI$ with two subsets of pairs $\Psi_1,\Psi_2$
such that $\Lambda_i\in C_{\Phi,\Psi_i}$ for $i=1,2$ and
$C_{\Phi,\Psi_1}\cap C_{\Phi,\Psi_2}=\emptyset$.
\end{lemma}
\begin{proof}
Suppose that $\Lambda_1=\Lambda_{M_1}$ and $\Lambda_2=\Lambda_{M_2}$
for $M_1,M_2\subset Z_\rmI$.
Because $M_1\neq M_2$ and $|(M_i)^*|\equiv |(M_i)_*|\pmod 2$ for $i=1,2$,
it is clear that we can find a pair $\Psi=\binom{s}{t}$ such that
one of $M_1,M_2$ contains exactly one of the two elements $s,t$ and the other set contains either
both $s,t$ or none, i.e.,
\begin{equation}\label{0603}
|M_1\cap\Psi|\not\equiv|M_2\cap\Psi|\pmod 2.
\end{equation}
Let $\Phi$ be any arrangement of $Z$ that contains $\Psi$ as a subset of pairs.
By (3) of Lemma~\ref{0303},
we know that $\Lambda_{M_1}\in C_{\Phi,\Psi_1}$ and
$\Lambda_{M_2}\in C_{\Phi,\Psi_2}$ for some subsets of pairs
$\Psi_1,\Psi_2$ of $\Phi$.
Then by Lemma~\ref{0304} and (\ref{0603}) we see that $\Psi_1\neq\Psi_2$.
Finally, by (2) of Lemma~\ref{0303},
we know that $C_{\Phi,\Psi_1}\cap C_{\Phi,\Psi_2}=\emptyset$.
\qed
\end{proof}

We need stronger versions of Lemma~\ref{0610} and Lemma~\ref{0611}.

\begin{lemma}\label{0519}
Let $\Psi_0$ be a set of consecutive pairs in $Z_\rmI$.
For any given $\Lambda\in\cals_Z^{\Psi_0}$,
there exist two arrangements $\Phi_1,\Phi_2$ of $Z_\rmI$ with subsets of pairs $\Psi_1,\Psi_2$ respectively
such that $\Psi_0\leq\Psi_i\leq\Phi_i$ and
\[
C^\natural_{\Phi_1,\Psi_1}\cap C^\natural_{\Phi_2,\Psi_2}=\{\Lambda\}
\]
where $C^\natural_{\Phi_i,\Psi_i}=C_{\Phi_i,\Psi_i}\cap\cals_Z^{\Psi_0}$ for $i=1,2$.
\end{lemma}
\begin{proof}
Because $\Psi_0$ is a set of consecutive pairs in $Z_\rmI$,
the symbol $Z'$ given by $Z'=Z\smallsetminus\Psi_0$ is still a special symbol of the same defect
and $Z'_\rmI=Z_\rmI\smallsetminus\Psi_0$.
Because $\Lambda\in\cals_Z^{\Psi_0}$,
we can write $\Lambda=\Lambda'\cup\Psi_0$ (\cf.~Subsection~\ref{0233}) for a unique $\Lambda'\in\cals_{Z'}$.
Write $Z'_\rmI=\binom{s'_1,s'_2,\ldots,s'_{\delta_1+1}}{t'_1,t'_2,\ldots,t'_{\delta_1}}$ and define
\[
\textstyle
\Phi'_1=\{\binom{s'_1}{t'_1},\binom{s'_2}{t'_2},\cdots,\binom{s'_{\delta_1}}{t'_{\delta_1}},\binom{s'_{\delta_1+1}}{-}\}
\text{\quad and\quad }
\Phi'_2=\{\binom{s'_1}{-},\binom{s'_2}{t'_1},\binom{s'_3}{t'_2},\cdots,\binom{s'_{\delta_1+1}}{t'_{\delta_1}}\}.
\]
By Lemma~\ref{0610},
we know that there exist sets of pairs $\Psi_1',\Psi'_2$ of $\Phi_1',\Phi_2'$ respectively
such that
\[
C_{\Phi'_1,\Psi_1'}\cap C_{\Phi'_2,\Psi'_2}=\{\Lambda'\}.
\]
Now $\Psi_0$ itself can be regarded as an arrangement of itself, so
we have
\[
C_{\Psi_0,\Psi_0}=\{\,\Lambda_{N_2}\in\cals_{\Psi_0}\mid N_2\leq\Psi_0\,\}.
\]

Now let $\Phi_i=\Phi'_i\cup\Psi_0$, $\Psi_i=\Psi'_i\cup\Psi_0$ for $i=1,2$,
so we have $\Psi_0\leq\Psi_i\leq\Phi_i$ for $i=1,2$.
From Lemma~\ref{0515}, we can see that
\[
C_{\Phi_i,\Psi_i}=\{\,\Lambda_1\cup\Lambda_2\mid\Lambda_1\in C_{\Phi'_i,\Psi'_i},\ \Lambda_2\in C_{\Psi_0,\Psi_0}\,\}.
\]
Therefore
\[
C_{\Phi_1,\Psi_1}\cap C_{\Phi_2,\Psi_2}
=\{\,\Lambda'\cup\Lambda_2\mid\Lambda_2\in C_{\Psi_0,\Psi_0}\,\}
\]
and hence
\[
C^\natural_{\Phi_1,\Psi_1}\cap C^\natural_{\Phi_2,\Psi_2}
=C_{\Phi_1,\Psi_1}\cap C_{\Phi_2,\Psi_2}\cap\cals_Z^{\Psi_0}
=\{\Lambda'\cup\Psi_0\}
=\{\Lambda\}.
\]
\qed
\end{proof}

\begin{lemma}\label{0520}
Let $\Psi_0$ be a set of consecutive pairs in $Z_\rmI$.
Let $\Lambda_1,\Lambda_2$ be two distinct symbols in $\cals_Z^{\Psi_0}$.
There exists an arrangement $\Phi$ of $Z_\rmI$ with two subsets of pairs $\Psi_1,\Psi_2$
such that $\Psi_0\leq\Psi_i$ and $\Lambda_i\in C_{\Phi,\Psi_i}$ for $i=1,2$,
and $C_{\Phi,\Psi_1}\cap C_{\Phi,\Psi_2}=\emptyset$.
\end{lemma}
\begin{proof}
Let $Z'$ be define in the previous lemma.
Then we known that $\Lambda_i=\Lambda_i'\cup\Psi_0$ for $\Lambda'_i\in\cals_{Z'}$.
Clearly, $\Lambda'_1,\Lambda'_2$ are distinct.
Then by Lemma~\ref{0611},
we know that there is an arrangement $\Phi'$ of $Z'$ with subsets of pairs $\Psi'_1,\Psi'_2$
such that $\Lambda'_i\in C_{\Phi',\Psi'_i}$ for $i=1,2$ and
$C_{\Phi',\Psi'_1}\cap C_{\Phi',\Psi'_2}=\emptyset$.
Let $\Phi_i=\Phi'_i\cup\Psi_0$, $\Psi_i=\Psi'_i\cup\Psi_0$ for $i=1,2$.
Then as in the proof of the previous lemma,
we can see that
\[
C_{\Phi,\Psi_1}\cap C_{\Phi,\Psi_2}
=\{\,\Lambda_1\cup\Lambda_2\mid\Lambda_1\in C_{\Phi'_1,\Psi'_1}\cap C_{\Phi'_2,\Psi'_2},\ \Lambda_2\in C_{\Psi_0,\Psi_0}\,\}
=\emptyset.
\]
\qed
\end{proof}

It is clear that if $\Psi_0=\emptyset$,
then Lemma~\ref{0519} and Lemma~\ref{0520} are reduced to Lemma~\ref{0610} and Lemma~\ref{0611} respectively.

\subsection{Cells for even orthogonal groups}
In this subsection, let $\bfG=\rmO_{2n}^\epsilon$, $Z$ a special symbol of rank $n$ and defect $0$,
$\Phi$ an arrangement of $Z_\rmI$, and $\Psi\leq\Phi$.

\begin{lemma}\label{0614}
If\/ $\Phi\smallsetminus\Psi$ consists of even number of pairs, then $C_{\Phi,\Psi}\subset\cals^+_Z$;
on the other hand,
if\/ $\Phi\smallsetminus\Psi$ consists of odd number of pairs, then $C_{\Phi,\Psi}\subset\cals^-_Z$.
\end{lemma}
\begin{proof}
Suppose that $\Lambda_M\in C_{\Phi,\Psi}$ for some $M\subset Z_\rmI$.
Then from the condition before Remark~\ref{0608},
we know that $M$ contains exactly one element from each pair in $\Phi\smallsetminus\Psi$
and contains either none or two elements in each pair in $\Psi$.
This implies that $|M^*|\not\equiv |M_*|\pmod 2$ if $\Phi\smallsetminus\Psi$ consists of odd number of pairs;
and $|M^*|\equiv |M_*|\pmod 2$ if $\Phi\smallsetminus\Psi$ consists of even number of pairs.
Hence the lemma follows from (\ref{0234}).
\qed
\end{proof}

\begin{example}
Suppose that $Z=\binom{5,3,1}{4,2,0}$, $\Phi=\{\binom{5}{4},\binom{3}{2},\binom{1}{0}\}$,
and $\Psi=\{\binom{5}{4},\binom{1}{0}\}$.
To construct a subset $M$ of $Z_\rmI$,
we need to choose one element from each pair in $\Phi\smallsetminus\Psi=\{\binom{3}{2}\}$ and
choose a subset of pairs of $\Psi$.
Hence we have $8$ possible subsets $M$, namely,
$\binom{3}{-}$, $\binom{3,1}{0}$, $\binom{5,3}{4}$, $\binom{5,3,1}{4,0}$,
$\binom{-}{2}$, $\binom{1}{2,0}$, $\binom{5}{4,2}$, $\binom{5,1}{4,2,0}$.
Hence $C_{\Phi,\Psi}$ consists of the following $8$ elements:
$\binom{5,1}{4,3,2,0}$, $\binom{5,0}{4,3,2,1}$, $\binom{4,1}{5,3,2,0}$, $\binom{4,0}{5,3,2,1}$,
$\binom{5,3,2,1}{4,0}$, $\binom{5,3,2,0}{4,1}$, $\binom{4,3,2,1}{5,0}$, $\binom{4,3,2,0}{5,1}$.
Note that $\Phi\smallsetminus\Psi$ consists of one pair,
so $C_{\Phi,\Psi}\subset\cals_Z^-$.
\end{example}

\begin{lemma}\label{0606}
Let $Z$ be a special symbol of defect $0$ and degree $\delta\geq 1$,
$\Phi$ a fixed arrangement of $Z_\rmI$, and $\Psi,\Psi'$ subsets of pairs of\/ $\Phi$.
Then
\begin{enumerate}
\item[(i)] $\Lambda\in C_{\Phi,\Psi}$ if and only if $\Lambda^\rmt\in C_{\Phi,\Psi}$;

\item[(ii)] $|C_{\Phi,\Psi}|=2^\delta$;

\item[(iii)] $C_{\Phi,\Psi}\cap C_{\Phi,\Psi'}=\emptyset$ if\/ $\Psi\neq\Psi'$;

\item[(iv)] we have
\[
\cals^+_Z=\bigcup_{\Psi\leq\Phi,\ |\Phi\smallsetminus\Psi|\text{\ even}}C_{\Phi,\Psi}\quad\text{and}\quad
\cals^-_Z=\bigcup_{\Psi\leq\Phi,\ |\Phi\smallsetminus\Psi|\text{\ odd}}C_{\Phi,\Psi}
\]
where $|\Phi\smallsetminus\Psi|$ means the number of pairs in $\Phi\smallsetminus\Psi$.
\end{enumerate}
\end{lemma}
\begin{proof}
Suppose that $\Lambda\in C_{\Phi,\Psi}$ and write $\Lambda=\Lambda_M$ for some $M\subset Z_\rmI$.
Then it is easy to check that $\Lambda^\rmt=\Lambda_{Z_\rmI\smallsetminus M}$.
It is clear that $M$ satisfies the condition that it consists of exactly one element from each pair in $\Phi\smallsetminus\Psi$
and a subset of pairs of $\Psi$ if and only if $Z_\rmI\smallsetminus M$ satisfies the same condition.
Hence (i) is proved.

From the conditions before Remark~\ref{0608},
we can write $M=M_1\sqcup M_2$ where $M_1$ consists of exactly one element from each pair
of $\Phi\smallsetminus\Psi$,
and $M_2$ consists of some pairs from $\Psi$.
Suppose that $\Psi$ contains $\delta'$ pairs for some $\delta'\leq\delta$.
So we have $2^{\delta'}$ possible choices for $M_2$ and $2^{\delta-\delta'}$ choices for $M_1$.
Thus the total choices for $M$ is $2^{\delta'}\cdot 2^{\delta-\delta'}=2^\delta$ and hence (ii) is proved.

The proof of (iii) is similar to that of Lemma~\ref{0303}.

For any fixed arrangement $\Phi$ of $Z_\rmI$, we have
\[
\cals_Z^+\sqcup\cals_Z^-=\bigsqcup_{\Psi\leq\Phi}C_{\Phi,\Psi}
\]
by the same argument of the proof of Lemma~\ref{0303}.
Then (iv) follows from Lemma~\ref{0614} immediately.
\qed
\end{proof}

A subset of pairs $\Psi$ of an arrangement $\Phi$ of $Z_\rmI$ is called \emph{admissible}
if $|\Phi\smallsetminus\Psi|$ is even when $\epsilon=+$ and
$|\Phi\smallsetminus\Psi|$ is odd when $\epsilon=-$.

\begin{proposition}\label{0612}
Let $\bfG=\rmO^\epsilon_{2n}$, $Z$ a special symbol of rank $n$ and defect $0$,
$\Phi$ an arrangement of $Z_\rmI$ with an admissible subset of pairs $\Psi$.
Then
\[
R^{\rmO^\epsilon}_{\underline c(Z,\Phi,\Psi)}
=\sum_{\Lambda\in C_{\Phi,\Psi}}\rho^{\rmO^\epsilon}_\Lambda
\]
\end{proposition}
\begin{proof}
If $\epsilon=+$ and $Z$ is of degree $0$, i.e., $Z$ is degenerate,
then it is clear that $\Phi=\Psi=\emptyset$, $\underline c(Z,\Phi,\Psi)=C_{\Phi,\Psi}=\{Z\}$ and
$R_Z^{\rmO^+}=\rho_Z^{\rmO^+}$.
If $\epsilon=-$ and $Z$ degenerate,
then $\underline c(Z,\Phi,\Psi)=C_{\Phi,\Psi}=\emptyset$.
So the proposition holds if $Z$ is degenerate.

Now suppose that $Z$ is of degree greater than or equal to $1$.
Let $C_{\Phi,\Psi}^{\rm SO^\epsilon}$ be a subset of $C_{\Phi,\Psi}$ such that
$C_{\Phi,\Psi}^{\rm SO^\epsilon}$ contains exactly one element from each pair $\{\Lambda,\Lambda^\rmt\}\subset C_{\Phi,\Psi}$.
Therefore $|C_{\Phi,\Psi}^{\rm SO^\epsilon}|=2^{\delta-1}$ by (ii) of Lemma~\ref{0606}.
By the argument in the proof of Proposition~\ref{0601},
we can show that $\langle\rho_\Lambda,R^{\SO^\epsilon}_{\underline c(Z,\Phi,\Psi)}\rangle_{\SO^\epsilon}=1$
for every $\Lambda\in C_{\Phi,\Psi}$.
Moreover, we know that $R^{\SO^\epsilon}_{\underline c(Z,\Phi,\Psi)}$ is a sum of $2^{\delta-1}$ distinct
irreducible characters of $\SO^\epsilon(q)$ by \cite{lg-orthogonal} proposition 3.13.
Thus we have
\[
R^{\rm SO^\epsilon}_{\underline c(Z,\Phi,\Psi)}
=\sum_{\Lambda\in C^{\rm SO^\epsilon}_{\Phi,\Psi}}\rho^{\rm SO^\epsilon}_\Lambda.
\]
Therefore
\begin{align*}
R^{\rmO^\epsilon}_{\underline c(Z,\Phi,\Psi)}
=\Ind_{\rm SO^\epsilon}^{\rm O^\epsilon}R^{\rm SO^\epsilon}_{\underline c(Z,\Phi,\Psi)}
=\sum_{\Lambda\in C^{\rm SO^\epsilon}_{\Phi,\Psi}}\Ind_{\rm SO^\epsilon}^{\rm O^\epsilon}\rho^{\rm SO^\epsilon}_\Lambda
&=\sum_{\Lambda\in C^{\rm SO^\epsilon}_{\Phi,\Psi}}(\rho^{\rmO^\epsilon}_\Lambda+\rho^{\rmO^\epsilon}_{\Lambda^\rmt})\\
&=\sum_{\Lambda\in C_{\Phi,\Psi}}\rho^{\rmO^\epsilon}_\Lambda.
\end{align*}
\qed
\end{proof}

\begin{lemma}\label{0516}
Let $\Lambda_1,\Lambda_2$ be two symbols in $\cals^\epsilon_Z$ such that
$\Lambda_1\neq\Lambda_2,\Lambda_2^\rmt$.
There exists an arrangement $\Phi$ of $Z_\rmI$ with admissible subsets of pairs $\Psi_1,\Psi_2$
such that $\Lambda_i,\Lambda_i^\rmt\in C_{\Phi,\Psi_i}$ for $i=1,2$ and
$C_{\Phi,\Psi_1}\cap C_{\Phi,\Psi_2}=\emptyset$.
\end{lemma}
\begin{proof}
Suppose that $\Lambda_1=\Lambda_{M_1}$ and $\Lambda_2=\Lambda_{M_2}$ for $M_1,M_2\subset Z_\rmI$.
The assumption that $\Lambda_1\neq\Lambda_2,\Lambda_2^\rmt$ means that
$M_1\neq M_2$ and $M_1\neq (Z_\rmI\smallsetminus M_2)$.
Then it is clear that we can find a pair $\Psi=\binom{s}{t}$ in $Z_\rmI$ such that
one of $M_1,M_2$ contains exactly one of the two elements $s,t$ and the other set contains either
both $s,t$ or none, i.e.,
\begin{equation}\label{0607}
|M_1\cap\Psi|\not\equiv|M_2\cap\Psi|\pmod 2.
\end{equation}
Let $\Phi$ be any arrangement of $Z_\rmI$ that contains $\Psi$ as a subset of pairs.
By (iv) of Lemma~\ref{0606},
we know that $\Lambda_{M_1}\in C_{\Phi,\Psi_1}$ and
$\Lambda_{M_2}\in C_{\Phi,\Psi_2}$ for some subsets of pairs
$\Psi_1,\Psi_2$ of $\Phi$.
Then by Lemma~\ref{0604} and (\ref{0607}) we see that $\Psi_1\neq\Psi_2$.
Finally, by (iii) of Lemma~\ref{0606},
we know that $C_{\Phi,\Psi_1}\cap C_{\Phi,\Psi_2}=\emptyset$.
\qed
\end{proof}

\begin{lemma}\label{0518}
Let $\Psi_0$ be a set of consecutive pairs in $Z_\rmI$.
Let $\Lambda_1,\Lambda_2$ be two symbols in $\cals^{\epsilon,\Psi_0}_Z$ such that
$\Lambda_1\neq\Lambda_2,\Lambda_2^\rmt$.
There exists an arrangement $\Phi$ of $Z_\rmI$ with subsets $\Psi_1,\Psi_2$
such that $\Psi_0\leq\Psi_i$ and $\Lambda_i,\Lambda_i^\rmt\in C_{\Phi,\Psi_i}$ for $i=1,2$, and
$C_{\Phi,\Psi_1}\cap C_{\Phi,\Psi_2}=\emptyset$.
\end{lemma}
\begin{proof}
The proof of the lemma is similar to that of Lemma~\ref{0520}
except that we need to apply Lemma~\ref{0516} instead of Lemma~\ref{0611}.
\qed
\end{proof}

If $\Psi_0=\emptyset$,
then $\cals_Z^{\epsilon,\Psi_0}=\cals_Z^\epsilon$ and Lemma~\ref{0518} is reduced to Lemma~\ref{0516}.


\section{A System of Linear Equations}\label{1004}

\subsection{Decomposition with respect special symbols}
Let $(\bfG,\bfG')=(\Sp_{2n},\rmO^\epsilon_{2n'})$,
and let $Z,Z'$ be special symbols of ranks $n,n'$, degrees $\delta,\delta'$, and defects $1,0$ respectively.
Let $\omega_{Z,Z'}$ denote the orthogonal projection of $\omega_{\bfG,\bfG',1}$ over $\calv_Z\otimes\calv^\epsilon_{Z'}$.
Then by Proposition~\ref{0306} and Proposition~\ref{0319} we have
\[
\omega_{\bfG,\bfG',1}=\sum_{Z,Z'}\omega_{Z,Z'}\quad\text{ and }\quad
\omega_{\bfG,\bfG',1}^\sharp=\sum_{Z,Z'}\omega^\sharp_{Z,Z'}
\]
where $Z,Z'$ run over all special symbols of rank $n,n'$ and defect $1,0$ respectively.
Moreover, because $\calb_{\bfG,\bfG'}=\bigsqcup_{Z,Z'}\calb^\epsilon_{Z,Z'}$, we have
\[
\sum_{(\Lambda,\Lambda')\in\calb_{\bfG,\bfG'}}\rho_\Lambda\otimes\rho_{\Lambda'}
=\sum_{Z,Z'}\sum_{(\Lambda,\Lambda')\in\calb^\epsilon_{Z,Z'}}\rho_\Lambda\otimes\rho_{\Lambda'}.
\]
Now Proposition~\ref{0310} implies that
\begin{equation}\label{0615}
\omega_{Z,Z'}^\sharp
=\sum_{(\Lambda,\Lambda')\in\calb^\epsilon_{Z,Z'}}(\rho_\Lambda\otimes\rho_{\Lambda'})^\sharp.
\end{equation}
Then, for any uniform class function $f$ of $G\times G'$,
we have
\begin{align}\label{1022}
\langle f,\omega_{Z,Z'}\rangle
=\langle f,\omega_{Z,Z'}^\sharp\rangle
=\biggl\langle f,
\sum_{(\Lambda,\Lambda')\in\calb^\epsilon_{Z,Z'}}(\rho_\Lambda\otimes\rho_{\Lambda'})^\sharp\biggr\rangle
=\biggl\langle f,\sum_{(\Lambda,\Lambda')\in\calb^\epsilon_{Z,Z'}}\rho_\Lambda\otimes\rho_{\Lambda'}\biggr\rangle.
\end{align}
Now the candidates of the uniform class functions will be those construct from cells
described in Section~\ref{0609}.

\subsection{Case for $\cald_{Z,Z'}$ one-to-one, $\deg(Z')=\deg(Z)+1$}\label{1010}
In this subsection,
we assume that $\cald_{Z,Z'}$ is one-to-one.
Now $\delta,\delta'$ are the degrees of $Z,Z'$ respectively,
so we can write
\begin{equation}\label{1007}
Z_\rmI=\binom{a_1,a_2,\ldots,a_{\delta+1}}{b_1,b_2,\ldots,b_\delta},\qquad
Z'_\rmI=\binom{c_1,c_2,\ldots,c_{\delta'}}{d_1,d_2,\ldots,d_{\delta'}}.
\end{equation}
Then from the remark~7.17 in \cite{pan-uniform},
we know that either $\delta'=\delta+1$ or $\delta'=\delta$.
In this subsection, we consider the case that $\delta'=\delta+1$.

Define
\begin{align}\label{0501}
\begin{split}
\theta\colon\{a_1,\ldots,a_{\delta+1}\}\sqcup\{b_1,\ldots,b_\delta\} &\rightarrow \{c_1,\ldots,c_{\delta+1}\}\sqcup\{d_1,\ldots, d_{\delta+1}\} \\
a_i &\mapsto d_i \\
b_i &\mapsto c_{i+1}
\end{split}
\end{align}
for each $i$.
Note that $c_1$ is not in the image of $\theta$.
Then $\theta$ induces a mapping $\theta^\epsilon\colon\cals_Z\rightarrow\cals^\epsilon_{Z'}$
given by
\[
\Lambda_M\mapsto\begin{cases}
\Lambda_{\theta(M)}, & \text{if $\epsilon=+$};\\
\Lambda_{\binom{c_1}{-}\cup\theta(M)}, & \text{if $\epsilon=-$}
\end{cases}
\]
where $M\subset Z_\rmI$.

If $\Phi=\{\binom{s_1}{t_1},\ldots,\binom{s_\delta}{t_\delta},\binom{s_{\delta+1}}{-}\}$ is an arrangement of $Z_\rmI$
(i.e., $\{s_1,\ldots,s_{\delta+1}\}$ (resp.\ $\{t_1,\ldots,t_\delta\}$) is a permutation of $\{a_1,\ldots,a_{\delta+1}\}$
(resp.~$\{b_1,\ldots,b_\delta\}$)),
then
\[
\theta(\Phi)=\left\{\binom{\theta(t_1)}{\theta(s_1)},\ldots,\binom{\theta(t_\delta)}{\theta(s_\delta)},\binom{c_1}{\theta(s_{\delta+1})}\right\}
\]
is an arrangement of $Z'_\rmI$.
If $\Psi=\{\binom{s_{i_1}}{t_{i_1}},\ldots,\binom{s_{i_k}}{t_{i_k}}\}$ is a subset of pairs of $\Phi$,
we define $\theta(\Psi)$ as follows:
\begin{enumerate}
\item[(1)] if either $\epsilon=-$ and $|\Phi\smallsetminus\Psi|$ is odd, or $\epsilon=+$ and $|\Phi\smallsetminus\Psi|$ is even,
let
\[
\theta(\Psi)=\Biggl\{\binom{\theta(t_{i_1})}{\theta(s_{i_1})},\ldots,\binom{\theta(t_{i_k})}{\theta(s_{i_k})}\Biggr\};
\]

\item[(2)] if either $\epsilon=-$ and $|\Phi\smallsetminus\Psi|$ is even, or $\epsilon=+$ and $|\Phi\smallsetminus\Psi|$ is odd,
let
\[
\theta(\Psi)=\Biggl\{\binom{\theta(t_{i_1})}{\theta(s_{i_1})},\ldots,\binom{\theta(t_{i_k})}{\theta(s_{i_k})},
\binom{c_1}{\theta(s_{\delta+1})}\Biggr\}.
\]
\end{enumerate}
Then $\theta(\Psi)$ is a subset of pairs of $\theta(\Phi)$ such that $|\theta(\Phi)\smallsetminus\theta(\Psi)|$ is even if $\epsilon=+$;
$|\theta(\Phi)\smallsetminus\theta(\Psi)|$ is odd if $\epsilon=-$.

\begin{lemma}\label{1008}
Suppose that $\cald_{Z,Z'}$ is one-to-one and $\deg(Z')=\deg(Z)+1$.
Let $\Phi$ be an arrangement of $Z_\rmI$ and $\Psi$ a subset of pairs of\/ $\Phi$.
Then
\[
C_{\theta(\Phi),\theta(\Psi)}
=\{\,\theta^\epsilon(\Lambda),\theta^\epsilon(\Lambda)^\rmt\mid
\Lambda\in C_{\Phi,\Psi}\,\}
\]
where $C_{\Phi,\Psi}$ is defined in (\ref{0302}).
\end{lemma}
\begin{proof}
As above, let $s_{\delta+1}$ be the isolated element in $\Phi$.
Suppose that $\Lambda=\Lambda_M\in C_{\Phi,\Psi}$ for some $M\subset Z_\rmI$ and $\theta^\epsilon(\Lambda)=\Lambda_{M'}$.
From the rules before Remark~\ref{0608}, we know that $M$ contains exactly one element from each pair of $\Phi\smallsetminus\Psi$
and contains some subset of pairs in $\Psi$.
Moreover, $M$ contains the isolated element $s_{\delta+1}$ if and only if $|\Phi\smallsetminus\Psi|$ is odd.
Then
\begin{enumerate}
\item[(1)] if $\epsilon=+$ and $|\Phi\smallsetminus\Psi|$ is even,
then $s_{\delta+1}\not\in M$ and $M'=\theta(M)$;

\item[(2)] if $\epsilon=+$ and $|\Phi\smallsetminus\Psi|$ is odd,
then $s_{\delta+1}\in M$ and $M'=\theta(M)$;

\item[(3)] if $\epsilon=-$ and $|\Phi\smallsetminus\Psi|$ is even,
then $s_{\delta+1}\not\in M$ and $M'=\binom{c_1}{-}\cup\theta(M)$;

\item[(4)] if $\epsilon=-$ and $|\Phi\smallsetminus\Psi|$ is odd,
then $s_{\delta+1}\in M$ and $M'=\binom{c_1}{-}\cup\theta(M)$.
\end{enumerate}
It is easy to see from the definition above that for each case above
$M'$ consists of exactly one element from each pair in
$\theta(\Phi)\smallsetminus\theta(\Psi)$ and a subset of $\theta(\Psi)$, i.e.,
$\theta^\epsilon(\Lambda)=\Lambda_{M'}\in C_{\theta(\Phi),\theta(\Psi)}$.

From (1) of Lemma~\ref{0606},
we know that
\[
\theta^\epsilon(\Lambda)\in C_{\theta(\Phi),\theta(\Psi)}\text{ if and only if }
\theta^\epsilon(\Lambda)^\rmt\in C_{\theta(\Phi),\theta(\Psi)}.
\]
So we have
\begin{equation}\label{1023}
\{\,\theta^\epsilon(\Lambda),\theta^\epsilon(\Lambda)^\rmt\mid
\Lambda\in C_{\Phi,\Psi}\,\}\subseteq C_{\theta(\Phi),\theta(\Psi)}.
\end{equation}
Now $Z$ is of degree $\delta$, so $|C_{\Phi,\Psi}|=2^\delta$ by Lemma~\ref{0303};
$Z'$ is of degree $\delta+1$, so $|C_{\theta(\Phi),\theta(\Psi)}|=2^{\delta+1}$ by Lemma~\ref{0606}.
Hence both sets in (\ref{1023}) have the same cardinality $2^{\delta+1}$, they must be the same.
\end{proof}

\begin{proposition}\label{1021}
Let $(\bfG,\bfG')=(\Sp_{2n},\rmO^\epsilon_{2n'})$,
and let $Z,Z'$ be special symbols of ranks $n,n'$ and defects $1,0$ respectively.
Suppose that $\cald_{Z,Z'}$ is one-to-one and $\deg(Z')=\deg(Z)+1$.
Then
\[
\omega_{Z,Z'}=\sum_{\Lambda\in\cals_Z}\rho_\Lambda\otimes\rho_{f(\Lambda)}
\]
where $f(\Lambda)$ is either equal to $\theta^\epsilon(\Lambda)$ or $\theta^\epsilon(\Lambda)^\rmt$.
\end{proposition}
\begin{proof}
Because we assume that $\cald_{Z,Z'}$ is one-to-one and $\deg(Z')=\deg(Z)+1$,
we know that $\calb_{Z,Z'}^\epsilon=\{\,(\Lambda,\theta^\epsilon(\Lambda))\mid\Lambda\in\cals_Z\,\}$
and (\ref{0615}) becomes
\begin{equation}\label{0616}
\omega_{Z,Z'}^\sharp=\sum_{\Lambda\in\cals_Z}(\rho_\Lambda\otimes\rho_{\theta^\epsilon(\Lambda)})^\sharp.
\end{equation}
For $\Lambda,\Lambda'\in\cals_Z$,
define
\[
x_{\Lambda,\Lambda'}
=\langle\rho_\Lambda\otimes\rho_{\theta^\epsilon(\Lambda')}+\rho_\Lambda\otimes\rho_{\theta^\epsilon(\Lambda')^\rmt},
\omega_{Z,Z'}\rangle,
\]
the sum of multiplicities of $\rho_\Lambda\otimes\rho_{\theta^\epsilon(\Lambda')}$ and
$\rho_\Lambda\otimes\rho_{\theta^\epsilon(\Lambda')^\rmt}$ in $\omega_{Z,Z'}$.
So we need to show that $x_{\Lambda,\Lambda'}=1$ if $\Lambda=\Lambda'$ and $x_{\Lambda,\Lambda'}=0$ otherwise.

Now suppose that $\Phi,\Phi'$ are any two arrangements of $Z_\rmI$,
and $\Psi,\Psi'$ are subsets of pairs of $\Phi,\Phi'$ respectively.
Then by Proposition~\ref{0601} and Proposition~\ref{0612},
the class function
\[
\sum_{\Lambda\in C_{\Phi,\Psi}}\sum_{\Lambda_1\in C_{\theta(\Phi'),\theta(\Psi')}}\rho_\Lambda\otimes\rho_{\Lambda_1}
\]
on $G\times G'$ is unform.
Then by Lemma~\ref{1008},
we have
\[
\sum_{\Lambda\in C_{\Phi,\Psi}}\sum_{\Lambda_1\in C_{\theta(\Phi'),\theta(\Psi')}}\rho_\Lambda\otimes\rho_{\Lambda_1}
=\sum_{\Lambda\in C_{\Phi,\Psi}}\sum_{\Lambda'\in C_{\Phi',\Psi'}}
(\rho_\Lambda\otimes\rho_{\theta^\epsilon(\Lambda')}+\rho_\Lambda\otimes\rho_{\theta^\epsilon(\Lambda')^\rmt}).
\]

Then, by (\ref{0616}), we have
\begin{align*}
& \sum_{\Lambda\in C_{\Phi,\Psi}}\sum_{\Lambda'\in C_{\Phi',\Psi'}} x_{\Lambda,\Lambda'} \\
& =\biggl\langle \sum_{\Lambda\in C_{\Phi,\Psi}}\sum_{\Lambda'\in C_{\Phi',\Psi'}}
(\rho_\Lambda\otimes\rho_{\theta^\epsilon(\Lambda')}+\rho_\Lambda\otimes\rho_{\theta^\epsilon(\Lambda')^\rmt}),
\omega_{Z,Z'}\biggr\rangle \\
& =\biggl\langle \sum_{\Lambda\in C_{\Phi,\Psi}}\sum_{\Lambda'\in C_{\Phi',\Psi'}}
(\rho_\Lambda\otimes\rho_{\theta^\epsilon(\Lambda')}+\rho_\Lambda\otimes\rho_{\theta^\epsilon(\Lambda')^\rmt}),
\omega_{Z,Z'}^\sharp\biggr\rangle \\
& =\biggl\langle \sum_{\Lambda\in C_{\Phi,\Psi}}\sum_{\Lambda'\in C_{\Phi',\Psi'}}
(\rho_\Lambda\otimes\rho_{\theta^\epsilon(\Lambda')}+\rho_\Lambda\otimes\rho_{\theta^\epsilon(\Lambda')^\rmt}),
\biggl(\sum_{\Lambda''\in\cals_Z}\rho_{\Lambda''}\otimes\rho_{\theta^\epsilon(\Lambda'')}\biggr)^\sharp\biggr\rangle \\
& =\biggl\langle \sum_{\Lambda\in C_{\Phi,\Psi}}\sum_{\Lambda'\in C_{\Phi',\Psi'}}
(\rho_\Lambda\otimes\rho_{\theta^\epsilon(\Lambda')}+\rho_\Lambda\otimes\rho_{\theta^\epsilon(\Lambda')^\rmt}),
\sum_{\Lambda''\in\cals_Z}\rho_{\Lambda''}\otimes\rho_{\theta^\epsilon(\Lambda'')}\biggr\rangle.
\end{align*}
From the definition of $\theta^\epsilon$ we know that $\theta^\epsilon(\Lambda'')\neq\theta^\epsilon(\Lambda')^\rmt$
for any $\Lambda'',\Lambda'\in\cals_Z$.
For a symbol $\Lambda''\in\cals_Z$ to contribute a multiplicity in the above identity,
we need $\Lambda''=\Lambda$ and $\Lambda''=\Lambda'$ for some $\Lambda\in C_{\Phi,\Psi}$ and
some $\Lambda'\in C_{\Phi',\Psi'}$, i.e.,
$\Lambda''$ must be in the intersection $C_{\Phi,\Psi}\cap C_{\Phi',\Psi'}$.
Therefore, the above equation becomes
\begin{equation}\label{1104}
\sum_{\Lambda\in C_{\Phi,\Psi}}\sum_{\Lambda'\in C_{\Phi',\Psi'}} x_{\Lambda,\Lambda'}
 =|C_{\Phi,\Psi}\cap C_{\Phi',\Psi'}|
\end{equation}
for any two arrangements $\Phi,\Phi'$ of $Z_\rmI$ with $\Psi\leq\Phi$ and $\Psi'\leq\Phi'$.

Suppose that $\Lambda_1,\Lambda_2$ are distinct symbols in $\cals_Z$.
Then by Lemma~\ref{0611}, there exist an arrangement $\Phi$ with two subsets of pairs $\Psi_1,\Psi_2$
such that $\Lambda_i\in C_{\Phi,\Psi_i}$ for $i=1,2$ and $C_{\Phi,\Psi_1}\cap C_{\Phi,\Psi_2}=\emptyset$.
Because each $x_{\Lambda,\Lambda'}$ is a non-negative integer,
from equation (\ref{1104}) we conclude that $x_{\Lambda_1,\Lambda_2}=0$ for any distinct $\Lambda_1,\Lambda_2\in\cals_Z$.

For any $\Lambda\in\cals_Z$, by Lemma~\ref{0610},
there exist two arrangements $\Phi_1,\Phi_2$ of $Z_\rmI$ with subsets of pairs $\Psi_1,\Psi_2$
respectively such that $C_{\Phi_1,\Psi_1}\cap C_{\Phi_2,\Psi_2}=\{\Lambda\}$.
Because we know $x_{\Lambda,\Lambda'}=0$ if $\Lambda\neq\Lambda'$, equation (\ref{1104}) is reduced to $x_{\Lambda,\Lambda}=1$.
Therefore, exactly one of $\rho_\Lambda\otimes\rho_{\theta^\epsilon(\Lambda)},\rho_\Lambda\otimes\rho_{\theta^\epsilon(\Lambda)^\rmt}$
occurs in $\omega_{Z,Z'}$.

For $\Lambda\in\cals_Z$,
let $f(\Lambda)$ be either $\theta^\epsilon(\Lambda)$ or $\theta^\epsilon(\Lambda)^\rmt$
such that  $\rho_\Lambda\otimes\rho_{f(\Lambda)}$ occurs in $\omega_{Z,Z'}$.
Then we just show the character $\overline\omega_{Z,Z'}$ defined by
\[
\overline\omega_{Z,Z'}=\sum_{\Lambda\in\cals_Z}\rho_\Lambda\otimes\rho_{f(\Lambda)}.
\]
is a sub-character of $\omega_{Z,Z'}$.
Note that $\rho_{\theta^\epsilon(\Lambda)}$ and $\rho_{\theta^\epsilon(\Lambda)^\rmt}$ have the same degree
for any $\Lambda\in\cals_Z$.
Therefore $\overline\omega_{Z,Z}$ and $\sum_{\Lambda\in\cals_Z}\rho_\Lambda\otimes\rho_{\theta^\epsilon(\Lambda)}$
have the same degree.
By (\ref{0616}) and (\ref{0330}),
we see that $\omega_{Z,Z'}$ and $\overline\omega_{Z,Z'}$ have the same degree.
Therefore $\overline\omega_{Z,Z'}=\omega_{Z,Z'}$ and the proposition is proved.
\end{proof}

\subsection{Case for $\cald_{Z,Z'}$ one-to-one, $\deg(Z')=\deg(Z)$}\label{1011}
Write $Z_\rmI,Z'_\rmI$ as in (\ref{1007})
and we assume that $\cald_{Z,Z'}$ is one-to-one and $\delta'=\delta$.
Define
\begin{align}\label{0502}
\begin{split}
\theta\colon\{c_1,\ldots,c_\delta\}\sqcup\{d_1,\ldots,d_\delta\} &\rightarrow \{a_1,\ldots,a_{\delta+1}\}\sqcup\{b_1,\ldots, b_\delta\} \\
c_i &\mapsto b_i \\
d_i &\mapsto a_{i+1}
\end{split}
\end{align}
for each $i$.
Note that $a_1$ is not in the image of $\theta$.
Then $\theta$ induces a mapping $\theta^\epsilon\colon\cals^\epsilon_{Z'}\rightarrow\cals_Z$
given by
\[
\Lambda_N\mapsto\begin{cases}
\Lambda_{\theta(N)}, & \text{if $\epsilon=+$};\\
\Lambda_{\binom{a_1}{-}\cup\theta(N)}, & \text{if $\epsilon=-$}
\end{cases}
\]
where $N\subset Z'_\rmI$.

If $\Phi'=\Bigl\{\binom{s_1'}{t_1'},\ldots,\binom{s_{\delta'}'}{t_{\delta'}'}\Bigr\}$ is an arrangement of $Z'_\rmI$,
then
\begin{equation}\label{1012}
\theta(\Phi')=\Biggl\{\binom{a_1}{-},\binom{\theta(t_1')}{\theta(s_1')},\ldots,\binom{\theta(t_{\delta'}')}{\theta(s_{\delta'}')}\Biggr\}
\end{equation}
is an arrangement of $Z_\rmI$.
If $\Psi'=\Bigl\{\binom{s'_{i_1}}{t'_{i_1}},\ldots,\binom{s'_{i_k}}{t'_{i_k}}\Bigr\}$ is a subset of pairs of $\Phi'$,
then
\begin{equation}\label{1013}
\theta(\Psi')=\Biggl\{\binom{\theta(t'_{i_1})}{\theta(s'_{i_1})},\ldots,\binom{\theta(t'_{i_k})}{\theta(s'_{i_k})}\Biggr\}
\end{equation}
is a subset of pairs in $\theta(\Phi')$.

\begin{lemma}\label{1009}
Suppose that $\cald_{Z,Z'}$ is one-to-one and $\deg(Z')=\deg(Z)$.
Let $\Phi'$ be an arrangement of $Z'_\rmI$ and let $\Psi'$ be a subset of pairs in $\Phi'$ such that
$|\Phi'\smallsetminus\Psi'|$ is even if $\epsilon=+$, and $|\Phi'\smallsetminus\Psi'|$ is odd if $\epsilon=-$.
Then
\[
C_{\theta(\Phi'),\theta(\Psi')}
=\{\,\theta^\epsilon(\Lambda)\mid\Lambda\in C_{\Phi',\Psi'}\,\}.
\]
\end{lemma}
\begin{proof}
Suppose that $\Lambda=\Lambda_{M'}\in C_{\Phi',\Psi'}$ for some $M'\subset Z'_\rmI$,
and $\theta(\Lambda)=\Lambda_M$.
We know that $M'=M'_1\sqcup M'_2$ where $M'_1$ consists of exactly one element from each pair in $\Phi'\smallsetminus\Psi'$
and $M'_2$ is a subset of pairs of $\Psi'$.

\begin{enumerate}
\item[(1)] First suppose that $\epsilon=+$.
Then we have $|M'^*|\equiv |M'_*|\pmod 2$.
From the above definition,
we see that $M=\theta(M')$.
Hence $M$ contains one element from each pair in $\theta(\Phi')\smallsetminus\theta(\Psi')$ and contains a subset of pairs in $\theta(\Psi')$.
Moreover, $|M^*|\equiv |M_*|\pmod 2$.
Therefore $\theta^\epsilon(\Lambda)$ is in $C_{\theta(\Phi'),\theta(\Psi')}$.

\item[(2)] Next suppose that $\epsilon=-$.
Now $|M'^*|\not\equiv |M'_*|\pmod 2$ and we see that $M=\binom{a_1}{-}\cup\theta(M')$.
Now again $M$ consists one element from each pair in $\theta(\Phi')\smallsetminus\theta(\Psi')$ and contains a subset of pairs in $\theta(\Psi')$.
Moreover, $|M^*|=|M'^*|+1$ and $|M_*|=|M'_*|$, hence $|M^*|\equiv |M_*|\pmod 2$.
Therefore $\theta^\epsilon(\Lambda)$ is in $C_{\theta(\Phi'),\theta(\Psi')}$, again.
\end{enumerate}
Now we conclude that  $\{\,\theta^\epsilon(\Lambda)\mid\Lambda\in C_{\Phi',\Psi'}\,\}\subseteq C_{\theta(\Phi'),\theta(\Psi')}$.
Since both sets have the same cardinality $2^\delta$, they must be equal.
\end{proof}

\begin{proposition}\label{1005}
Let $(\bfG,\bfG')=(\Sp_{2n},\rmO^\epsilon_{2n'})$,
and let $Z,Z'$ be special symbols of ranks $n,n'$ and defects $1,0$ respectively.
Suppose that $\cald_{Z,Z'}$ is one-to-one and $\deg(Z')=\deg(Z)$.
Then
\[
\omega_{Z,Z'}=\sum_{\Lambda'\in\cals^\epsilon_{Z'}}\rho_{\theta^\epsilon(\Lambda')}\otimes\rho_{f'(\Lambda')}
\]
where $f'(\Lambda')$ is equal to either $\Lambda'$ or $\Lambda'^\rmt$.
\end{proposition}
\begin{proof}
The proof is similar to that of Proposition~\ref{1021}.
Because we assume that $\cald_{Z,Z'}$ is one-to-one and $\deg(Z')=\deg(Z)$,
we know that $\calb_{Z,Z'}^\epsilon=\{\,(\theta^\epsilon(\Lambda'),\Lambda')\mid\Lambda'\in\cals^\epsilon_{Z'}\,\}$
and (\ref{0615}) becomes
\begin{equation}\label{0617}
\omega_{Z,Z'}^\sharp=\sum_{\Lambda'\in\cals^\epsilon_{Z'}}(\rho_{\theta^\epsilon(\Lambda')}\otimes\rho_{\Lambda'})^\sharp.
\end{equation}

For $\Lambda\in\cals_Z$ and $\Lambda'\in\cals^\epsilon_{Z'}$, define
\[
x_{\Lambda,\Lambda'}=\langle\rho_\Lambda\otimes\rho_{\Lambda'},\omega_{Z,Z'}\rangle.
\]
Then each $x_{\Lambda,\Lambda'}$ is a non-negative integer.
For an arrangement $\Phi$ of $Z_\rmI$ with a subset of pairs $\Psi$,
and an arrangement $\Phi'$ of $Z'_\rmI$ with an admissible subset of pairs $\Psi'$,
the class function
\[
\sum_{\Lambda\in C_{\Phi,\Psi}}\sum_{\Lambda'\in C_{\Phi',\Psi'}}\rho_{\Lambda}\otimes\rho_{\Lambda'}
\]
on $G\times G'$ is uniform by Proposition~\ref{0601} and Proposition~\ref{0612}.
Then, by (\ref{0617}),
we have
\begin{align*}
\sum_{\Lambda\in C_{\Phi,\Psi}}\sum_{\Lambda'\in C_{\Phi',\Psi'}} x_{\Lambda,\Lambda'}
& =\biggl\langle\sum_{\Lambda\in C_{\Phi,\Psi}}\sum_{\Lambda'\in C_{\Phi',\Psi'}}
\rho_\Lambda\otimes\rho_{\Lambda'},\omega_{Z,Z'}\biggr\rangle \\
& =\biggl\langle\sum_{\Lambda\in C_{\Phi,\Psi}}\sum_{\Lambda'\in C_{\Phi',\Psi'}}\rho_\Lambda\otimes\rho_{\Lambda'},
\sum_{\Lambda''\in\cals^\epsilon_{Z'}}\rho_{\theta^\epsilon(\Lambda'')}\otimes\rho_{\Lambda''}\biggr\rangle.
\end{align*}
For a symbol $\Lambda''\in\cals^\epsilon_{Z'}$ to contribute a multiplicity, we need
both $\Lambda''=\Lambda'$ for some $\Lambda'\in C_{\Phi',\Psi'}$,
i.e., $\theta^\epsilon(\Lambda'')=\theta^\epsilon(\Lambda')$
for some $\theta^\epsilon(\Lambda')\in C_{\theta(\Phi'),\theta(\Psi')}$,
and $\theta^\epsilon(\Lambda'')=\Lambda$ for some $\Lambda\in C_{\Phi,\Psi}$, i.e.,
we need $\theta^\epsilon(\Lambda'')$ to be in the intersection $C_{\Phi,\Psi}\cap C_{\theta(\Phi'),\theta(\Psi')}$
by Lemma~\ref{1009}.
Therefore,
\begin{equation}\label{1106}
\sum_{\Lambda\in C_{\Phi,\Psi}}\sum_{\Lambda'\in C_{\Phi',\Psi'}} x_{\Lambda,\Lambda'}
=|C_{\Phi,\Psi}\cap C_{\theta(\Phi'),\theta(\Psi')}|
\end{equation}
for any arrangements $\Phi$ of $Z_\rmI$ with subset of pairs $\Psi$,
and any arrangement $\Phi'$ of $Z'_\rmI$ with admissible subset of pairs $\Psi'$.

Now let $\Phi=\theta(\Phi'')$ and $\Psi=\theta(\Psi'')$ for some arrangement $\Phi''$ of $Z'_\rmI$
with some subset of pairs $\Psi''$.
Then (\ref{1106}) becomes
\begin{equation}\label{1014}
\sum_{\Lambda''\in C_{\Phi'',\Psi''}}\sum_{\Lambda'\in C_{\Phi',\Psi'}} x_{\theta^\epsilon(\Lambda''),\Lambda'}
=|C_{\theta(\Phi''),\theta(\Psi'')}\cap C_{\theta(\Phi'),\theta(\Psi')}|
=|C_{\Phi'',\Psi''}\cap C_{\Phi',\Psi'}|
\end{equation}
for any two arrangements $\Phi',\Phi''$ of $Z'_\rmI$ and some admissible subsets of pairs $\Psi',\Psi''$ respectively.
Suppose that $\Lambda_1,\Lambda_2$ are symbols in $\cals^\epsilon_{Z'}$ such that $\Lambda_1\neq\Lambda_2,\Lambda_2^\rmt$.
Then by Lemma~\ref{0516}, there exists an arrangement $\Phi$ of $Z'_\rmI$ with two subsets of pairs $\Psi_1,\Psi_2$
such that $\Lambda_i,\Lambda_i^\rmt\in C_{\Phi,\Psi_i}$ for $i=1,2$ and $C_{\Phi,\Psi_1}\cap C_{\Phi,\Psi_2}=\emptyset$.
Because each $x_{\theta^\epsilon(\Lambda''),\Lambda'}$ is a non-negative integer,
from equation (\ref{1014}) we conclude that
$x_{\theta^\epsilon(\Lambda_1),\Lambda_2}=x_{\theta^\epsilon(\Lambda_1),\Lambda_2^\rmt}=0$.

Clearly there is an arrangement $\Phi'$ of $Z'_\rmI$ such that $\theta(\Phi')$ in (\ref{1012}) is
the arrangement $\Phi_2$ in Lemma~\ref{0602}.
Moreover, any subset of pairs of $\Phi_2$ is of the form $\theta(\Psi')$ for some subset of pairs $\Psi'$ of $\Phi'$.
For any given $\Lambda'\in\cals^\epsilon_{Z'}$,
then $\theta^\epsilon(\Lambda')\in\cals_Z$,
and by Lemma~\ref{0610}, there exist an arrangements $\Phi_1$ of $Z_\rmI$ with a subset of pairs $\Psi_1$
and a subset of pairs $\Psi'$ of $\Phi'$ given above such that $C_{\Phi_1,\Psi_1}\cap C_{\theta(\Phi'),\theta(\Psi')}=\{\theta^\epsilon(\Lambda')\}$.
Because we know $x_{\theta^\epsilon(\Lambda_1),\Lambda_2}=0$ for any $\Lambda_1\neq\Lambda_2,\Lambda_2^\rmt$,
equation (\ref{1106}) is reduced to
\[
x_{\theta^\epsilon(\Lambda'),\Lambda'}+x_{\theta^\epsilon(\Lambda'),\Lambda'^\rmt}=1.
\]

For $\Lambda'\in\cals^\epsilon_{Z'}$,
let $f'(\Lambda')$ be either $\Lambda'$ or $\Lambda'^\rmt$ such that
$\rho_{\theta^\epsilon(\Lambda')}\otimes\rho_{f'(\Lambda')}$ occurs in $\omega_{Z,Z'}$.
We just show that the character $\sum_{\Lambda'\in\cals_{Z'}^\epsilon}\rho_{\theta^\epsilon(\Lambda')}\otimes\rho_{f'(\Lambda')}$
is a sub-character of $\omega_{Z,Z'}$.
By the same argument in the last paragraph of proof of Proposition~\ref{1021},
we conclude that
\[
\omega_{Z,Z'}=\sum_{\Lambda'\in\cals_{Z'}^\epsilon}\rho_{\theta^\epsilon(\Lambda')}\otimes\rho_{f'(\Lambda')}.
\]
\end{proof}

\subsection{The general case}\label{1020}
Now let $Z,Z'$ be special symbols of defect $1,0$ respectively.
Suppose that $\cald_{Z,Z'}\neq\emptyset$.
Let $\Psi_0,\Psi'_0$ be the cores of $\cald_{Z,Z'}$ in $Z_\rmI,Z'_\rmI$ respectively, i.e.,
\begin{align*}
D_{Z'}& :=\{\,\Lambda\in\cals_{Z,1}\mid(\Lambda,Z')\in\cald_{Z,Z'}\,\}=\cals_{Z,\Psi_0}; \\
D_Z&:=\{\,\Lambda'\in\cals_{Z',0}\mid(Z,\Lambda')\in\cald_{Z,Z'}\,\}=\cals_{Z',\Psi'_0}
\end{align*}
by proposition~6.4 in \cite{pan-uniform}.
Let $\Phi,\Phi'$ be arrangements of $Z_\rmI,Z'_\rmI$ with subsets of pairs $\Psi,\Psi'$ respectively such that
$\Psi_0\leq\Psi\leq\Phi$ and $\Psi'_0\leq\Psi'\leq\Phi'$.
Then from the proofs of Proposition~7.16 in \cite{pan-uniform},
we know that either $\deg(Z'\smallsetminus\Psi'_0)=\deg(Z\smallsetminus\Psi_0)$
or $\deg(Z'\smallsetminus\Psi'_0)=\deg(Z\smallsetminus\Psi_0)+1$.
\begin{enumerate}
\item[(1)] Suppose that $\deg(Z'\smallsetminus\Psi'_0)=\deg(Z\smallsetminus\Psi_0)+1$.
Then we have a mapping $\theta^\epsilon\colon\cals_Z^{\Psi_0}\rightarrow\cals_{Z'}^{\epsilon,\Psi'_0}$.
Now $\Phi\smallsetminus\Psi_0$ is an arrangement of $Z_\rmI\smallsetminus\Psi_0$,
and $\theta(\Phi\smallsetminus\Psi_0)$ is an arrangement of $Z'_\rmI\smallsetminus\Psi'_0$
where $\theta$ is given in (\ref{0501}).
Now $\Psi\smallsetminus\Psi_0$ is a subset of pairs of $\Phi\smallsetminus\Psi_0$,
and we define
\[
\bar\theta(\Phi)=\theta(\Phi\smallsetminus\Psi_0)\cup\Psi'_0,\qquad
\bar\theta(\Psi)=\theta(\Psi\smallsetminus\Psi_0)\cup\Psi'_0.
\]
It is easy to see that $\bar\theta(\Phi)$ is an arrangement of $Z'_\rmI$ and $\bar\theta(\Psi)$ is a
subset of pairs of $\bar\theta(\Phi)$.

\item[(2)] Suppose that $\deg(Z'\smallsetminus\Psi'_0)=\deg(Z\smallsetminus\Psi_0)$.
Then we have a mapping $\theta^\epsilon\colon\cals_{Z'}^{\epsilon,\Psi'_0}\rightarrow\cals_Z^{\Psi_0}$.
Then we define
\begin{equation}\label{0503}
\bar\theta(\Phi')=\theta(\Phi'\smallsetminus\Psi'_0)\cup\Psi_0,\qquad
\bar\theta(\Psi')=\theta(\Psi'\smallsetminus\Psi'_0)\cup\Psi_0
\end{equation}
where $\theta$ is given in (\ref{0502}).
Similarly, $\bar\theta(\Phi')$ is an arrangement of $Z_\rmI$ and $\bar\theta(\Psi')$ is a
subset of pairs of $\theta(\Phi')$.
\end{enumerate}

\begin{lemma}\label{1017}
Keep the above setting.
\begin{enumerate}
\item[(i)] Suppose that $\deg(Z'\smallsetminus\Psi'_0)=\deg(Z\smallsetminus\Psi_0)+1$.
Then
\[
C_{\bar\theta(\Phi),\bar\theta(\Psi)}=
\{\,\theta^\epsilon(\Lambda_1)+\Lambda_2,\theta^\epsilon(\Lambda_1)^\rmt+\Lambda_2
\mid\Lambda_1\in C_{\Phi,\Psi}\cap\cals_Z^{\Psi_0},\ \Lambda_2\in\cals_{Z',\Psi'_0}\,\}.
\]

\item[(ii)] Suppose that $\deg(Z'\smallsetminus\Psi'_0)=\deg(Z\smallsetminus\Psi_0)$.
Let $\Phi'$ be an arrangement of $Z'_\rmI$ and let $\Psi'$ be a subset of pairs in $\Phi'$ such that
$|\Phi'\smallsetminus\Psi'|$ is even if $\epsilon=+$, and $|\Phi'\smallsetminus\Psi'|$ is odd if $\epsilon=-$.
Then
\[
C_{\bar\theta(\Phi'),\bar\theta(\Psi')}
=\{\,\theta^\epsilon(\Lambda_1)+\Lambda_2\mid\Lambda_1\in C_{\Phi',\Psi'}\cap\cals_{Z'}^{\epsilon,\Psi'_0},\ \Lambda_2\in\cals_{Z,\Psi_0}\,\}.
\]
\end{enumerate}
\end{lemma}
\begin{proof}
First suppose that $\deg(Z'\smallsetminus\Psi'_0)=\deg(Z\smallsetminus\Psi_0)+1$.
Now $\Psi_0'\leq\theta(\Psi)\leq\theta(\Phi)$,
so by Lemma~\ref{0515}, we have
\[
C_{\bar\theta(\Phi),\bar\theta(\Psi)}=
\{\,\Lambda_1'+\Lambda_2\mid\Lambda'_1\in C_{\theta(\Phi),\theta(\Psi)}\cap\cals_{Z'}^{\epsilon,\Psi'_0},\ \Lambda_2\in\cals_{Z',\Psi'_0}\,\}.
\]
We know that the relation $\cald_{Z,Z'}\cap(\cals_Z^{\Psi_0}\times\cals_{Z'}^{\epsilon,\Psi_0'})$ is one-to-one,
so by Lemma~\ref{1008}, we see that
\[
C_{\bar\theta(\Phi),\bar\theta(\Psi)}\cap\cals_{Z'}^{\epsilon,\Psi'_0}
=\{\,\theta^\epsilon(\Lambda_1),\theta^\epsilon(\Lambda_1)^\rmt\mid
\Lambda_1\in C_{\Phi,\Psi}\cap\cals_Z^{\Psi_0}\,\}
\]
and hence the lemma is proved.

The proof of (ii) is similar.
\end{proof}

\begin{theorem}\label{0522}
Let $(\bfG,\bfG')=(\Sp_{2n},\rmO^\epsilon_{2n'})$,
and let $Z,Z'$ be special symbols of rank $n,n'$ and defect $1,0$ respectively.
Then $\rho_\Lambda\otimes\rho_{\Lambda'}$ or $\rho_\Lambda\otimes\rho_{\Lambda'^\rmt}$
occurs in $\omega_{Z,Z'}$ if and only if $(\Lambda,\Lambda')$ or $(\Lambda,\Lambda'^\rmt)$
occurs in $\calb^\epsilon_{Z,Z'}$.
\end{theorem}
\begin{proof}
Let $Z,Z'$ be special symbols of defects $1,0$ respectively.
Let $\Psi_0,\Psi'_0$ be the cores in $Z_\rmI,Z'_\rmI$ of $\cald_{Z,Z'}$,
and $\delta_0,\delta_0'$ the degrees of $\Psi_0,\Psi'_0$ respectively.
From \cite{pan-uniform} (6.4) and (8.1), we have
\begin{align*}
\sum_{(\Lambda,\Lambda')\in\calb^\epsilon_{Z,Z'}}\rho_\Lambda\otimes\rho_{\Lambda'}
&=\sum_{(\Lambda_1,\Lambda'_1)\in\calb^{\epsilon,\natural}_{Z,Z'}}
\sum_{\Lambda_2\in D_{Z'}}\sum_{\Lambda'_2\in D_Z} \rho_{\Lambda_1+\Lambda_2}\otimes\rho_{\Lambda'_1+\Lambda'_2}.
\end{align*}
Note that $D_{Z'}=\{\,\Lambda_M\mid M\leq\Psi_0\,\}$ and $D_Z=\{\,\Lambda_N\mid N\leq\Psi'_0\,\}$ by
proposition~6.4 in \cite{pan-uniform}.

If $\Phi$ is an arrangement of $Z_\rmI$ and $\Psi$ is a subset of pairs in $\Phi$ such that $\Psi_0\leq\Psi\leq\Phi$,
then by Lemma~\ref{0515} and Proposition~\ref{0601},
\[
\sum_{\Lambda\in C_{\Phi,\Psi}}\rho_\Lambda
=\sum_{\Lambda_1\in C_{\Phi,\Psi}^\natural}\sum_{\Lambda_2\in\cals_{Z,\Psi_0}}\rho_{\Lambda_1+\Lambda_2}
\]
is a uniform class function on $G$ where $C_{\Phi,\Psi}^\natural=C_{\Phi,\Psi}\cap\cals_Z^{\Psi_0}$.
Similarly,
\[
\sum_{\Lambda'\in C_{\Phi',\Psi'}}\rho_{\Lambda'}+\rho_{\Lambda'^\rmt}
=\sum_{\Lambda'_1\in C^\natural_{\Phi',\Psi'}}\sum_{\Lambda'_2\in\cals_{Z',\Psi_0'}}\rho_{\Lambda'_1+\Lambda'_2}+\rho_{\Lambda_1'^\rmt+\Lambda'_2}
\]
is a uniform class function on $G'$ where $\Phi'$ is an arrangement of $Z'_\rmI$,
$\Psi'$ an admissible subset of pairs in $\Phi'$ such that $\Psi'_0\leq\Psi'\leq\Phi'$,
and $C^\natural_{\Phi',\Psi'}=C_{\Phi',\Psi'}\cap\cals_{Z'}^{\epsilon,\Psi'_0}$.

\begin{enumerate}
\item[(1)] Suppose that $\deg(Z'\smallsetminus\Psi'_0)=\deg(Z\smallsetminus\Psi_0)+1$.
Then from the result in \cite{pan-uniform} subsection~7.4, we see that
there exists a mapping $\theta^\epsilon\colon\cals_Z^{\Psi_0}\rightarrow\cals^{\epsilon,\Psi'_0}_{Z'}$ and
\begin{equation}\label{1015}
\calb^{\epsilon,\natural}_{Z,Z'}
=\{\,(\Lambda,\theta^\epsilon(\Lambda))\mid\Lambda\in\cals_Z^{\Psi_0}\,\}.
\end{equation}
For an arrangement $\Phi'$ of $Z_\rmI$ with a subset of pairs $\Psi'$ such that $\Psi_0\leq\Psi'$,
by (1) of Lemma~\ref{1017},
the class function
\[
\sum_{\Lambda'\in C_{\bar\theta(\Phi'),\bar\theta(\Psi')}}\rho_{\Lambda'}
=\sum_{\Lambda'_1\in C^\natural_{\Phi',\Psi'}}\sum_{\Lambda'_2\in D_Z}
\rho_{\theta^\epsilon(\Lambda_1')+\Lambda'_2}+\rho_{\theta^\epsilon(\Lambda_1')^\rmt+\Lambda'_2}
\]
on $G'$ is uniform.
For $\Lambda_1,\Lambda_1'\in\cals_Z^{\Psi_0}$,
define
\[
x_{\Lambda_1,\Lambda_1'}
=\frac{1}{2^{\delta_0+\delta'_0}}\sum_{\Lambda_2\in D_{Z'}}\sum_{\Lambda'_2\in D_Z}
\langle\rho_{\Lambda_1+\Lambda_2}\otimes\rho_{\theta^\epsilon(\Lambda_1')+\Lambda'_2}
+\rho_{\Lambda_1+\Lambda_2}\otimes\rho_{\theta^\epsilon(\Lambda_1')^\rmt+\Lambda'_2},
\omega_{Z,Z'}\rangle.
\]
Note that $|D_{Z'}|=2^{\delta_0}$ and $|D_Z|=2^{\delta_0'}$.
Now by (\ref{1022}) and (\ref{1015}), we have
\begin{align*}
& \sum_{\Lambda_1\in C^\natural_{\Phi,\Psi}}\sum_{\Lambda'_1\in C^\natural_{\Phi',\Psi'}} x_{\Lambda_1,\Lambda'_1}\\
& = \frac{1}{2^{\delta_0+\delta'_0}}\biggl\langle \sum_{\Lambda_1\in C^\natural_{\Phi,\Psi}}\sum_{\Lambda'_1\in C^\natural_{\Phi',\Psi'}}
\sum_{\Lambda_2\in D_{Z'}}\sum_{\Lambda'_2\in D_Z}
(\rho_{\Lambda_1+\Lambda_2}\otimes\rho_{\theta^\epsilon(\Lambda_1')+\Lambda'_2}
+\rho_{\Lambda_1+\Lambda_2}\otimes\rho_{\theta^\epsilon(\Lambda_1')^\rmt+\Lambda'_2}), \\
&\qquad\qquad \sum_{\Lambda''_1\in\cals_Z^{\Psi_0}}
\sum_{\Lambda''_2\in D_{Z'}}\sum_{\Lambda'''_2\in D_Z} \rho_{\Lambda''_1+\Lambda''_2}\otimes\rho_{\theta^\epsilon(\Lambda''_1)+\Lambda'''_2}\biggr\rangle
\end{align*}
where $\Phi,\Phi'$ are arrangements of $Z_\rmI$ with subsets of pairs $\Psi,\Psi'$ respectively
such that $\Psi_0\leq\Psi$ and $\Psi_0\leq\Psi'$.
For a symbol $\Lambda''_1\in\cals_Z^{\Psi_0}$ to contribute a multiplicity in the above identity,
by Lemma~\ref{0521},
we need $\Lambda_1''=\Lambda_1$ and $\Lambda_1''=\Lambda_1'$ for some $\Lambda_1\in C^\natural_{\Phi,\Psi}$ and some
$\Lambda'_1\in C^\natural_{\Phi',\Psi'}$,
i.e., $\Lambda''_1$ must be in the intersection $C^\natural_{\Phi,\Psi}\cap C^\natural_{\Phi',\Psi'}$.
Then
\begin{align}\label{1016}
\begin{split}
\sum_{\Lambda_1\in C^\natural_{\Phi,\Psi}}\sum_{\Lambda'_1\in C^\natural_{\Phi',\Psi'}} x_{\Lambda_1,\Lambda'_1}
&= \frac{1}{2^{\delta_0+\delta_0'}}\cdot|C^\natural_{\Phi,\Psi}\cap C^\natural_{\Phi',\Psi'}|\cdot\sum_{\Lambda_2\in D_{Z'}}\sum_{\Lambda_2'\in D_Z}1 \\
&= |C^\natural_{\Phi,\Psi}\cap C^\natural_{\Phi',\Psi'}|.
\end{split}
\end{align}
for any arrangements $\Phi,\Phi'$ of $Z_\rmI$ with subsets of pairs $\Psi,\Psi'$ respectively such that
$\Psi_0\leq\Psi$ and $\Psi_0\leq\Psi'$.

Suppose that $\Lambda_1,\Lambda_2$ are distinct symbols in $\cals_Z^{\Psi_0}$.
By Lemma~\ref{0520}, there exist an arrangement $\Phi$ of $Z_\rmI$ with two subsets of pairs $\Psi_1,\Psi_2$
such that $\Psi_0\leq\Psi_i$, $\Lambda_i\in C_{\Phi,\Psi_i}$ for $i=1,2$ and $C_{\Phi,\Psi_1}\cap C_{\Phi,\Psi_2}=\emptyset$.
Then $C_{\Phi,\Psi_1}^\natural\cap C_{\Phi,\Psi_2}^\natural=\emptyset$.
Because each $x_{\Lambda,\Lambda'}$ is a non-negative integer,
from (\ref{1016}) we conclude that $x_{\Lambda_1,\Lambda_2}=0$ for any distinct $\Lambda_1,\Lambda_2\in\cals_Z^{\Psi_0}$.

Finally, for any $\Lambda\in\cals_Z^{\Psi_0}$, by Lemma~\ref{0519},
there exist two arrangements $\Phi_1,\Phi_2$ of $Z_\rmI$ with subsets of pairs $\Psi_1,\Psi_2$
respectively such that $\Psi_0\leq\Psi_i$ for $i=1,2$ and
\[
C_{\Phi_1,\Psi_1}^\natural\cap C_{\Phi_2,\Psi_2}^\natural
=C_{\Phi_1,\Psi_1}\cap C_{\Phi_2,\Psi_2}\cap\cals_Z^{\Psi_0}=\{\Lambda\}.
\]
Because we know $x_{\Lambda_1,\Lambda_2}=0$ if $\Lambda_1\neq\Lambda_2$,
equation (\ref{1016}) is reduced to $x_{\Lambda,\Lambda}=1$.

Suppose that $\Lambda\in\cals_Z$ and $\Lambda'\in\cals^\epsilon_{Z'}$ such that
$(\Lambda,\Lambda')\in\calb^\epsilon_{Z,Z'}$.
We can write $\Lambda=\Lambda_1+\Lambda_2$ for $\Lambda_1\in\cals_Z^{\Psi_0}$ and $\Lambda_2\in\cals_{Z,\Psi_0}$
and similarly write $\Lambda'=\Lambda'_1+\Lambda'_2$ for $\Lambda'_1\in\cals_{Z'}^{\epsilon,\Psi_0}$ and $\Lambda'_2\in\cals_{Z',\Psi'_0}$
such that $(\Lambda_1,\Lambda'_1)\in\calb_{Z,Z'}^{\epsilon,\natural}$.
Then we have shown that either
$\rho_{\Lambda_1+\Lambda_2}\otimes\rho_{\Lambda'_1+\Lambda'_2}$ or
$\rho_{\Lambda_1+\Lambda_2}\otimes\rho_{\Lambda'^\rmt_1+\Lambda'_2}$ occurs in $\omega_{Z,Z'}$.
It is not difficult to check that $\Lambda'^\rmt=(\Lambda_1'+\Lambda'_2)^\rmt=\Lambda_1'^\rmt+\Lambda'_2$
(\cf.~\cite{pan-uniform} lemma~2.1).
Hence we conclude that either $\rho_\Lambda\otimes\rho_{\Lambda'}$ or
$\rho_\Lambda\otimes\rho_{\Lambda'^\rmt}$ occurs in $\omega_{Z,Z'}$.

\item[(2)] Suppose that $\deg(Z'\smallsetminus\Psi'_0)=\deg(Z\smallsetminus\Psi_0)$.
From the result in \cite{pan-uniform} subsection~7.4, we see that
there exists a mapping $\theta^\epsilon\colon\cals_{Z'}^{\epsilon,\Psi'_0}\rightarrow\cals_Z^{\Psi_0}$ and
\[
\calb^{\epsilon,\natural}_{Z,Z'}
=\{\,(\theta^\epsilon(\Lambda'),\Lambda')\mid\Lambda'\in\cals_{Z'}^{\epsilon,\Psi'_0}\,\}.
\]
For an arrangement $\Phi$ of $Z_\rmI$ with a subset of pairs $\Psi$ such that $\Psi_0\leq\Psi$,
and an arrangement $\Phi'$ of $Z'_\rmI$ with an admissible subset of pairs $\Psi'$ such that $\Psi_0'\leq\Psi'$,
the class function
\begin{align*}
\sum_{\Lambda\in C_{\Phi,\Psi}}\sum_{\Lambda'\in C_{\Phi',\Psi'}}\rho_{\Lambda}\otimes\rho_{\Lambda'}
=\sum_{\Lambda_1\in C^\natural_{\Phi,\Psi}}\sum_{\Lambda'_1\in C^\natural_{\Phi',\Psi'}}
\sum_{\Lambda_2\in D_{Z'}}\sum_{\Lambda'_2\in D_Z}\rho_{\Lambda_1+\Lambda_2}\otimes\rho_{\Lambda_1'+\Lambda'_2}
\end{align*}
on $G\times G'$ is uniform by Proposition~\ref{0601} and Proposition~\ref{0612}.

For $\Lambda_1\in\cals_Z^{\Psi_0}$ and $\Lambda_1'\in\cals_{Z'}^{\epsilon,\Psi'_0}$,
define
\[
x_{\Lambda_1,\Lambda_1'}
=\frac{1}{2^{\delta_0+\delta'_0}}\sum_{\Lambda_2\in D_{Z'}}\sum_{\Lambda'_2\in D_Z}
\langle\rho_{\Lambda_1+\Lambda_2}\otimes\rho_{\Lambda_1'+\Lambda'_2},
\omega_{Z,Z'}\rangle.
\]
Then
\begin{align*}
\sum_{\Lambda_1\in C^\natural_{\Phi,\Psi}}\sum_{\Lambda'_1\in C^\natural_{\Phi',\Psi'}} x_{\Lambda_1,\Lambda'_1}
&= \frac{1}{2^{\delta_0+\delta'_0}}\biggl\langle \sum_{\Lambda_1\in C^\natural_{\Phi,\Psi}}\sum_{\Lambda'_1\in C^\natural_{\Phi',\Psi'}}
\sum_{\Lambda_2\in D_{Z'}}\sum_{\Lambda'_2\in D_Z}\rho_{\Lambda_1+\Lambda_2}\otimes\rho_{\Lambda_1'+\Lambda'_2}, \\
&\qquad\qquad \sum_{\Lambda''_1\in\cals_{Z'}^{\epsilon,\Psi_0'}}
\sum_{\Lambda''_2\in D_{Z'}}\sum_{\Lambda'''_2\in D_Z} \rho_{\theta^\epsilon(\Lambda''_1)+\Lambda''_2}\otimes\rho_{\Lambda''_1+\Lambda'''_2}
\biggr\rangle.
\end{align*}
By the same argument in (1) we conclude that
\begin{equation}\label{1018}
\sum_{\Lambda_1\in C^\natural_{\Phi,\Psi}}\sum_{\Lambda'_1\in C^\natural_{\Phi',\Psi'}} x_{\Lambda_1,\Lambda'_1}
= |C^\natural_{\Phi,\Psi}\cap C^\natural_{\bar\theta(\Phi'),\bar\theta(\Psi')}|
\end{equation}
for any arrangement $\Phi$ of $Z_\rmI$ with subset of pairs $\Psi$ such that $\Psi_0\leq\Psi$,
and any arrangement $\Phi'$ of $Z'_\rmI$ with admissible subset of pairs $\Psi'$ such that $\Psi'_0\leq\Psi'$.

Now let $\Phi=\bar\theta(\Phi'')$ and $\Psi=\bar\theta(\Psi'')$ for some arrangement $\Phi''$ of $Z'_\rmI$
with some admissible subset of pairs $\Psi''$.
Then (\ref{1018}) becomes
\begin{equation}\label{1019}
\sum_{\Lambda''_1\in C^\natural_{\Phi'',\Psi''}}\sum_{\Lambda'_1\in C^\natural_{\Phi',\Psi'}} x_{\theta^\epsilon(\Lambda''_1),\Lambda'_1}
=|C^\natural_{\bar\theta(\Phi''),\bar\theta(\Psi'')}\cap C^\natural_{\bar\theta(\Phi'),\bar\theta(\Psi')}|
=|C^\natural_{\Phi'',\Psi''}\cap C^\natural_{\Phi',\Psi'}|
\end{equation}
for any two arrangements $\Phi',\Phi''$ of $Z'_\rmI$ and some admissible subsets of pairs $\Psi',\Psi''$ respectively.
Suppose $\Lambda_1,\Lambda_2$ are symbols in $\cals^{\epsilon,\Psi_0'}_{Z'}$ such that $\Lambda_1\neq\Lambda_2,\Lambda_2^\rmt$.
Then by Lemma~\ref{0518}, there exist an arrangement $\Phi$ of $Z'_\rmI$ with two subsets of pairs $\Psi_1,\Psi_2$
such that $\Lambda_i,\Lambda_i^\rmt\in C_{\Phi,\Psi_i}$ for $i=1,2$ and $C_{\Phi,\Psi_1}\cap C_{\Phi,\Psi_2}=\emptyset$.
Because each $x_{\theta^\epsilon(\Lambda''),\Lambda'}$ is a non-negative integer,
from equation (\ref{1019}) we conclude that
\[
x_{\theta^\epsilon(\Lambda_1),\Lambda_2}=x_{\theta^\epsilon(\Lambda_1),\Lambda_2^\rmt}=0.
\]

Clearly there is an arrangement $\Phi'$ of $Z'_\rmI$ such that $\bar\theta(\Phi')$ in (\ref{0503}) is
the arrangement $\Phi_2$ in Lemma~\ref{0519}.
Moreover, any subset of pairs of $\Phi_2$ is of the form $\bar\theta(\Psi')$
for some admissible subset of pairs $\Psi'$ of $\Phi'$.
For any given $\Lambda'_1\in\cals^{\epsilon,\Psi_0'}_{Z'}$,
then $\theta^\epsilon(\Lambda'_1)\in\cals_Z^{\Psi_0}$,
and by Lemma~\ref{0610}, there exist an arrangements $\Phi_1$ of $Z_\rmI$ with a subset of pairs $\Psi_1$
and an admissible subset of pairs $\Psi'$ of $\Phi'$ given above such that
$\Psi_0\leq\Psi_1$, $\Psi_0'\leq\Psi'$,
$C^\natural_{\Phi_1,\Psi_1}\cap C^\natural_{\theta(\Phi'),\theta(\Psi')}=\{\theta^\epsilon(\Lambda'_1)\}$.
Because we know $x_{\theta^\epsilon(\Lambda_1),\Lambda_2}=0$ for any $\Lambda_1\neq\Lambda_2,\Lambda_2^\rmt$,
equation (\ref{1018}) is reduced to
\[
x_{\theta^\epsilon(\Lambda'_1),\Lambda'_1}+x_{\theta^\epsilon(\Lambda'_1),\Lambda'^\rmt_1}=1.
\]
By the same argument in the last paragraph of (1),
we conclude that if $(\Lambda,\Lambda')\in\calb^\epsilon_{Z,Z'}$,
then either $\rho_\Lambda\otimes\rho_{\Lambda'}$ or
$\rho_\Lambda\otimes\rho_{\Lambda'^\rmt}$ occurs in $\omega_{Z,Z'}$.
\end{enumerate}

For two cases, and for each $(\Lambda,\Lambda')\in\calb^\epsilon_{Z,Z'}$, define $\widetilde\Lambda'$ to be either $\Lambda'$
or $\Lambda'^\rmt$ such that $\rho_\Lambda\otimes\rho_{\widetilde\Lambda'}$ occurs in $\omega_{Z,Z'}$.
Therefore, $\sum_{(\Lambda,\Lambda')\in\calb^\epsilon_{Z,Z'}}\rho_\Lambda\otimes\rho_{\widetilde\Lambda'}$ is a sub-character of
$\omega_{Z,Z'}$.
Then by the same argument in the last paragraph of the proof of Proposition~\ref{1021},
we conclude that
\[
\omega_{Z,Z'}=\sum_{(\Lambda,\Lambda')\in\calb^\epsilon_{Z,Z'}}\rho_\Lambda\otimes\rho_{\widetilde\Lambda'}.
\]
\end{proof}

From the above proof, we conclude the following corollary:

\begin{corollary}\label{1024}
Let $(\bfG,\bfG')=(\Sp_{2n},\rmO_{2n'}^\epsilon)$,
$\rho_\Lambda\in\cale(G)_1$, $\rho_{\Lambda'}\in\cale(G')_1$.
\begin{enumerate}
\item[(i)] Both $(\Lambda,\Lambda')$ and $(\Lambda,\Lambda'^\rmt)$ occur in $\calb_{\bfG,\bfG'}$
if and only if both $\rho_\Lambda\otimes\rho_{\Lambda'}$ and $\rho_\Lambda\otimes\rho_{\Lambda'^\rmt}$
occur in the correspondence.

\item[(ii)] Exactly one of $(\Lambda,\Lambda'),(\Lambda,\Lambda'^\rmt)$ occurs in $\calb_{\bfG,\bfG'}$
if and only if exactly one of $\rho_\Lambda\otimes\rho_{\Lambda'},\rho_\Lambda\otimes\rho_{\Lambda'^\rmt}$
occurs in the correspondence.
\end{enumerate}
\end{corollary} 

\section{Symbol Correspondence and Parabolic Induction}\label{0618}

\subsection{Properties of the parametrization}\label{1006}
For a symbol
\begin{equation}\label{1115}
\Lambda=\binom{a_1,a_2,\ldots,a_{m_1}}{b_1,b_2,\ldots,b_{m_2}},
\end{equation}
we define $\Omega^+_\Lambda$ to be the set consisting of the following types of symbols:
\begin{enumerate}
\item[(I)] $\displaystyle\binom{a_1,\ldots,a_{i-1},a_i+1,a_{i+1},\ldots,a_{m_1}}{b_1,b_2,\ldots,b_{m_2}}$ for $i=1,\ldots,m_1$
such that $a_{i-1}>a_i+1$ (the condition is empty if $i=1$);

\item[(II)] $\displaystyle\binom{a_1,a_2,\ldots,a_{m_1}}{b_1,\ldots,b_{j-1},b_j+1,b_{j+1},\ldots,b_{m_2}}$ for $j=1,\ldots,m_2$
such that $b_{j-1}>b_j+1$ (the condition is empty if $j=1$);

\item[(III)] $\displaystyle\binom{a_1+1,a_2+1,\ldots,a_{m_1}+1,1}{b_1+1,b_2+1,\ldots,b_{m_2}+1,0}$ if $a_{m_1}\neq 0$;

\item[(IV)] $\displaystyle\binom{a_1+1,a_2+1,\ldots,a_{m_1}+1,0}{b_1+1,b_2+1,\ldots,b_{m_2}+1,1}$ if $b_{m_2}\neq 0$.
\end{enumerate}

Clearly, if $\Lambda'\in\Omega_\Lambda^+$,
then ${\rm rank}(\Lambda')={\rm rank}(\Lambda)+1$ and ${\rm def}(\Lambda')={\rm def}(\Lambda)$.
We also define
\[
\Omega^-_\Lambda=\{\,\Lambda'\mid\Lambda\in\Omega^+_{\Lambda'}\,\}.
\]
It is easy to see that $\Omega^-_\Lambda$ consists of symbols of the following types:
\begin{enumerate}
\item[(I')] $\displaystyle\binom{a_1,\ldots,a_{i-1},a_i-1,a_{i+1},\ldots,a_{m_1}}{b_1,b_2,\ldots,b_m}$ for $i=1,\ldots,m_1-1$
such that $a_i>a_{i+1}+1$;

\noindent $\displaystyle\binom{a_1,\ldots,a_{m_1-1},a_{m_1}-1}{b_1,b_2,\ldots,b_m}$ if $a_{m_1}\geq 1$ and $(a_{m_1},b_{m_2})\neq (1,0)$;

\item[(II')] $\displaystyle\binom{a_1,a_2,\ldots,a_{m_1}}{b_1,\ldots,b_{j-1},b_j-1,b_{j+1},\ldots,b_{m_2}}$ for $j=1,\ldots,m_2-1$ such that $b_j>b_{j+1}+1$;

\noindent $\displaystyle\binom{a_1,a_2,\ldots,a_{m_1}}{b_1,\ldots,b_{m_2-1},b_{m_2}-1}$ if $b_{m_2}\geq 1$ and $(a_{m_1},b_{m_2})\neq(0,1)$;

\item[(III')] $\displaystyle\binom{a_1-1,a_2-1,\ldots,a_{m_1-1}-1}{b_1-1,b_2-1,\ldots,b_{m_2-1}-1}$ if $(a_{m_1},b_{m_2})=(1,0),(0,1)$.
\end{enumerate}
Therefore, if $\Lambda'\in\Omega_\Lambda^-$,
then ${\rm rank}(\Lambda')={\rm rank}(\Lambda)-1$ and ${\rm def}(\Lambda')={\rm def}(\Lambda)$.

\begin{example}
Suppose that $\Lambda=\binom{4,2,1}{3,0}$.
Then
\begin{align*}
\Omega_\Lambda^+ &=\textstyle\{\binom{5,2,1}{3,0},\binom{4,3,1}{3,0},\binom{4,2,1}{4,0},\binom{4,2,1}{3,1},\binom{5,3,2,1}{4,1,0}\},\\
\Omega_\Lambda^- &=\textstyle\{\binom{3,2,1}{3,0},\binom{4,2,1}{2,0},\binom{3,1}{2}\}.
\end{align*}
\end{example}

Recall that in Subsection~\ref{0335} and Subsection~\ref{0331},
for each irreducible unipotent character $\rho$ of $\Sp_{2n}(q)$ or $\rmO^\epsilon_{2n'}(q)$
we can associate it a unique symbol $\Lambda$ and denote $\rho=\rho_\Lambda$.
The parametrization also satisfies the following conditions:
\begin{enumerate}
\item[(1)] Let $\zeta_m$ denote the unique unipotent cuspidal character of $\Sp_{2m(m+1)}(q)$.
Then by our convention in (\ref{0320}),
we have $\zeta_m=\rho_{\Lambda_m}$ where
\[
\Lambda_m=\begin{cases}
\binom{2m,2m-1,\ldots,1,0}{-}, & \text{if $m$ is even};\\
\binom{-}{2m,2m-1,\ldots,1,0}, & \text{if $m$ is odd}.
\end{cases}
\]
Let $\zeta_m^\rmI,\zeta_m^{\rm II}$ be the unipotent cuspidal characters of $\rmO^\epsilon_{2m^2}(q)$ (for $m\geq 1$)
where $\epsilon$ is the sign of $(-1)^m$
such that $\zeta_m\otimes\zeta_m^{\rm II}$ and $\zeta_m\otimes\zeta_{m+1}^\rmI$ occur in the Howe
correspondence (\cf.~\cite{adams-moy}).
Then, by our convention in Subsection~\ref{0331},
we have $\zeta_m^\rmI=\rho_{\Lambda'_m}$ and $\zeta_m^{\rm II}=\rho_{(\Lambda'_m)^\rmt}$
where
\[
\Lambda'_m=\begin{cases}
\binom{2m-1,2m-2,\ldots,1,0}{-}, & \text{if $m$ is even}; \\
\binom{-}{2m-1,2m-2,\ldots,1,0}, & \text{if $m$ is odd}.
\end{cases}
\]

\item[(2)] ${\bf 1}_{\rmO^+_2}=\rho_{\binom{1}{0}}$ and hence $\sgn_{\rmO^+_2}=\rho_{\binom{0}{1}}$.

\item[(3)] The following branching rule holds:
\begin{align}\label{1201}
\begin{split}
\Ind^{\Sp_{2(n+1)}(q)}_{\Sp_{2n}(q)\times\GL_1(q)}(\rho_\Lambda\otimes{\bf 1})
&=\sum_{\Lambda''\in\Omega_\Lambda^+}\rho_{\Lambda''}; \\
\Ind^{\rmO^+_{2(n'+1)}(q)}_{\rmO_{2n'}^+(q)\times\GL_1(q)}(\rho_{\Lambda'}\otimes{\bf 1})
&=\sum_{\Lambda'''\in\Omega_{\Lambda'}^+}\rho_{\Lambda'''}.
\end{split}
\end{align}
In particular, the defects are preserved by parabolic induction.
\end{enumerate}

\begin{proposition}\label{1107}
Let $(\bfG,\bfG')=(\Sp_{2n},\rmO_{2n'}^\epsilon)$,
$\rho_\Lambda\in\cale(G)_1$, $\rho_{\Lambda'}\in\cale(G')_1$.
Suppose that ${\rm def}(\Lambda')\neq 0$.
Then $\rho_\Lambda\otimes\rho_{\Lambda'}$ occurs in $\omega_{\bfG,\bfG',1}$ if and only if
$(\Lambda,\Lambda')\in\calb_{\bfG,\bfG'}$.
\end{proposition}
\begin{proof}
If $(\Lambda,\Lambda')\in\calb_{\bfG,\bfG'}$,
then from the definition of $\calb^\epsilon_{Z,Z'}$ in (\ref{0223}) and (\ref{0227}) we know that
\begin{equation}\label{1108}
{\rm def}(\Lambda')=\begin{cases}
-{\rm def}(\Lambda)+1, & \text{if $\epsilon=+$};\\
-{\rm def}(\Lambda)-1, & \text{if $\epsilon=-$}.
\end{cases}
\end{equation}
We also know that ${\rm def}(\Lambda'^\rmt)=-{\rm def}(\Lambda')\neq{\rm def}(\Lambda')$
since we assume that ${\rm def}(\Lambda')\neq 0$.
So we only need to show that the defects of $\Lambda,\Lambda'$ satisfies (\ref{1108})
if $\rho_\Lambda\otimes\rho_{\Lambda'}$ occurs in $\omega_{\bfG,\bfG',1}$.
Now suppose that $\rho_\Lambda\otimes\rho_{\Lambda'}$ occurs in $\omega_{\bfG,\bfG',1}$.
\begin{enumerate}
\item[(1)] First suppose that both $\rho_\Lambda,\rho_{\Lambda'}$ are cuspidal.
Then we know that $n=m(m+1)$ and $n'=m^2,(m+1)^2$ for some $m$.

\begin{enumerate}
\item Suppose that $n=m(m+1)$ and $n'=(m+1)^2$ (and $\epsilon$ is the sign of $(-1)^{m+1}$).
Then we know that $\zeta_m\otimes\zeta_{m+1}^\rmI$ occurs in $\omega_{\bfG,\bfG',1}$ from \cite{adams-moy} theorem 5.2.
Now $\zeta_m=\rho_{\Lambda_m}$, $\zeta^\rmI_{m+1}=\rho_{\Lambda'_{m+1}}$ from our convention above
and we have
\[
\Lambda'_{m+1}=\begin{cases}
\binom{-}{2m+1}\cup(\Lambda_m)^\rmt, & \text{if $m$ is even};\\
\binom{2m+1}{-}\cup(\Lambda_m)^\rmt, & \text{if $m$ is odd}.
\end{cases}
\]
Hence ${\rm def}(\Lambda'_{m+1})=-{\rm def}(\Lambda_m)+1$ if $\epsilon=+$;
and ${\rm def}(\Lambda'_{m+1})=-{\rm def}(\Lambda_m)-1$ if $\epsilon=-$.

\item Suppose that $n=m(m+1)$ and $n'=m^2$ (and $\epsilon$ is the sign of $(-1)^m$).
Then we know that $\zeta_m\otimes\zeta_m^{\rm II}$ occurs in $\omega_{\bfG,\bfG',1}$ from \cite{adams-moy} theorem 5.2.
Now $\zeta_m=\rho_{\Lambda_m}$, $\zeta^{\rm II}_m=\rho_{(\Lambda'_m)^\rmt}$ from our convention above
and
\[
(\Lambda'_m)^\rmt=\begin{cases}
(\Lambda_m)^\rmt\smallsetminus\binom{-}{2m}, & \text{if $m$ is even (hence $\epsilon=+$)};\\
(\Lambda_m)^\rmt\smallsetminus\binom{2m}{-}, & \text{if $m$ is odd (hence $\epsilon=-$)}.
\end{cases}
\]
Hence ${\rm def}((\Lambda'_m)^\rmt)=-{\rm def}(\Lambda_m)+1$ if $\epsilon=+$;
and ${\rm def}((\Lambda'_m)^\rmt)=-{\rm def}(\Lambda_m)-1$ if $\epsilon=-$.
\end{enumerate}
Therefore, the corollary is proved if both $\rho_\Lambda,\rho_{\Lambda'}$ are cuspidal.

\item[(2)] Next, suppose that $\Lambda$ is not cuspidal and ${\rm def}(\Lambda)=4k+1$ for some integer $k$.
If $\epsilon=+$, we also assume that $k\neq 0$.
This assumption implies that ${\rm def}(\Lambda')\neq 0$.
Then $\rho_\Lambda$ lies in the parabolic induced character
$\Ind^{\Sp_{2n}(q)}_{\Sp_{2m_0(m_0+1)}(q)\times T}(\rho_{\Lambda_0}\otimes{\bf 1})$
where $T$ is a split torus of rank $n-m_0(m_0+1)$ and $\Lambda_0$ is the cuspidal symbol of rank $m_0(m_0+1)$
for some non-negative integer $m_0$ such that $m_0(m_0+1)<n$.

It is well-known that the theta correspondence is compatible with the parabolic induction (\cf.~\cite{amr} th\'eor\`eme 3.7),
so $\rho_{\Lambda'}$ lies in the parabolic induced character
$\Ind^{\rmO_{2n'}^\epsilon(q)}_{\rmO_{2m'^2_0}^\epsilon(q)\times T'}(\rho_{\Lambda_0'}\otimes{\bf 1})$
where $m_0'$ is equal to $m_0$ or $m_0+1$ depending on $\epsilon$, $ T'$ is a split torus
of rank $2n'-2m_0'^2$, and $\Lambda'_0$ is a cuspidal symbol of rank $m_0'^2$.
Now the defects of $\Lambda_0,\Lambda_0'$ satisfy (\ref{1108}) by (1).
Moreover, we know that ${\rm def}(\Lambda)={\rm def}(\Lambda_0)$ and ${\rm def}(\Lambda')={\rm def}(\Lambda'_0)$ from
the remark before the proposition.
Therefore the defects of $\Lambda,\Lambda'$ satisfy (\ref{1108}),
and hence the corollary is proved.
\end{enumerate}
\end{proof}

\subsection{Branching rule and symbol correspondence}\label{1113}
From now on we consider the case which is not settled by Proposition~\ref{1107}, i.e.,
$\epsilon=+$, ${\rm def}(\Lambda)=1$ and ${\rm def}(\Lambda')=0$.

For a symbol $\Lambda$ of defect $1$ and a set $\Omega'$ of symbols of defect $0$,
we define
\begin{align}
\begin{split}
\Theta_\Lambda(\Omega') &=\{\,\Lambda'\in\Omega'\mid(\Lambda,\Lambda')\in\calb^+\,\},\\
\Theta_\Lambda^*(\Omega') &=\{\,\Lambda'\in\Omega'\mid(\Lambda,\Lambda'^\rmt)\in\calb^+\text{ and }(\Lambda,\Lambda')\not\in\calb^+\,\}.
\end{split}
\end{align}
Similarly, for a symbol $\Lambda'$ of defect $0$ and a set $\Omega$ of symbols of defect $1$,
we define
\begin{align}
\begin{split}
\Theta_{\Lambda'}(\Omega) &=\{\,\Lambda\in\Omega\mid(\Lambda,\Lambda')\in\calb^+\,\},\\
\Theta_{\Lambda'}^*(\Omega) &=\{\,\Lambda\in\Omega\mid(\Lambda,\Lambda'^\rmt)\in\calb^+\text{ and }(\Lambda,\Lambda')\not\in\calb^+\,\}.
\end{split}
\end{align}

\begin{example}
Let $\Lambda=\binom{8,5,1}{6,2}$ and $\Lambda'=\binom{7,4,1}{8,5,1}$.
By Lemma~\ref{0210}, it is easy to check that $(\Lambda,\Lambda')\in\calb^+$ and
\[
\textstyle
\Theta_\Lambda(\Omega^+_{\Lambda'})=
\{\binom{8,4,1}{8,5,1},\binom{7,5,1}{8,5,1},\binom{7,4,2}{8,5,1}\},\qquad
\Theta^*_\Lambda(\Omega^+_{\Lambda'})=
\{\binom{7,4,1}{9,5,1},\binom{7,4,1}{8,6,1},\binom{7,4,1}{8,5,2}\}.
\]
\end{example}

\begin{lemma}\label{1112}
Let $\Lambda,\Lambda'$ be symbols of sizes $(m+1,m),(m',m')$ respectively such that $(\Lambda,\Lambda')\in\calb^+$.
Then
\[
|\Theta_\Lambda(\Omega_{\Lambda'}^+)|=1+|\Theta_{\Lambda'}(\Omega_\Lambda^-)|\quad\text{and}\quad
|\Theta_{\Lambda'}(\Omega_\Lambda^+)|=1+|\Theta_\Lambda(\Omega_{\Lambda'}^-)|.
\]
\end{lemma}
\begin{proof}
Write $\Lambda=\binom{a_1,a_2,\ldots,a_{m+1}}{b_1,b_2,\ldots,b_m}$,
$\Lambda'=\binom{c_1,c_2,\ldots,c_{m'}}{d_1,d_2,\ldots,d_{m'}}$,
and suppose that $(\Lambda,\Lambda')\in\calb^+$.
We know that $m'=m,m+1$ by Lemma~\ref{0213}.
It is clear that the symbol
\[
\Lambda_0'':=\binom{c_1+1,c_2,\ldots,c_{m'}}{d_1,d_2,\ldots,d_{m'}}
\]
is in $\Omega_{\Lambda'}^+$ and $(\Lambda,\Lambda_0'')\in\calb^+$ by Lemma~\ref{0210}.

First suppose that $m'=m+1$.
\begin{enumerate}
\item[(1)] Suppose that
\[
\Lambda_1=\binom{a_1,\ldots,a_{i-1},a_i-1,a_{i+1},\ldots,a_{m+1}}{b_1,b_2,\ldots,b_m}
\]
for some $i$ such that $a_i>a_{i+1}+1$ is in $\Theta_{\Lambda'}(\Omega^-_\Lambda)$, i.e.,
$\Lambda_1\in\Omega^-_\Lambda$ and $(\Lambda_1,\Lambda')\in\calb^+$.
Then we have $d_{i-1}-1>a_i-1\geq d_i$ by Lemma~\ref{0210},
i.e., $d_{i-1}-1>d_i$, and so
\[
\Lambda_1':=\binom{c_1,c_2,\ldots,c_{m'}}{d_1,\ldots,d_{i-1},d_i+1,d_{i+1},\ldots,d_{m'}}
\]
is a symbol in $\Omega^+_{\Lambda'}$.
Moreover, the condition $a_i-1\geq d_i$ implies that $(\Lambda,\Lambda'_1)\in\calb^+$.

\item[(2)] Suppose that
\[
\Lambda_1=\binom{a_1,\ldots,a_{m+1}}{b_1,\ldots, b_{i-1},b_i-1,b_{i+1},\ldots,b_m}
\]
for some $i$ such that $b_i>b_{i+1}+1$ is in $\Theta_{\Lambda'}(\Omega^-_\Lambda)$.
Then we have $c_i-1>b_i-1\geq c_{i+1}$, i.e., $c_i-1>c_{i+1}$,
and so
\[
\Lambda_1':=\binom{c_1,\ldots,c_i,c_{i+1}+1,c_{i+2},\ldots,c_{m'}}{d_1,d_2,\ldots,d_{m'}}
\]
is a symbol in $\Omega^+_{\Lambda'}$.
Moreover, the condition $b_i-1\geq c_{i+1}$ implies that $(\Lambda,\Lambda'_1)\in\calb^+$.
\end{enumerate}
It is easy to see that the mapping $\Lambda_1\mapsto\Lambda_1'$ from $\Theta_{\Lambda'}(\Omega^-_{\Lambda})$
to $\Theta_\Lambda(\Omega^+_{\Lambda'})$ is injective with one extra element $\Lambda_0''$ not in the image.
Hence the lemma is proved for the case that $m'=m+1$.

Next suppose that $m'=m$.
Except for the situations similar to those considered above,
there is another possibility, i.e., the case that $(a_{m+1},b_m)=(1,0),(0,1)$.
Let
\[
\Lambda_1=\binom{a_1-1,\ldots,a_m-1}{b_1-1,\ldots,b_{m-1}-1}.
\]
Then we have $\Lambda_1\in\Omega^-_\Lambda$ and $(\Lambda_1,\Lambda')\in\calb^+$.
Note that now $d_m\geq a_{m+1}$ and $c_m\geq b_m$ since $(\Lambda,\Lambda')\in\calb^+$.
So $d_m\geq 1$ if $(a_{m+1},b_m)=(1,0)$; $c_m\geq 1$ if $(a_{m+1},b_m)=(0,1)$.
Let
\[
\Lambda_1'=\begin{cases}
\displaystyle\binom{c_1+1,\ldots,c_{m'}+1,0}{d_1+1,\ldots,d_{m'}+1,1}, & \text{if $(a_{m+1},b_m)=(1,0)$};\\
\displaystyle\binom{c_1+1,\ldots,c_{m'}+1,1}{d_1+1,\ldots,d_{m'}+1,0}, & \text{if $(a_{m+1},b_m)=(0,1)$}.
\end{cases}
\]
It is easy to check that $\Lambda_1'\in\Omega^+_{\Lambda'}$ and $(\Lambda,\Lambda_1')\in\calb^+$.
Again, we still have an injective mapping from $\Theta_{\Lambda'}(\Omega^-_{\Lambda})$
to $\Theta_\Lambda(\Omega^+_{\Lambda'})$ given by $\Lambda_1\mapsto\Lambda_1'$ with one extra element $\Lambda_0''$ not in the image.
Hence the lemma is proved.
\end{proof}

\begin{lemma}\label{1105}
Let $\Lambda=\binom{a_1,\ldots,a_{m+1}}{b_1,\ldots,b_m}$,
$\Lambda'=\binom{c_1,\ldots,c_{m'}}{d_1,\ldots,d_{m'}}$
be symbols of sizes $(m+1,m),(m',m')$ respectively such that $(\Lambda,\Lambda')\in\calb^+$.
\begin{enumerate}
\item[(i)] If $m'=m+1$, then $\Theta_\Lambda(\Omega^-_{\Lambda'})\neq\emptyset$.

\item[(ii)] If $m'=m$ and $\Theta_{\Lambda'}(\Omega^-_\Lambda)=\emptyset$,
then $m=0$ and $\Lambda=\binom{0}{-}$
\end{enumerate}
\end{lemma}
\begin{proof}
First suppose that $m'=m+1$.
We can define
\[
\Lambda_1'=\begin{cases}
\displaystyle\binom{c_1-1,\ldots,c_m-1}{d_1-1,\ldots,d_m-1}, & \text{if $(c_{m+1},d_{m+1})=(1,0),(0,1)$};\\
\displaystyle\binom{c_1,\ldots,c_m,c_{m+1}-1}{d_1,\ldots,d_{m+1}}, & \text{if $c_{m+1}\geq 1$ and $(c_{m+1},d_{m+1})\neq(1,0)$};\\
\displaystyle\binom{c_1,\ldots,c_{m+1}}{d_1,\ldots,d_m,d_{m+1}-1}, & \text{if $d_{m+1}\geq 1$ and $(c_{m+1},d_{m+1})\neq(0,1)$}.
\end{cases}
\]
It is easy to see that $\Lambda_1'\in\Omega^-_{\Lambda'}$.
Moreover, the assumption that $(\Lambda,\Lambda')\in\calb^+$ implies that
$(\Lambda,\Lambda_1')\in\calb^+$ by Lemma~\ref{0210}.
Thus (i) is proved.

Next suppose that $m'=m$ and $m\geq 1$.
Then we can define
\[
\Lambda_1=\begin{cases}
\displaystyle\binom{a_1-1,\ldots,a_m-1}{b_1-1,\ldots,b_{m-1}-1}, & \text{if $(a_{m+1},b_m)=(1,0),(0,1)$};\\
\displaystyle\binom{a_1,\ldots,a_m,a_{m+1}-1}{b_1,\ldots,b_m}, & \text{if $a_{m+1}\geq 1$ and $(a_{m+1},b_m)\neq (1,0)$};\\
\displaystyle\binom{a_1,\ldots,a_{m+1}}{b_1,\ldots,b_{m-1},b_m-1}, & \text{if $b_m\geq 1$ and $(a_{m+1},b_m)\neq(0,1)$}.
\end{cases}
\]
It is easy to see that $\Lambda_1\in\Omega^-_{\Lambda}$.
Moreover, the assumption that $(\Lambda,\Lambda')\in\calb^+$ implies that
$(\Lambda_1,\Lambda')\in\calb^+$ by Lemma~\ref{0210}.
Therefore we conclude that if $m\geq 1$,
then $\Theta_{\Lambda'}(\Omega^-_\Lambda)\neq\emptyset$.
Next suppose that $m'=m$ and $m=0$, i.e., $\Lambda=\binom{a_1}{-}$ and $\Lambda'=\binom{-}{-}$ for some $a_1\geq 0$.
If $a_1\geq 1$, then $\Lambda_1=\binom{a_1-1}{-}\in\Theta_{\Lambda'}(\Omega^-_{\Lambda})$
and hence $\Theta_{\Lambda'}(\Omega^-_\Lambda)\neq \emptyset$.
\end{proof}

\begin{example}\label{1114}
Let $\Lambda=\binom{a_1}{-}$ and $\Lambda_1'=\binom{c_1}{d_1}$ be symbols of sizes $(1,0)$ and $(1,1)$ respectively
such that $(\Lambda,\Lambda_1')\in\calb^+$.
So now we have $a_1\geq d_1$ by Lemma~\ref{0210}.
Clearly, $\binom{c_1+1}{d_1}\in\Theta_\Lambda(\Omega^+_{\Lambda_1'})$.
If $\Lambda'\in\Omega^+_{\Lambda'_1}$ and $\Lambda'$ is of type (III) or (IV) in Subsection~\ref{1006},
then $\Lambda'$ is of size $(2,2)$ and hence $(\Lambda,\Lambda')\not\in\calb^+$.

Suppose that $\Theta^*_\Lambda(\Omega^+_{\Lambda_1'})\neq\emptyset$.
Then we must have $\Lambda'':=\binom{c_1}{d_1+1}\in\Theta^*_\Lambda(\Omega^+_{\Lambda_1'})$,
i.e., $(\Lambda,\Lambda''^\rmt)\in\calb^+$ and $(\Lambda,\Lambda'')\not\in\calb^+$,
which imply $a_1=d_1$ and $a_1\geq c_1$.
Let
\[
\Lambda_2'=\begin{cases}
\binom{d_1+1}{c_1-1}, & \text{if $c_1\geq 1$};\\
\binom{d_1}{c_1}, & \text{if $c_1=0$}.
\end{cases}
\]
Clearly, $\Lambda''^\rmt\in\Omega^+_{\Lambda_2'}$ and $(\Lambda,\Lambda''^\rmt)\in\calb^+$.
\begin{enumerate}
\item[(1)] If $c_1\geq 1$, then $\Theta_\Lambda(\Omega^+_{\Lambda_2'})=\{\binom{d_1+2}{c_1-1},\binom{d_1+1}{c_1}\}$.

\item[(2)] If $c_1=0$, then $d_1\geq 1$ (since we always consider reduced symbols) and hence $a_1=d_1\geq c_1+1$.
Therefore $\Theta_\Lambda(\Omega^+_{\Lambda_2'})=\{\binom{d_1+1}{c_1},\binom{d_1}{c_1+1}\}$.
\end{enumerate}
For both cases, we conclude that $\Theta^*_\Lambda(\Omega^+_{\Lambda_2'})=\emptyset$.
\end{example}

Now we want to show that the above example is a general phenomena:

\begin{lemma}\label{1101}
Let $\Lambda=\binom{a_1,\ldots,a_{m+1}}{b_1,\ldots,b_m}$,
$\Lambda_1'=\binom{c_1,\ldots,c_{m+1}}{d_1,\ldots,d_{m+1}}$
be symbols of size $(m+1,m),(m+1,m+1)$ respectively such that $(\Lambda,\Lambda_1')\in\calb^+$.
Suppose that $\Lambda''\in\Theta^*_\Lambda(\Omega^+_{\Lambda_1'})$.
Then there exists a symbol $\Lambda'_2$ such that
$\Lambda''^\rmt\in\Omega^+_{\Lambda'_2}$,
$(\Lambda,\Lambda_2')\in\calb^+$ and
$\Theta^*_\Lambda(\Omega^+_{\Lambda_2'})=\emptyset$.
\end{lemma}
\begin{proof}
Because now $(\Lambda,\Lambda_1')\in\calb^+$,
by Lemma~\ref{0210} we have
\begin{align}
\begin{split}
& a_1\geq d_1>a_2\geq d_2>\cdots > a_{m+1}\geq d_{m+1}, \\
& c_1>b_1\geq c_2>b_2\geq \cdots \geq c_m>b_m\geq c_{m+1}.
\end{split}
\end{align}
Let $\Lambda''\in\Theta^*_\Lambda(\Omega^+_{\Lambda_1'})$, i.e,
$\Lambda''\in\Omega^+_{\Lambda_1'}$, $(\Lambda,\Lambda''^\rmt)\in\calb^+$ and $(\Lambda,\Lambda'')\not\in\calb^+$.
If $\Lambda''$, as an element of $\Omega^+_{\Lambda_1'}$, is of type (III) or (IV) in Subsection~\ref{1113},
then $\Lambda'',\Lambda''^\rmt$ are of size $(m+2,m+2)$ and hence $(\Lambda,\Lambda''^\rmt)\not\in\calb^+$ by Lemma~\ref{0213}.
Therefore, $\Lambda''$ must be of type (I) or (II):
\begin{enumerate}
\item[(1)] Suppose that
\[
\Lambda''=\binom{c_1,\ldots,c_{k-1},c_k+1,c_{k+1},\ldots,c_{m+1}}{d_1,\ldots,d_{m+1}}
\]
for some $k$ such that $c_{k-1}>c_k+1$.
If $k=1$, then $(\Lambda,\Lambda_1')\in\calb^+$ implies $(\Lambda,\Lambda'')\in\calb^+$ and we get a contradiction.
So we must have $k\geq 2$.
Because now $(\Lambda,\Lambda'')\not\in\calb^+$ and $(\Lambda,\Lambda''^\rmt)\in\calb^+$,
we have
\begin{itemize}
\item $b_{k-1}=c_k$

\item $a_i\geq c_i>a_{i+1}$ for $i\neq k$ and $a_k\geq c_k+1>a_{k+1}$; $d_j>b_j\geq d_{j+1}$ for each $j$.
\end{itemize}
Then $d_{k-1}>a_k\geq c_k+1=b_{k-1}+1\geq d_k+1$ and hence
\[
\Lambda_2':=\binom{d_1,\ldots,d_{k-2},d_{k-1}-1,d_k,\ldots,d_{m+1}}{c_1,\ldots,c_{k-1},c_k+1,c_{k+1},\ldots,c_{m+1}}
\]
is a symbol.
It is easy to see that $\Lambda''^\rmt\in\Omega^+_{\Lambda_2'}$ and $(\Lambda,\Lambda_2')\in\calb^+$.
Moreover, for any
\[
\Lambda'=\binom{c'_1,\ldots,c'_{m+1}}{d'_1,\ldots,d'_{m+1}}\in\Omega^+_{\Lambda_2'},
\]
we have $d'_k\geq c_k+1$.
Because now $b_{k-1}=c_k<d'_k$, we have $(\Lambda,\Lambda'^\rmt)\not\in\calb^+$.
Therefore $\Theta^*_\Lambda(\Omega^+_{\Lambda_2'})=\emptyset$.

\item[(2)] Suppose that
\[
\Lambda''=\binom{c_1,\ldots,c_{m+1}}{d_1,\ldots,d_{l-1},d_l+1,d_{l+1},\ldots,d_{m+1}}
\]
for some $l\geq 2$ such that $d_{l-1}>d_l+1$.
Because now $(\Lambda,\Lambda'')\not\in\calb^+$ and $(\Lambda,\Lambda''^\rmt)\in\calb^+$,
we have
\begin{itemize}
\item $a_l=d_l$;

\item $a_i\geq c_i>a_{i+1}$ for each $i$; $b_j\geq d_{j+1}>b_{j+1}$ for $j\neq l-1$ and $b_{l-1}\geq d_l+1>b_l$.
\end{itemize}
Then $c_{l-1}>b_{l-1}\geq d_l+1=a_l+1\geq c_l+1$ and hence
\[
\Lambda_2':=\binom{d_1,\ldots,d_{l-1},d_l+1,d_{l+1},\ldots,d_{m+1}}{c_1,\ldots,c_{l-2},c_{l-1}-1,c_l,\ldots,c_{m+1}}
\]
is a symbol.
It is easy to see that $\Lambda''^\rmt\in\Omega^+_{\Lambda_2'}$ and $(\Lambda,\Lambda_2')\in\calb^+$.
Moreover, by the same argument in (1), we can see that $\Theta^*_{\Lambda}(\Omega^+_{\Lambda_2'})=\emptyset$.

\item[(3)] Suppose that
\[
\Lambda''=\binom{c_1,\ldots,c_{m+1}}{d_1+1,d_2,\ldots,d_{m+1}}.
\]
Then we have
\begin{itemize}
\item $a_1=d_1$

\item $a_i\geq c_i>a_{i+1}$ and $b_i\geq d_{i+1}>b_{i+1}$ for each $i$, and $d_1+1>b_1$.
\end{itemize}
If $m=0$, this is just Example~\ref{1114}.
So now we assume that $m\geq 1$.
Note that $c_{m+1},d_{m+1}$ can not both be $0$.
We can define
\[
\Lambda_2'=\begin{cases}
\displaystyle\binom{d_1+1,d_2,\ldots,d_m,d_{m+1}-1}{c_1,\ldots,c_{m+1}},
& \text{if $d_{m+1}\geq 1$ and $(c_{m+1},d_{m+1})\neq (0,1)$};\\
\displaystyle\binom{d_1+1,d_2,\ldots,d_{m+1}}{c_1,\ldots,c_m,c_{m+1}-1},
& \text{if $c_{m+1}\geq 1$ and $(c_{m+1},d_{m+1})\neq (1,0)$};\\
\displaystyle\binom{d_1,d_2-1\ldots,d_m-1}{c_1-1,\ldots,c_m-1},
& \text{if $(c_{m+1},d_{m+1})=(1,0),(0,1)$}.
\end{cases}
\]
For above three cases,
it is easy to see that $\Lambda''^\rmt\in\Omega^+_{\Lambda_2'}$ and $(\Lambda,\Lambda_2')\in\calb^+$.
By the similar argument in (1) we can see that $\Theta^*_\Lambda(\Omega^+_{\Lambda_2'})=\emptyset$.
\end{enumerate}
\end{proof}

\begin{example}
Let $\Lambda=\binom{8,5,1}{6,2}$ and $\Lambda_1'=\binom{7,4,1}{8,3,0}$.
Then $(\Lambda,\Lambda_1')\in\calb^+$,
\[
\textstyle
\Theta_\Lambda(\Omega^+_{\Lambda_1'})=\{\binom{8,4,1}{8,3,0},\binom{7,5,1}{8,3,0},
\binom{7,4,2}{8,3,0},\binom{7,4,1}{8,4,0},\binom{7,4,1}{8,3,1}\}\quad\text{and}\quad
\Theta^*_\Lambda(\Omega^+_{\Lambda_1'})=\{\binom{7,4,1}{9,3,0}\}.
\]
Now $\Lambda''=\binom{7,4,1}{9,3,0}$.
So let $\Lambda_2'=\binom{8,2}{6,3}$ and we have $(\Lambda,\Lambda_2')\in\calb^+$,
\[
\textstyle
\Theta_\Lambda(\Omega^+_{\Lambda_2'})=\{\binom{9,2}{6,3},\binom{8,3}{6,3},
\binom{8,2}{7,3},\binom{8,2}{6,4},\binom{9.3,1}{7,4,0},\binom{9,3,0}{7,4,1}\}\quad\text{and}\quad
\Theta^*_\Lambda(\Omega^+_{\Lambda_2'})=\emptyset.
\]
Note that $\Lambda''^\rmt\in\Theta_\Lambda(\Omega^+_{\Lambda_2'})$.
\end{example}

\begin{lemma}\label{1103}
Let $\Lambda=\binom{a_1,\ldots,a_{m+1}}{b_1,\ldots,b_m}$,
$\Lambda_1'=\binom{c_1,\ldots,c_m}{d_1,\ldots,d_m}$
be symbols of size $(m+1,m),(m,m)$ respectively such that $m\geq 1$ and $(\Lambda,\Lambda_1')\in\calb^+$.
Suppose that $\Lambda''\in\Theta^*_\Lambda(\Omega^+_{\Lambda_1'})$.
Then there exists a symbol $\Lambda'_2$ such that
$\Lambda''^\rmt\in\Omega^+_{\Lambda'_2}$,
$(\Lambda,\Lambda_2')\in\calb^+$ and
$\Theta^*_\Lambda(\Omega^+_{\Lambda_2'})=\emptyset$.
\end{lemma}
\begin{proof}
Because now $(\Lambda,\Lambda_1')\in\calb^+$, by Lemma~\ref{0210} we have
\begin{align}
\begin{split}
& a_1> d_1\geq a_2> d_2\geq \cdots\geq a_m>d_m\geq a_{m+1}, \\
& c_1\geq b_1>c_2\geq b_2>\cdots>c_m\geq b_m.
\end{split}
\end{align}
Let $\Lambda''\in\Theta^*_\Lambda(\Omega^+_{\Lambda_1'})$.
\begin{enumerate}
\item[(1)] Suppose that
\[
\Lambda''=\binom{c_1,\ldots,c_{k-1},c_k+1,c_{k+1},\ldots,c_m}{d_1,\ldots,d_m}
\]
for some $k$ such that $c_{k-1}>c_k+1$.
The proof for this case is similar to (1) in the proof of Lemma~\ref{1101}.

\item[(2)] Suppose that
\[
\Lambda''=\binom{c_1,\ldots,c_m}{d_1,\ldots,d_{l-1},d_l+1,d_{l+1},\ldots,d_m}
\]
for some $l$ such that $d_{l-1}>d_l+1$.
The proof for this case is similar to (2) in the proof of Lemma~\ref{1101}.

\item[(3)] Suppose that
\[
\Lambda''=\binom{c_1+1,\ldots,c_m+1,1}{d_1+1,\ldots,d_m+1,0}.
\]
So we have $c_m\geq 1$.
The assumptions $(\Lambda,\Lambda'')\not\in\calb^+$ and $(\Lambda,\Lambda''^\rmt)\in\calb^+$ imply that
\begin{itemize}
\item $b_m=0$;

\item $a_i>c_i\geq a_{i+1}$ for each $i$;
and $d_j\geq b_j>d_{j+1}$ for $j=1,\ldots,m-1$.
\end{itemize}
Because $\Lambda$ is reduced and now $b_m=0$,
we must have $a_{m+1}\neq 0$.
Hence $d_m\geq a_{m+1}\geq 1$.
Let
\[
\Lambda_2'=\binom{d_1+1,\ldots,d_{m-1}+1,d_m,0}{c_1+1,\ldots,c_m+1,1}.
\]
It is easy to see that $\Lambda''^\rmt\in\Omega^+_{\Lambda_2'}$ and $(\Lambda,\Lambda_2')\in\calb^+$.
Moreover, for any
\[
\Lambda'=\binom{c'_1,\ldots,c'_{m+1}}{d'_1,\ldots,d'_{m+1}}\in\Omega^+_{\Lambda_2'},
\]
we have $d'_{m+1}\geq 1$.
Because now $b_m=0<d'_{m+1}$, we have $(\Lambda,\Lambda'^\rmt)\not\in\calb^+$.
Therefore $\Theta^*_\Lambda(\Omega^+_{\Lambda_2'})=\emptyset$.

\item[(4)] Suppose that
\[
\Lambda''=\binom{c_1+1,\ldots,c_m+1,0}{d_1+1,\ldots,d_m+1,1}.
\]
The proof is similar to (3) above.
\end{enumerate}
\end{proof}

\begin{lemma}\label{1102}
Let $\Lambda_1=\binom{a_1,\ldots,a_{m+1}}{b_1,\ldots,b_m}$,
$\Lambda'=\binom{c_1,\ldots,c_{m'}}{d_1,\ldots,d_{m'}}$
be symbols of size $(m+1,m),(m',m')$ respectively such that $(\Lambda_1,\Lambda')\in\calb^+$.
Suppose that $\Lambda''\in\Theta^*_{\Lambda'}(\Omega^+_{\Lambda_1})$.
Then there exists a symbol $\Lambda_2$ such that
$(\Lambda_2,\Lambda')\in\calb^+$,
$\Lambda''\in\Omega^+_{\Lambda_2}$ and
$\Theta^*_{\Lambda'^\rmt}(\Omega^+_{\Lambda_2})=\emptyset$.
\end{lemma}
\begin{proof}
We know that $m'=m,m+1$.
Then the proof is similar to those of Lemma~\ref{1101} and Lemma~\ref{1103}.
\end{proof}

\subsection{Branching rule and Howe correspondence}\label{1118}
For $\rho\in\cale(\Sp_{2n}(q))$, let $\Omega_\rho^+$ denote the set
of irreducible constituents of the induced character
$\Ind^{\Sp_{2(n+1)}(q)}_{\Sp_{2n}(q)\times\GL_1(q)}(\rho\otimes{\bf 1})$,
then we also define
\[
\Omega_\rho^-=\{\,\rho_1\in\cale(\Sp_{2(n-1)}(q))\mid \rho\in\Omega^+_{\rho_1}\,\}.
\]
The analogous definition also applies to split even orthogonal groups.

For $\rho\in\cale(\Sp_{2n}(q))$ and a set $\Omega'\subset\cale(\rmO^+_{2n'}(q))$,
we define
\begin{equation}
\Theta_\rho(\Omega')=\{\,\rho'\in\Omega'\mid\rho\otimes\rho'\text{ occurs in the Howe correspondence}\,\}.
\end{equation}
Similarly, for $\rho'\in\cale(\rmO^+_{2n'}(q))$ and a set $\Omega\subset\cale(\Sp_{2n}(q))$,
we define
\begin{equation}
\Theta_{\rho'}(\Omega)=\{\,\rho\in\Omega\mid\rho\otimes\rho'\text{ occurs in the Howe correspondence}\,\}.
\end{equation}

The following lemma is extracted from the proof of \cite{amr} th\'eor\`eme~3.7:

\begin{lemma}\label{1109}
Let $\rho\in\cale(\Sp_{2n}(q))$ and $\rho'\in\cale(\rmO^+_{2n'}(q))$
such that $\rho\otimes\rho'$ occurs in the correspondence.
Then
\[
|\Theta_\rho(\Omega_{\rho'}^+)|=1+|\Theta_{\rho'}(\Omega_\rho^-)|\quad\text{and}\quad
|\Theta_{\rho'}(\Omega_\rho^+)|=1+|\Theta_\rho(\Omega_{\rho'}^-)|.
\]
\end{lemma}
\begin{proof}
Suppose that $\rho\in\cale(\Sp_{2n}(q))$, $\rho'\in\cale(\rmO^+_{2n'}(q))$ and
$\rho\otimes\rho'$ occurs in the correspondence.
Note that now the Howe correspondence is of multiplicity one (\cf.~\cite{mvw} p.97).
By Frobenius reciprocity, we have
\begin{align*}
|\Theta_\rho(\Omega_{\rho'}^+)|
& =\Bigl\langle\omega_{\Sp_{2n},\rmO^+_{2(n'+1)}},\rho\otimes\Ind_{\rmO^+_{2n'}(q)\times\GL_1(q)}^{\rmO^+_{2(n'+1)}(q)}(\rho'\otimes1)\Bigr\rangle \\
& =\langle\omega_{\Sp_{2n},\rmO^+_{2(n'+1)}}|_{\Sp_{2n}(q)\times\rmO^+_{2n'}(q)\times\GL_1(q)},\rho\otimes\rho'\otimes 1\rangle.
\end{align*}
From \cite{amr} p.382, we know that
\begin{multline*}
\omega_{\Sp_{2n},\rmO^+_{2(n'+1)}}|_{\Sp_{2n}(q)\times\rmO^+_{2n'}(q)\times\GL_1(q)} \\
=\omega_{\Sp_{2n},\rmO^+_{2n'}}\otimes{\bf 1}
+\Ind^{\Sp_{2n}(q)\times\rmO^+_{2n'}(q)\times\GL_1(q)}_{\Sp_{2(n-1)}(q)\times\rmO^+_{2n'}(q)\times\GL_1(q)\times\GL_1(q)}
(\omega_{\Sp_{2(n-1)},\rmO^+_{2n'}}\otimes R_{\GL_1})
\end{multline*}
where $R_{\GL_1}$ denotes the character of the representation of $\GL_1(q)\times\GL_1(q)$ on $\bbC(\GL_1(q))$.
By our assumption, we have
\[
\langle\omega_{\Sp_{2n},\rmO^+_{2n'}}\otimes{\bf 1},\rho\otimes\rho'\otimes{\bf 1}\rangle=1.
\]
Moreover, by Frobenius reciprocity again,
\begin{align*}
& \Bigl\langle\Ind^{\Sp_{2n}(q)\times\rmO^+_{2n'}(q)\times\GL_1(q)}_{\GL_1(q)\times\GL_1(q)\times\Sp_{2(n-1)}(q)\times\rmO^+_{2n'}(q)}
(\omega_{\Sp_{2(n-1)},\rmO^+_{2n'}}\otimes R_{\GL_1}),\rho\otimes\rho'\otimes{\bf 1}\Bigr\rangle \\
& =\langle\omega_{\Sp_{2(n-1)},\rmO^+_{2n'}},(\rho|_{\Sp_{2(n-1)}(q)})\otimes\rho'\rangle \\
& =|\Theta_{\rho'}(\Omega_\rho^-)|.
\end{align*}
Hence the first equality is proved.
The proof of the second equality is similar.
\end{proof}

\begin{lemma}\label{1111}
Suppose that $\rho_\Lambda\otimes\rho_{\Lambda_1'}$ occurs in the correspondence,
$(\Lambda,\Lambda_1')\in\calb^+$, and $\Theta^*_\Lambda(\Omega^+_{\Lambda'_1})=\emptyset$.
For any $\Lambda'\in\Omega_{\Lambda_1'}^+$,
if $\rho_\Lambda\otimes\rho_{\Lambda'}$ occurs in the correspondence,
then $(\Lambda,\Lambda')\in\calb^+$.
\end{lemma}
\begin{proof}
Suppose that $\Lambda'\in\Omega_{\Lambda_1'}^+$ and $\rho_\Lambda\otimes\rho_{\Lambda'}$ occurs
in the correspondence.
Then we know that $(\Lambda,\Lambda')$ or $(\Lambda,\Lambda'^\rmt)$ is in  $\calb^+$ by Theorem~\ref{0522}.
But $\Theta^*_\Lambda(\Omega^+_{\Lambda_1'})=\emptyset$ means that there is no $\Lambda''\in\Omega^+_{\Lambda_1'}$
such that $(\Lambda,\Lambda''^\rmt)\in\calb^+$ and $(\Lambda,\Lambda'')\not\in\calb^+$.
Hence we have $(\Lambda,\Lambda')\in\calb^+$.
\end{proof}

\begin{lemma}
Suppose that $\rho_{\Lambda_1}\otimes\rho_{\Lambda'}$ occurs in the correspondence,
$(\Lambda_1,\Lambda')\in\calb^+$, and $\Theta^*_{\Lambda'}(\Omega^+_{\Lambda_1})=\emptyset$.
For any $\Lambda\in\Omega_{\Lambda_1}^+$, if $\rho_\Lambda\otimes\rho_{\Lambda'}$
occurs in the correspondence,
then $(\Lambda,\Lambda')\in\calb^+$.
\end{lemma}
\begin{proof}
The proof is similar to that of Lemma~\ref{1111}.
\end{proof}

\begin{proposition}\label{1116}
Let $(\bfG,\bfG')=(\Sp_{2n},\rmO_{2n'}^+)$,
$\rho_\Lambda\in\cale(G)_1$, $\rho_{\Lambda'}\in\cale(G')_1$
for some symbols $\Lambda,\Lambda'$ of defects $1,0$ respectively.
Then $\rho_\Lambda\otimes\rho_{\Lambda'}$ occurs in $\omega_{\bfG,\bfG',1}$ if and only if
$(\Lambda,\Lambda')\in\calb_{\bfG,\bfG'}$.
\end{proposition}
\begin{proof}
For the case that ${\rm def}(\Lambda')>0$,
by Proposition~\ref{1107},
we know that $\rho_\Lambda\otimes\rho_{\Lambda'}$
occurs in the correspondence if and only if $(\Lambda,\Lambda')\in\calb_{\bfG,\bfG'}$.
From Theorem~\ref{0310} we know that
\[
\omega_{\bfG,\bfG',1}^\sharp=\sum_{(\Lambda,\Lambda')\in\calb_{\bfG,\bfG'}}
\rho_\Lambda^\sharp\otimes\rho_{\Lambda'}^\sharp.
\]
Therefore, to prove the proposition, it suffices to show that if ${\rm def}(\Lambda')=0$ and $(\Lambda,\Lambda')\in\calb_{\bfG,\bfG'}$,
then $\rho_\Lambda\otimes\rho_{\Lambda'}$ occurs in the correspondence.

We now want to prove the above assertion by induction on $n+n'$.
First consider the case that $n'=0$.
The Howe correspondence for the dual pair $(\Sp_{2n}(q),\rmO_0^+(q))$ is given by
${\bf 1}_{\Sp_{2n}}\otimes{\bf 1}_{\rmO^+_0}$ and we know that
${\bf 1}_{\Sp_{2n}}=\rho_{\binom{n}{-}}$, ${\bf 1}_{\rmO_0^+}=\rho_{\binom{-}{-}}$.
It is clear that
\[
\calb_{\Sp_{2n},\rmO^+_0}=\{(\textstyle\binom{n}{-},\binom{-}{-})\}.
\]
Hence the proposition holds for $n'=0$ (and any non-negative integer $n$).

Suppose that $(\Lambda,\Lambda')\in\calb_{\Sp_{2n},\rmO^+_{2n'}}$ and ${\rm def}(\Lambda')=0$ for some $n,n'$
such that $n'>0$ and write
$\Lambda=\binom{a_1,\ldots,a_{m+1}}{b_1,\ldots,b_m}$,
$\Lambda'=\binom{c_1,\ldots,c_{m'}}{d_1,\ldots,d_{m'}}$
for some nonnegative integers $m,m'$.
It is known that $m'=m,m+1$ by Lemma~\ref{0213}.

\begin{enumerate}
\item[(1)] Suppose that $m'=m+1$.
By Lemma~\ref{1105}, we know that $\Theta_\Lambda(\Omega^-_{\Lambda'})\neq\emptyset$.
Let $\Lambda_1'\in\Theta_\Lambda(\Omega^-_{\Lambda'})$, i.e.,
$\Lambda'\in\Omega^+_{\Lambda'_1}$ and $(\Lambda,\Lambda_1')\in\calb^+$.
Now $\Lambda_1'$ is a symbol of rank $n'-1$,
then, by induction hypothesis,
$\rho_\Lambda\otimes\rho_{\Lambda_1'}$ occurs in the correspondence and we have
\begin{equation}\label{1119}
|\Theta_{\rho_{\Lambda_1'}}(\Omega^-_{\rho_\Lambda})|=|\Theta_{\Lambda_1'}(\Omega^-_\Lambda)|.
\end{equation}
Now by Lemma~\ref{1109} and Lemma~\ref{1112},
we have the equality
\begin{equation}\label{1117}
|\Theta_{\rho_\Lambda}(\Omega^+_{\rho_{\Lambda_1'}})|=|\Theta_{\Lambda}(\Omega^+_{\Lambda_1'})|.
\end{equation}
\begin{enumerate}
\item Suppose that $\Lambda_1'$ is of type (I') or (II') in Subsection~\ref{1006}.
Then $\Lambda_1'$ is of size $(m+1,m+1)$.
\begin{enumerate}
\item[(i)]
First suppose that $\Theta^*_\Lambda(\Omega^+_{\Lambda_1'})=\emptyset$.
For any $\Lambda''\in\Omega^+_{\Lambda_1'}$,
by Lemma~\ref{1111},
if $\rho_\Lambda\otimes\rho_{\Lambda''}$ occurs in the correspondence,
then $(\Lambda,\Lambda'')\in\calb^+$.

\item[(ii)]
Next suppose that $\Theta^*_\Lambda(\Omega^+_{\Lambda_1'})\neq\emptyset$.
For any $\Lambda'''\in\Theta^*_\Lambda(\Omega^+_{\Lambda_1'})$,
we know that $(\Lambda,\Lambda'''^\rmt)\in\calb^+$ and
$(\Lambda,\Lambda''')\not\in\calb^+$ by the definition.
By Lemma~\ref{1101}, we can find $\Lambda_2'$ of rank $n'-1$
such that $\Lambda'''^\rmt\in\Omega^+_{\Lambda_2'}$, $(\Lambda,\Lambda_2')\in\calb^+$
and $\Theta^*_\Lambda(\Omega^+_{\Lambda_2'})=\emptyset$.
By (i), with $\Lambda'_1$ replaced by $\Lambda'_2$,
we conclude that $\rho_\Lambda\otimes\rho_{\Lambda'''^\rmt}$ occurs in the correspondence,
and hence $\rho_\Lambda\otimes\rho_{\Lambda'''}$ does not occur in the correspondence
by Corollary~\ref{1024}.
So now for any $\Lambda''\in\Omega^+_{\Lambda'_1}$,
if $\rho_\Lambda\otimes\rho_{\Lambda''}$ occurs in the correspondence
and $(\Lambda,\Lambda'')\not\in\calb^+$,
we must have $(\Lambda,\Lambda''^\rmt)\in\calb^+$, i.e.,
$\Lambda''\in\Theta^*_\Lambda(\Omega^+_{\Lambda_1'})$ and we get a contradiction.
Therefore, for any $\Lambda''\in\Omega^+_{\Lambda'_1}$,
the occurrence of $\rho_\Lambda\otimes\rho_{\Lambda''}$ in the correspondence
implies that $(\Lambda,\Lambda'')\in\calb^+$.
\end{enumerate}

Hence for both (i) and (ii), by (\ref{1117}),
the condition $(\Lambda,\Lambda'')\in\calb^+$ also implies that
$\rho_\Lambda\otimes\rho_{\Lambda''}$ occurs in the correspondence.
In particular,
since $\Lambda'\in\Omega^+_{\Lambda_1'}$ and $(\Lambda,\Lambda')\in\calb^+$,
we have that $\rho_\Lambda\otimes\rho_{\Lambda'}$ occurs in the correspondence.

\item Suppose that $\Lambda_1'$ is of type (III') in Subsection~\ref{1006}.
Then $\Lambda_1'$ is of size $(m,m)$.
\begin{enumerate}
\item[(i)] First suppose that $\Theta^*_\Lambda(\Omega^+_{\Lambda_1'})=\emptyset$.
The proof is exactly the same as in (a.i).

\item[(ii)] Suppose that $\Theta^*_\Lambda(\Omega^+_{\Lambda_1'})\neq\emptyset$.
First suppose that $m=0$.
This means that $\Lambda=\binom{n}{-}$ for some $n\geq 0$, $\Lambda_1'=\binom{-}{-}$, $\Lambda'=\binom{1}{0}$
and $\Theta^*_\Lambda(\Omega^+_{\Lambda_1'})=\{\binom{0}{1}\}$.
It is well known that $\rho_{\binom{n}{-}}={\bf 1}_{\Sp_{2n}}$, $\rho_{\binom{1}{0}}={\bf 1}_{\rmO^+_2}$,
${\bf 1}_{\Sp_{2n}}\otimes{\bf 1}_{\rmO^+_2}$ occurs in the correspondence,
and $(\binom{n}{-},\binom{1}{0})\in\calb^+$, i.e., the proposition is true for this case.
Next suppose that $m\geq 1$.
Then the proof is similar to that of (a.ii).
The only difference is that we need to apply Lemma~\ref{1103} instead of Lemma~\ref{1101}.
\end{enumerate}
\end{enumerate}

\item[(2)] Suppose that $m'=m$.
Since the case that $m=0$ is just the case that $n'=0$,
we assume that $m\geq 1$.
By Lemma~\ref{1105}, we know that $\Theta_{\Lambda'}(\Omega^-_\Lambda)\neq\emptyset$.
Let $\Lambda_1\in\Theta_{\Lambda'}(\Omega^-_\Lambda)$, i.e.,
$\Lambda\in\Omega^+_{\Lambda_1}$ and $(\Lambda_1,\Lambda')\in\calb^+$.
Then $\rho_{\Lambda_1}\otimes\rho_{\Lambda'}$ occurs in the correspondence by induction hypothesis.
The remaining proof is similar to that of Case (1).
Note that we need to apply Lemma~\ref{1102} instead of Lemma~\ref{1101} and Lemma~\ref{1103}.
\end{enumerate}
\end{proof}

\begin{proof}[Proof of Theorem~\ref{0334}]
The theorem is just the combination of Proposition~\ref{1107} and Proposition~\ref{1116}.
\end{proof}

\bibliography{refer}

\providecommand{\bysame}{\leavevmode\hbox to3em{\hrulefill}\thinspace}
\providecommand{\MR}{\relax\ifhmode\unskip\space\fi MR }
\providecommand{\MRhref}[2]{%
  \href{http://www.ams.org/mathscinet-getitem?mr=#1}{#2}
}
\providecommand{\href}[2]{#2}
\begin{thebibliography}{MVW87}

\bibitem[AM93]{adams-moy}
J.~Adams and A.~Moy, \emph{Unipotent representations and reductive dual pairs
  over finite fields}, Trans. Amer. Math. Soc. \textbf{340} (1993), 309--321.

\bibitem[AMR96]{amr}
A.-M. Aubert, J.~Michel, and R.~Rouquier, \emph{Correspondance de {H}owe pour
  les groupes r\'eductifs sur les corps finis}, Duke Math. J. \textbf{83}
  (1996), 353--397.

\bibitem[Car85]{carter-finite}
R.~Carter, \emph{Finite groups of {L}ie type, conjugacy classes and complex
  characters}, John Wiley \& Sons, England, 1985.

\bibitem[G{\'e}r77]{gerardin}
P.~G{\'e}rardin, \emph{Weil representations associated to finite fields}, J.
  Algebra \textbf{46} (1977), 54--101.

\bibitem[GP00]{Geck-Pfeiffer}
M.~Geck and G.~Pfeiffer, \emph{Characters of finite {C}oxeter groups and
  {I}wahori-{H}ecke algebras}, Oxford University Press, New York, 2000.

\bibitem[KS05]{Kable-Sanat}
A.~Kable and N.~Sanat, \emph{The exterior and symmetric square of the
  reflection representation of ${A}_n(q)$ and ${D}_n(q)$}, J. Alegbra
  \textbf{288} (2005), 409--444.

\bibitem[Lus77]{lg}
G.~Lusztig, \emph{Irreducible representations of finite classical groups},
  Invent. Math. \textbf{43} (1977), 125--175.

\bibitem[Lus81]{lg-symplectic}
\bysame, \emph{Unipotent characters of the symplectic and odd orthogonal groups
  over a finite field}, Invent. Math. \textbf{64} (1981), 263--296.

\bibitem[Lus82]{lg-orthogonal}
\bysame, \emph{Unipotent characters of the even orthogonal groups over a finite
  field}, Trans. Amer. Math. Soc. \textbf{272} (1982), 733--751.

\bibitem[MVW87]{mvw}
C.~M{\oe}glin, M.-F. Vign\'eras, and J.-L. Waldspurger, \emph{Correspondances
  de {H}owe sur un corps p-adiques}, Lecture Notes in Math., vol. 1291,
  Springer-Verlag, Berlin-Heidelberg-New York, 1987.

\bibitem[Pan19a]{pan-Lusztig-correspondence}
S.-Y. Pan, \emph{Lusztig correspondence and {H}owe correspondence for finite
  reductive dual pairs}, arXiv:1906.01158 (2019).

\bibitem[Pan19b]{pan-chain01}
\bysame, \emph{Supercuspidal representations and preservation principle of
  theta correspondence}, J. Reine Angew. Math. \textbf{750} (2019), 1--52.

\bibitem[Pan20a]{pan-uniform}
\bysame, \emph{Decomposition of the uniform projection of the {W}eil
  character}.

\bibitem[Pan20b]{pan-eta}
\bysame, \emph{On theta and eta correspondences for finite
  symplectic/orthogonal dual pairs}, arXiv:2006.06241 (2020).

\end{thebibliography}
\bibliographystyle{amsalpha}

\end{document}